\documentclass[a4paper,english,10pt]{article}

\usepackage{etex}
\usepackage{amsmath}
\usepackage{amssymb}
\usepackage{amsthm}
\usepackage[english]{babel}
\usepackage{graphics}
\usepackage[all,cmtip]{xy}
\usepackage{xcolor}
\usepackage{hyperref}
\usepackage{epsfig}

\newcommand{\Adj}{\operatorname{Adj}}

\newcommand{\Card}{\operatorname{Card}}
\newcommand{\Cl}{\operatorname{Cl}}

\newcommand{\Cone}{\operatorname{Cone}}
\newcommand{\Conv}{\operatorname{Conv}}
\newcommand{\Div}{\operatorname{Div}}
\newcommand{\Eff}{\operatorname{Eff}}
\newcommand{\Hom}{\operatorname{Hom}}
\newcommand{\Mob}{\operatorname{Mob}}

\newcommand{\mult}{\operatorname{mult}}
\newcommand{\NE}{\overline{\operatorname{NE}}}
\newcommand{\Nef}{\operatorname{Nef}}

\newcommand{\rk}{\operatorname{rank}}
\newcommand{\Span}{\operatorname{Span}}
\newcommand{\Spec}{\operatorname{Spec}}
\newcommand{\Star}{\operatorname{Star}}

\renewcommand{\div}{\operatorname{div}}

\newcommand{\V}{\textup{\textbf{V}}}

\renewcommand{\AA}{\mathbb A}

\newcommand{\PP}{\mathbb P}
\newcommand{\NN}{\mathbb N}
\newcommand{\ZZ}{\mathbb Z}
\newcommand{\QQ}{\mathbb Q}
\newcommand{\RR}{\mathbb R}
\newcommand{\CC}{\mathbb C}
\newcommand{\GG}{\mathbb G}
\newcommand{\FF}{\mathbb F}

\newcommand{\AAA}{\mathcal A}
\newcommand{\BBB}{\mathcal B}

\newcommand{\KKK}{\mathcal K}

\newcommand{\RRR}{\mathcal R}

\renewcommand{\bar}[1]{\overline{#1}}
\newcommand{\barphi}[1]{\overline{#1}^{\phi}}
\newcommand{\vv}[2]{\left(\begin{smallmatrix}#1 \\ #2 \end{smallmatrix}\right)} 
\newcommand{\sslash}{\mathbin{/\! /}}

\renewcommand{\geq}{\geqslant}
\renewcommand{\leq}{\leqslant}

\newtheorem*{theorem*}{Theorem}
\newtheorem{theo}{Theorem}[section]
\newtheorem{defi}[theo]{Definition}
\newtheorem{cor}[theo]{Corollary}
\newtheorem{lem}[theo]{Lemma}
\newtheorem{lemdef}[theo]{Lemma-Definition}
\newtheorem{prop}[theo]{Proposition}
\newtheorem*{conj*}{Conjecture}
\newtheorem{nota}[theo]{Notation}
\newtheorem*{question*}{Question}
\newtheorem{exemple}[theo]{Example}
\newtheorem{remark}[theo]{Remark}
\newtheorem{remarks}[theo]{Remarks}

\newenvironment{preuve}{\medbreak \noindent {\bf Proof~---}}
                       {\hfill $\square$ \medbreak}

\newenvironment{demo}[1]{\medbreak \noindent {\bf #1~---}}
                        {\hfill $\square$ \medbreak}

\def\cqf{\mathrel{%
    \mathchoice{\CQF}{\CQF}{\scriptsize\CQF}{\tiny\CQF}%
}}
\def\CQF{{%
    \setbox0\hbox{$\square$}%
    \rlap{\hbox to \wd0{\raisebox{0.04cm}{$\sim$}}}\box0 
}}

\newenvironment{ex}{\begin{exemple} \normalfont}
                        {\end{exemple}\hfill \Large{$\diamond$} \medbreak}
\newenvironment{Remark}{\begin{remark} \normalfont}
                   {\end{remark}}
\newenvironment{Remarks}{\begin{remarks} \normalfont}
                   {\end{remarks}}

\renewcommand*{\eqref}[1]{%
  \begingroup
    \hypersetup{
      linkcolor=blue,
      linkbordercolor=blue,
    }%
    \ref{#1}%
  \endgroup
}

\usepackage[disable]{todonotes} % notes not showed
\begin{document}

\thispagestyle{empty}

\begin{center}
\LARGE 
Low degree hypersurfaces of projective toric 
varieties defined over $C_1$ fields have a rational point

\vspace{0.2cm}
\normalsize{(Long version)}

\vspace{0.3cm}

\Large{\textsc{Robin Guilbot}}
\end{center}

%\begin{abstract}
%Quasi algebraically closed fields, or $C_1$ fields, are defined in terms of a low degree condition. Namely, the field $K$ is $C_1$ if every degree $d$ hypersurface of the projective space $\mathbb{P}_K^n$ contains a $K$-point as soon as $d\leq n$. 

%In this article we define a notion of low toric degree generalizing this condition for hypersurfaces of simplicial projective split toric varieties.
%This allows us to prove a particular case of the $C_1$ conjecture of Koll\'ar, Lang and Manin~: any smooth separably rationnally connected variety that can be embedded as such a hypersurface over a $C_1$ field $K$ has a rational point.

%Our results are based on the fact that the ambient toric varieties are Mori Dream Spaces~: they are naturally endowed with homogeneous coordinates and their Minimal Model Program works in all cases.
%\end{abstract}

\setcounter{tocdepth}{2}
\tableofcontents

\section{Introduction}
 
The definition of quasi algebraically closed fields, or $C_1$ fields, relies on a notion of low degree for hypersurfaces of a projective space~:

\begin{defi}\label{defC1}
The field $K$ is a \textit{$C_1$ field} if every degree $d$ hypersurface of the projective space $\PP_K^n$ contains a $K$-rational point as soon as $d\leq n$. 
\end{defi}

J.-L. Colliot-Th\'{e}l\`{e}ne gives in \cite{Colliot} a list of classical results of existence of rational points in various classes of algebraic varieties over given fields. The list ends on the following general question~:
\begin{question*}\label{Colliot}
Sur un type donn\'{e} de corps, y a-t-il une classe naturelle de vari\'{e}t\'{e}s
alg\'{e}briques qui sur un tel corps ont automatiquement un point rationnel ?
\footnote{For a given type of field, does there exist a natural class of algebraic varieties having automatically a rational point over such a field ?}

\end{question*}

For $C_1$ fields, Fano varieties seem to be good candidates for such a natural class, according to conjectures by S. Lang and by Y.I. Manin and J. Koll\'{a}r (\cite{Lang00} and \cite[Rem.~2.6]{Manin}). 
But during the years 2000, strong results suggested to enlarge a bit the scope of these expectations and give a new formulation of the so-called ``$C_1$ conjecture''~:

\begin{conj*}[J. Koll\'{a}r, S. Lang, Y.I. Manin]\label{Folk}
Every proper smooth separably rationally connected variety over a $C_1$ field has a rational point.
\end{conj*} 
\smallskip

A variety $X$ over a field $K$ is rationally connected if every pair of points of $X$ can be joined by a rational curve over any algebraically closed extension of $K$ (see section \ref{RC} for more definitions). The notion of rational connectedness was introduced independently by F.~Campana and by J.~Koll\'ar, Y.~Miyaoka and S.~Mori in 1992. It arises naturally from Mori's program of classification of projective algebraic varieties as a good higher-dimensional analogue of rational curves and rational surfaces.  

The $C_1$ conjecture has already been proved for almost every known type of $C_1$ field :
\begin{itemize} 
\item For a finite field $\FF_q$, by H. Esnault in \cite{Esnault}. Her result is in fact much stronger since it needs only rational chain-connectedness (see Subsection \ref{Ltd=>RCC} for the definition).
\item For a field of functions $\bar{k}(C)$ of a curve $C$ over an algebraically closed field $\bar{k}$, by T. Graber, J. Harris, A.J. de Jong and J. Starr in \cite{GHS03} and \cite{deJS03}. Their result is more geometric since it is stated in terms of sections of a fibration on the curve $C$.
\item For a field of fractions of formal series with coefficients in an algebraically closed field $\bar{k}((t))$. This is a consequence of the previous case, see \cite[Th. 7.5]{Colliot}.  
\end{itemize}

As there are not many other known $C_1$ fields (besides the pseudo-algebraically closed ones for which the conjecture is trivial, see \cite[§1.1.3]{Witt}), the $C_1$ conjecture remains open essentially for 
the maximal unramified extension of the $p$-adic numbers $\QQ_{p}^{nr}$.
\bigskip

The aim of this article is to present a new approach to the $C_1$ conjecture, somewhat orthogonal to the one of the articles cited above : instead of considering \textit{every} separably rationally connected variety over \textit{a given} $C_1$ field, as is done in those papers, we prove the conjecture \textit{only} for hypersurfaces of a simplicial projective split toric variety endowed with homogeneous coordinates (see \ref{TorC1} and \ref{Homcoor} for precisions) but over an arbitrary $C_1$ field admitting normic forms of any degree, hence over \textit{every known} $C_1$ field. 
\smallskip

We proceed in two steps. First we define geometrically the notion of low toric degree for a hypersurface $D$ of a simplicial projective split toric variety $X_\Sigma$ (Theorem \ref{gral} and Definition \ref{ltd}) and show that it implies the existence of a rational point in the hypersurface $D$ when $X_\Sigma$ and $D$ are defined over a $C_1$ field whose characteristic is prime to the order of the torsion of the class group $\Cl(X_\Sigma)$ (Corollary \ref{ratpt}). 
This takes the major part of the paper. 
Then we show that any smooth and separably rationally connected hypersurface of $X_\Sigma$ has low toric degree, yielding the main result of the paper~:

\begin{theorem*}[\ref{RCpoint}]\label{RCpoint0}
Let $X_{\Sigma}$ be a simplicial projective toric variety defined over a $C_1$ field $K$ such that $K$ admits normic forms of arbitrary degree, the dense torus of $X_\Sigma$ is split over $K$ and the order of the torsion part of $\Cl(X_\Sigma)$ is prime to the characteristic of $K$. Let $D$ be a hypersurface of $X_{\Sigma}$ defined over $K$.  
If $D$ is smooth and separably rationally connected then it has low toric degree and in particular admits a rational point over $K$. 
\end{theorem*}
\smallskip

Along the article we define in fact three different notions of low toric degree. Roughly speaking the first is local, the second is global and the third is a combination of the two others, involving several birational modifications of the ambient variety $X_\Sigma$ and the toric morphism defined by the linear system of $D$.
\medskip

Let us give some details. 

\subsubsection*{Restricted low toric degree.}

The first step in defining a low degree condition similar to the one of Definition \ref{defC1} for more general ambient varieties was made by J. Koll\'ar in \cite{Kollar}. It draws a notion of low degree for an intersection of hypersurfaces in a weighted projective space by generalizing \cite{Lang}, Th.~4~: 

\begin{theo}[Th. 6.7 p.232 \cite{Kollar}]\label{Kollar}
Let $K$ be a $C_1$ field admitting normic forms of arbitrary degree. Let $K[x_0, \ldots, x_n]$ be the polynomial algebra graded by $\deg(x_i)=a_i \in \NN^*$. 
Let $f_1, \ldots, f_r \in K[x_0, \ldots, x_n]$ be homogeneous polynomials for this graduation. 
If
$$
\sum_{j=1}^{r} \deg(f_j) < \sum_{i=0}^{n} a_i,
$$
then the system of equations $f_1 = \cdots = f_r=0$ has a non trivial solution in $K^{n+1}$.
\end{theo}
\smallskip

This result applies also to hypersurfaces in a fake weighted projective space, that is, the quotient of a well formed one by a finite abelian group, provided this ambient space admits homogeneous coordinates (i.e. provided the order of the finite group is prime to the characteristic of $K$ in order for the quotient to be a good geometric one, see \ref{Homcoor} ).
It thus defines a notion of low weighted degree for (intersections of) hypersurfaces of any complete simplicial toric variety of Picard number 1 with homogeneous coordinates.   
\smallskip

We say that an effective divisor $D$ of a (normal separated) toric variety $X_\Sigma$ 
has restricted low toric degree if there exists a complete simplicial toric subvariety $W \subseteq  X_\Sigma$ of Picard number 1 such that the restriction $D_{|W}$ as low weighted degree in the sense of Theorem \ref{Kollar}.
\smallskip

The notion of restricted low toric degree is a local notion since it involves nothing but a specific toric subvariety of the ambient $X_\Sigma$. In particular it does not require $X_\Sigma$ either to be  simplicial or to be projective. 
\smallskip

One crucial observation is that complete simplicial toric subvarieties $W \subseteq  X_\Sigma$ 
of Picard number 1 appear naturally as fibers of extremal contractions, i.e., contractions of an extremal ray of the Mori cone $\NE(X_\Sigma)$ (see \ref{TorMMP}). 
In fact every extremal ray of $\NE(X_\Sigma)$ is generated by the class of a line belonging to a family that covers such a fiber $W$. 
In particular (see Lemma \ref{extrgral}) the divisor $D$ has restricted low toric degree as soon as there is an extremal class $[C]$ such that 
\begin{equation}\label{titdegextr}
0 < C \cdot D < C \cdot L_W,
\end{equation}
 where $L_W$ is any divisor whose restriction to  $W$ is the anticanonical divisor $-K_W$.   

An important property of the restricted low toric degree condition for $D$ (Lemma \ref{dltdbar}) is that it holds as soon as the same condition holds for one of its pullback $\phi^*D$ by a contracting toric morphism $\phi$ (see Definition \ref{tormap}).

\subsubsection*{Global low toric degree.}

Let $X_\Sigma$ be a complete simplicial toric variety with homogeneous coordinates (i.e. such that the order of the torsion part of the class group is prime to the characteristic of $K$, see \ref{Homcoor} ).
By Cox's construction $X_\Sigma$ is isomorphic to the good geometric quotient \begin{equation*}  \left( \AA^{\Sigma(1)} \setminus Z(\Sigma) \right) \sslash G
\end{equation*}
where $\Sigma(1)$ is the set of 1-dimensional cones of the fan $\Sigma$, $G=\Hom(\Cl(X_{\Sigma}),\GG_m)$ is the group whose character group is the class group $\Cl(X_{\Sigma})$ and $Z(\Sigma)$ is the exceptional set (see \ref{Homcoor}). 

In particular every hypersurface $D$ of $X_\Sigma$ admits a homogeneous equation~: there exists a polynomial $f$ in the $[D]$-homogeneous piece of the Cox ring
$K[x_{\rho}, \rho \in \Sigma(1)]$, (which is graded by $\Cl(X_{\Sigma})$) such that
\begin{equation*}
D=\V(f):=\lbrace \pi(x) \in X_{\Sigma} \mid f(x)=0 \rbrace, 
\end{equation*} 
where $\pi : \AA^{\Sigma(1)} \setminus Z(\Sigma) \rightarrow X_{\Sigma}$ is the  quotient by $G$.
\medskip

With this setting, we can find $K$-rational points in $D$ as roots $\alpha \in K^{\Sigma(1)}$ of $f$ over $K$ that are not contained in the exceptional set $Z(\Sigma)$, so that $\pi(\alpha) \in D$. 
\bigskip

On the one hand, since $f$ is homogeneous for the graduation by the class group $\Cl(X_\Sigma)$, it is also homogeneous for any graduation by $\ZZ$ induced by the intersection numbers with a given 1-cycle $C$, i.e. defined by putting $\deg(x_\rho)= m C \cdot D_\rho$ for a well chosen $m\in\NN^*$.

By Theorem \ref{Kollar}, we can hence deduce the existence of a root $\alpha \in K^{\Sigma(1)}$ of $f$ from some low degree condition along a 1-cycle $C$ of $X_\Sigma$, that is, from inequalities of the form 
\begin{equation}\label{titdeg0}
0< C \cdot D <  \sum_{\rho \in I} C \cdot D_{\rho}, 
\end{equation}
where $C \cdot D_{\rho}>0$ for all $\rho \in I$ (Lemma \ref{fund}).
\medskip

On the other hand, to get a rational point we have to ensure that such a root is relevant, i.e., does not belong to $Z(\Sigma)$. This is a priori a quite restrictive condition since it requires all the zero coordinates of $\alpha$ to be contained in a single cone of $\Sigma$~: 
\begin{equation*}
\exists \sigma \in \Sigma, \ 
\lbrace \rho \in \Sigma(1) \mid \alpha_\rho = 0 \rbrace \subseteq \sigma(1). 
\end{equation*}

We then distinguish three cases~:
\begin{enumerate}
\item[$1^{\text{st}}$ case] The effective divisor $D$ is not Cartier and nef. Then it contains automatically a torus-invariant point, that is a $\ZZ$-rational point (Proposition \ref{trivial}). We may thus always assume that $D$ is a nef Cartier divisor.   
\item[$2^{\text{nd}}$ case] The divisor $D$ is ample. If there is an inequality of the form \eqref{titdeg0} with $C$ effective, then writing the class $[C]$ as a linear combination of extremal classes with positive coefficient, we can deduce an inequality of the form \eqref{titdegextr} (Proposition \ref{tresbeau}), which implies that $D$ has restricted low toric degree. 
\item[$3^{\text{rd}}$ case] The divisor $D$ is Cartier and nef but not ample. 
In this case it is quite natural to consider the proper toric morphism $\phi : X_\Sigma \to X_{\bar{\Sigma}_{[D]}}$ associated to the complete linear system of $D$. It contracts all the curves that have zero intersection with $D$ and presents the class $[D]$ as the pullback of an ample class $[\bar{D}]$ in $X_{\bar{\Sigma}_{[D]}}$. 
The fan $\bar{\Sigma}_{[D]}$ is obtained from $\Sigma$ first by removing all the walls $\tau \in \Sigma(n-1)$ such that $C_\tau \cdot D=0$. This yields a fan $\Sigma_{[D]}$ that is degenerate if $D$ is not big. In this case we quotient by the linear space $U_0$ contained in every cone of $\Sigma_{[D]}$ to get a genuine fan $\bar{\Sigma}_{[D]}$ 
of lower dimension.
 
The key property is that the fan $\Sigma_{[D]}$ is the normal fan of a polytope that can be identified with the Newton polytope of the homogeneous polynomial $f$ defining $D$ (see \ref{Polytopes}). 
The consequence (Proposition \ref{easierpoint}) is that a root $\alpha \in K^{\Sigma(1)}$ of $f$ yields a rational point of $D$ as soon as it does not belong to the exceptional set associated to $\Sigma_{[D]}$, that is as soon as  
\begin{equation}\label{relevant'}
\exists \sigma' \in \Sigma_{[D]}, \ 
\lbrace \rho \in \Sigma(1) \mid \alpha_\rho = 0 \rbrace \subseteq \sigma'. 
\end{equation}

We say that $D$ has global low toric degree if there exist a 1-cycle $C$ on $X_\Sigma$ and a subset $I \subseteq \Sigma(1)$ such that inequalities \eqref{titdeg0} are verified and any root of $f$ given by Lemma \ref{fund} satisfies condition \eqref{relevant'}.   
\end{enumerate}

An important consequence of the appearance of $\Sigma_{[D]}$ in the definition is that the global low toric degree condition can be pulled back by toric birational contractions, i.e. equivariant birational maps surjective in codimension 1 (Proposition \ref{open2}).  

\subsubsection*{Low toric degree.}

The general notion of low toric degree comes as a combination of the two other notions. 
It is defined by three conditions whose equivalence constitutes the main result of the paper.  

\begin{theorem*}[\ref{gral}]
Let $X_{\Sigma}$ be a simplicial projective toric variety.
Let $D$ be a nef Cartier divisor on $X_{\Sigma}$.
Let $\phi : X_{\Sigma} \rightarrow X_{\bar{\Sigma}_{[D]}}$ be the toric morphism of Theorem \ref{phi}, contracting all the rational curves having zero intersection with $D$
and let $\bar{D}$ be an ample divisor on $X_{\bar{\Sigma}_{[D]}}$ such that $D$ is linearly equivalent to the pullback $\phi^*\bar{D}$.

The following conditions are equivalent~:
\begin{enumerate}
\item[(i)] The ample divisor $\bar{D}$ has restricted low toric degree. 

\item[(ii)] There exist a subvariety $V \subset X_{\Sigma}$  
such that the restriction $D_{|V}$ has global low toric degree.
 
\item[(iii)] There exists a desingularization $\eta : X_{\tilde{\Sigma}} \rightarrow X_{\Sigma}$ such that the pullback $\tilde{D}=\eta^* D$ has global low toric degree.
\end{enumerate}
\end{theorem*}
\medskip

The proof relies entirely on tools from the toric Minimal Model Program.
A sketch of the proof is given short after the statement in Section \ref{skpr}.
\bigskip

Here is the outline of the paper~:
\smallskip

We begin by reviewing in Section 2 known results of toric geometry, with some explanations on how they apply over an arbitrary field.  

In Section 3 we show that a divisor $D$ has automatically rational points over any field if it is not Cartier and nef. Then we define the notion of restricted low toric degree and show that it implies the existence of rational points over $C_1$ fields. We also point out its relation to low intersection numbers with extremal curves.

We define the notion of global low toric degree in Section 4, giving some of its important properties, notably the fact that it implies the existence of rational points over $C_1$ fields and its behavior under toric birational contractions.   

Section 5 contains the statement and the proof of the core result of the article yielding the general definition of low toric degree, together with the proof that it ensures the existence of rational points over $C_1$ fields.

In Section 6 we state some general and easy to handle conditions implying or equivalent to the low toric degree condition. 
Then we discuss briefly the interpretation of low degree conditions in terms of proximity to the boarder of classical cones of divisor classes such as the nef cone or the effective cone. 

Finally in Section 7 we show that hypersurfaces that are smooth and separably rationally connected have low toric degree, proving that the $C_1$ conjecture is true for hypersurfaces of simplicial projective split toric varieties with homogeneous coordinates.

\subsubsection*{Acknowledgements}
The author wishes to thank M. Perret for providing this beautiful subject and for his help all along the Ph.D, S. Lamy for helpful conversations about the geometric meaning of the earlier results and C. Araujo for her suggestions and help during the final period of this work.
This work began during the author's Ph.D Thesis in France, and was finished with the support of the CNPq, Conselho Nacional de Desenvolvimento Cient\'{i}fico e Tecnol\'{o}gico, Brasil.

\section{Preliminaries}

\subsection{General setting : Toric varieties over a $C_1$ field} \label{TorC1}

In this article, we work over an arbitrary $C_1$ field $K$. The only additional hypothesis we make on $K$ is satisfied by every known $C_1$ field~:

\begin{equation*} \tag{H1}
K \ \text{admits normic forms of arbitrary degree.} \
\end{equation*}

A normic form of degree $d$ is a homogeneous polynomial in $d$ variables without non trivial root. This hypothesis is used in Theorem \ref{Kollar} which is the fundamental building block of all the notions of low toric degree. It is not known whether the definition of $C_1$ field implies this hypothesis or not. 
Basic facts on $C_1$ fields, and more generally $C_i$ fields, can be found in \cite{Lang} and \cite[Ch.1]{Witt}. 
\bigskip

Many aspects of toric geometry can be studied regardless of the choice a particular ground field. This is what is done for instance by V. Danilov in the classical \cite{Danilov}. 
However most references among which the (also very classical) \cite{Fulton} and \cite{CLS} work exclusively over the field of complex numbers $\CC$.
In fact, a normal separated toric variety $X$ defined over a field $K$ can be handled almost in the same way as in the complex case as soon as the following hypothesis is satisfied~:
\begin{equation*} \tag{H2}
\ \text{The dense torus} \ T \subset X \ \text{acting on} \ X \ \text{is split over} \ K. \ 
\end{equation*}
This means that we have an isomorphism $T \simeq_K \GG_m^n$ over $K$, where $n$ is the dimension of $X$.
 Indeed in this case the usual construction of $X$ as a gluing of toric affine varieties work as well as if $K$ were algebraically closed (see \cite[Ch. 3]{CLS}). 
More specifically this hypothesis is a necessary and sufficient condition for the character group $M$ of $T$ and its group of 1-parameter subgroups $N$ to be isomorphic to $\ZZ^n$. 
Such a toric variety is said to be split over $K$.
For a discussion of what happens when the torus is not split, see \cite{Huruguen} and \cite{ELST14}.
\bigskip

Throughout the article, by a toric variety we mean a normal separated split toric variety $X_\Sigma$ defined by a fan $\Sigma$ contained in the real span $N_\RR$ of the lattice $N$.
We fix an isomorphism $T \simeq_K \GG_m^n$, inducing another $T_{\bar{K}} \simeq (\bar{K}^*)^n$ over the algebraic closure of $K$. 
This in turn induces an isomorphism $M \simeq \ZZ^n$ given by $m\in\ZZ^n \leftrightarrow (\chi^m : (t_1, \ldots, t_n) \mapsto t_1^{m_1}  \cdots t_n^{m_n})$.
By abuse we identify $M$ and $N$ to $\ZZ^n$ and write $m=(m_1, \ldots, m_n)\in M$, $u=(u_1, \ldots, u_n)\in N$.

The fan $\Sigma$ is contained in the space $N_\RR = N \otimes \RR \simeq\RR^n$.
For $k=1,\ldots,n$, we denote by $\Sigma(k)$ the set of $k$-dimensional cones of $\Sigma$. 
We denote by $\beta \preceq \gamma$ the fact that the cone $\beta \in \Sigma$ is a face of $\gamma \in \Sigma$. 
To each cone $\gamma\in\Sigma(k)$ corresponds a toric subvariety $\V(\gamma)$ that is a subvariety of dimension $n-k$ invariant under the action of $T$. 
In particular $D_\rho$ denotes the toric divisor associated to $\rho \in \Sigma(1)$ and $C_\tau$ denotes the toric curve associated to $\tau \in \Sigma(n-1)$. 
The minimal generator of the 1-dimensional cone $\rho \in \Sigma(1)$, is denoted by $u_\rho$. 
We say that the fan $\Sigma$ (or the toric variety $X_\Sigma$) is simplicial, resp. complete, if each cone $\gamma\in\Sigma(k)$ is generated by exactly $k$ 1-dimensional cones, resp. if the cones of $\Sigma$ cover the whole space $N_\RR$. 

The class group of $X_\Sigma$ is denoted by $\Cl(X_\Sigma)$ and the space of numerical classes of divisors by $N^1(X_\Sigma)=\Cl(X_\Sigma) \otimes \RR$.   
Its dual, the space of numerical classes of curves, is denoted by $N_1(X_\Sigma)$.
Linear equivalence of divisors as well as numerical equivalence of 1-cycles are denoted with brackets~: $[D]$ and $[C]$.

\subsection{Homogeneous coordinates} \label{Homcoor}

Giving an algebro-geometric meaning of earlier analytic results of Audin, Delzant and Kerwan, D. Cox published in 1995 the famous article \cite{CoxJAG} showing that any (normal separated) complex toric variety is naturally isomorphic to a GIT quotient, this quotient being geometric if and only if the variety is simplicial. 

It turns out that this result holds in full generality over an arbitrary field of characteristic zero and also in characteristic $p$ for toric varieties satisfying the following additional hypothesis~:

\begin{equation*} \tag{H3}
\ \text{The order of the torsion part of the class group} \ \Cl(X) \ \text{is prime to} \ p.  
\end{equation*}

Let us give more details.

D. Cox showed that in the complex case, a simplicial toric variety $X_\Sigma$ is isomorphic to the good geometric quotient  

\begin{equation*} 
 \left( \AA^{\Sigma(1)} \setminus Z(\Sigma) \right) \sslash G
\end{equation*}
where 
$$
Z(\Sigma)=\V(x^{\hat{\sigma}} \mid \sigma \in \Sigma)
$$ 
is the exceptional set, that is the subvariety of $\AA^{\Sigma(1)}=\Spec(K[x_{\rho}, \rho \in \Sigma(1)])$ defined by the vanishing of every monomial of the form $x^{\hat{\sigma}}= \prod_{\rho \notin \sigma(1)} x_\rho$ for $\sigma \in \Sigma$, and 
$$
G=\Hom(\Cl(X_{\Sigma}),\GG_m)
$$ 
is the group acting on $\AA^{\Sigma(1)} \setminus Z(\Sigma)$ whose character group is the class group $\Cl(X_{\Sigma})$ (see \cite[Ch. 5]{CLS}). 
In fact, the original proof of the existence of such a quotient 
(Th. 2.1 of \cite{CoxJAG}) relies on Theorem 1.1 of \cite{GIT} and its amplification 1.3. These results only depend on the reductivity of the group $G$ by which we quotient. 
If we denote by $F$ the torsion part of $\Cl(X_{\Sigma})$, we have 
$
\Cl(X_{\Sigma}) \simeq  \ZZ^{s-n} \times F , 
$ 
where $s$ is the cardinality of $\Sigma(1)$ and $n$ the dimension of $X_{\Sigma}$, so that
$$
G \simeq (\GG_m)^{s-n}  \times \Hom(F,\GG_m). 
$$ 
The factor on the left is of course a reductive torus. The one on the right is a finite group scheme and hence is reductive as soon as (H3) is verified. 
\bigskip

Under the hypotheses (H2) and (H3), subvarieties of $X_{\Sigma}$ are thus given by homogeneous ideals of the total coordinate ring (or Cox ring) of $X_{\Sigma}$, which is the polynomial algebra generated over $K$ by one variable for each 1-dimensional cone $\rho \in \Sigma(1)$, and graded by the class group~: 
$$
S = K[x_{\rho}, \rho \in \Sigma(1)] = \bigoplus_{\beta \in \Cl(X_{\Sigma})} S_\beta.
$$
In particular, any hypersurface $D$ is the zero locus in $X_\Sigma$ of a polynomial $f$ belonging to the homogeneous piece $S_{[D]}$~:
\begin{equation*}
D=\V(f):=\lbrace \pi(x) \in X_{\Sigma} \mid f(x)=0 \rbrace, 
\end{equation*} 
where $\pi : \AA^{\Sigma(1)} \setminus Z(\Sigma) \rightarrow X_{\Sigma}$ is the good geometric quotient.
\medskip

In this setting, our results give criteria of existence of rational points of the form $\pi(\alpha)$ for $\alpha \in K^{\Sigma(1)}$ a root of $f$ which is not in the exceptional set $Z(\Sigma)$, that is to say such that the set $\lbrace \rho \in \Sigma(1) \mid \alpha_{\rho}=0 \rbrace$ generates a cone of the simplicial fan $\Sigma$.

\subsection{Toric intersection theory}\label{TorInt}

Intersection theory on a normal toric variety $X_\Sigma$ is made quite easy by the fact that the Chow homology group $A_k(X_\Sigma)$ of classes of $k$-cycles on $X_\Sigma$ are generated by the classes $[\V(\gamma)]$ of $k$-dimensional toric subvarieties ($\gamma \in \Sigma(n-k)$).
Most of the computations thus reduce to combinatorics in the fan and linear algebra.

\subsubsection{Chow groups.}\label{Chow}

The class group is computed through the following exact sequence~:
$$
0 \longrightarrow M \stackrel{\div}{\longrightarrow} \ZZ^{\Sigma(1)}  
\stackrel{\deg}{\longrightarrow} A^1(X_{\Sigma})=\Cl(X_\Sigma) \longrightarrow 0,
$$
where $\div(m)=\left( \langle m, u_{\rho} \rangle \right)_{\rho \in \Sigma(1)}$ and $\deg(D)=[D]$.

Tensoring by $\RR$, this gives
\begin{equation}\label{exact}
0 \longrightarrow M_{\RR}  \stackrel{\div}{\longrightarrow} \RR^{\Sigma(1)}  \stackrel{\deg}{\longrightarrow} N^1(X_{\Sigma}) \longrightarrow 0,
\end{equation}

and taking the dual sequence we get
$$
0 \longrightarrow N_1(X_{\Sigma}) \stackrel{\deg^*}{\longrightarrow} \RR^{\Sigma(1)} 
 \stackrel{\div^*}{\longrightarrow} N_\RR \longrightarrow 0,
$$
where $\deg^*([C]) = (C\cdot D_{\rho})_{\rho\in \Sigma(1)}$  and 
$\div^* \left( (b_{\rho})_{\rho\in \Sigma(1)} \right) = \sum_{\rho\in \Sigma(1)} b_\rho u_\rho$.

This last exact sequence identifies the space of numerical classes of 1-cycles $N_1(X_{\Sigma})$ with the space of linear relations between minimal generators of the 1-dimensional cones of $\Sigma$. This allows to consider the numerical class $[C]$ as given by the relation 
\begin{equation}\label{relC}
\div^* \left( \deg^*([C])\right) = 0, \quad \text{i.e.} \quad
\sum_{\rho\in \Sigma(1)} C\cdot D_{\rho} u_\rho = 0.
\end{equation}

\subsubsection{Intersection numbers.}

The exact sequences of the previous paragraph show that the numerical classes of divisors (resp. curves) are determined by their intersection numbers with the toric curves (resp. divisors).

In particular for curves we introduce the following notations :

\begin{nota}
$J_C^\pm = \lbrace \rho \in \Sigma(1) \mid C \cdot D_{\rho} \gtrless 0 \rbrace$ and $J_C=J_C^+ \cup J_C^-$.
\end{nota} 

The intersection matrix $(C_\tau \cdot D_\rho)_{\tau \in \Sigma(n-1),\rho \in \Sigma(1)}$ \, is easily computed using the relation \eqref{relC} and the fact that for a couple $(\tau ,\rho) \in \Sigma(n-1) \times \Sigma(1)$ such that $\rho \nprec \tau$ we have 
\begin{equation}\label{CtauDrho}
C_\tau \cdot D_\rho = \begin{cases} 
\dfrac{\mult(\tau) 
}{\mult(\sigma)} & \text{if } \tau + \rho = \sigma \in \Sigma(n) \\
0 & \text{if } \tau + \rho \notin \Sigma,
\end{cases}
\end{equation}
where the multiplicity of a cone $\gamma=\Cone(u_1, \ldots, u_k) \in \Sigma(k)$ is the index of the lattice generated by its minimal generators $u_1, \ldots, u_k$ in the maximal sublattice $N_\gamma=N \cap \Span(\gamma)$ of $N$ contained in the vector space spanned by the cone $\gamma$~:
$$
\mult(\gamma) = \bigl[ N_\gamma : \bigoplus_{i=1}^{k}\ZZ u_i \bigr].
$$  

\subsubsection{Pushforwards and pullbacks.}

Pushforwards and pullbacks of cycle classes by toric morphisms can be described explicitly in terms of fans and linear algebra. However the general formulas are rather complicated (see \cite{Park} and \cite{Fulsturm}).
\medskip

We focus on cycles of dimension or codimension 1 and describe their behaviour under some rational maps. We first fix some terminology.

\begin{defi}\label{tormap}
A torus-equivariant rational map between toric varieties is called a \textit{toric rational map}. To such a map $\varphi : X_{\Sigma} \dashrightarrow X_{\Sigma'}$ corresponds a unique linear map $\bar{\varphi} : N_\RR \rightarrow N_\RR'$ between the latticed spaces containing the fans $\Sigma$ and $\Sigma'$, and $\bar{\varphi}$ is compatible with the fans (i.e. for all $\sigma \in \Sigma$ there exists $\sigma' \in \Sigma'$ such that $\bar{\varphi}(\sigma)\subset \sigma'$) if and only if $\varphi$ is a \textit{toric morphism}.

We say that $\varphi$ is a \textit{toric rational contraction} if the restriction $\bar{\varphi} : N \to N'$ is surjective and $\varphi$ is surjective in codimension 1, i.e. we have $\Sigma'(1) \subseteq \bar{\varphi}(\Sigma(1))$. A regular toric rational contraction is called a \textit{contracting toric morphism}. 
\end{defi}

Let us describe some properties of pushforwards and pullbacks by a toric birational contraction $\xi : X_\Sigma \dashrightarrow X_{\Sigma'}$ between simplicial varieties.  
Since $\xi$ is birational we may without loss of generality assume that the fans $\Sigma$ and $\Sigma'$ belong to the same space $N_\RR$ and that the linear map $\bar{\xi}$ is the identity of $N_\RR$.
Since moreover $\xi$ is surjective in codimension 1 we have $\Sigma'(1) \subseteq \Sigma(1)$, and for $\rho \in \Sigma(1)$, we have 
$\xi_*D_\rho = \begin{cases} D'_{\rho} & \text{if} \ \rho \in \Sigma'(1), \\
0 & \text{otherwise}, \end{cases}$ 
where $D'_{\rho}$ is the toric divisor of $X_\Sigma'$ associated to $\rho \in \Sigma'(1)$.
This induces a  surjective map between the groups of torus-invariant divisors as well as a surjective pushforward map $\xi_* : N^1(X_\Sigma) \rightarrow N^1(X_\Sigma')$ for numerical classes. 

In general, the pullback of divisor classes are easy to describe in terms of their associated polytopes or the support functions of some torus-invariant representative, which we do not address here, the interested reader may refer to \cite[§6.1 and 6.2]{CLS}.
Let us just mention that in our case the pullback of divisor classes \mbox{$\xi^* : N^1(X') \to N^1(X)$} is a section of the pushforward $\xi_* : N^1(X) \to N^1(X')$. In particular it is injective and the composition $\xi_* \circ \xi^*$ is the identity on $N^1(X')$. 
\smallskip

By duality, we can deduce a pullback map for numerical classes of curves. It is called numerical pullback of curves, was introduced by Batyrev in \cite{Bat92} in a more general context, and have been neatly described by C. Araujo in \cite{Ara10}. 

\begin{defi}\label{numpullback}
Let $\xi : X_\Sigma \dashrightarrow X_\Sigma'$ be a toric birational contraction. 
The \textit{numerical pullback of curves} 
$\xi^*_{num} : N_1(X') \to N_1(X)$ is the dual linear map of 
$\xi_* : N^1(X) \to N^1(X')$. It is the unique injective
linear map verifying the following projection formula~:
$$
\xi^*_{num}[C'] \cdot [D] = [C'] \cdot \xi_*[D] \quad \ \text{for all} \
\; [C'] \in  N_1(X') \ \text{and} \ [D] \in N^1(X).
$$
\end{defi}

In the toric context, the numerical pullback $\xi^*_{num}[C']$ is characterized by the equalities 
$$
\xi^*_{num}[C']\cdot [D_\rho] = \begin{cases} 
D'_{\rho} & \text{if} \ \bar{\xi}(\rho)=\rho' \in \Sigma'(1), \\
0 & \text{otherwise}.
\end{cases}
$$ 
In particular when $[C]=\xi^*_{num}[C']$ we have $J_C^\pm = J_{C'}^\pm$.   

Finally the pushforward of curve classes may either be defined directly on invariant curves or by duality via the intersection pairing. In particular it can as well be described in terms of relations between minimal generators, and if $\xi$ is birational, then it admits $\xi^*_{num}$ as a linear section. 

Notice that the usual projection formula always holds~:  
$$
\xi_*[C] \cdot [D'] = [C] \cdot \xi^*[D'] \quad \ \text{for all} \
[C] \in  N_1(X) \ \text{and} \ [D'] \in N^1(X').
$$ 

\subsubsection{Intersections in a toric subvariety.}\label{InterSubvar}

The notions of low toric degree we define in the following sections are based on inequalities between intersection numbers. 
Since they are designed to ensure the existence of a rational point in a divisor $D$, one would expect they obey the principle~: 
\begin{center}
For any subvariety $V  \subset X_\Sigma$,  $D_{|V}$ has low degree $\Rightarrow$ $D$ has low degree.
\end{center}

It turns out that it is not that simple because for a proper toric subvariety $\V(\gamma) \subset X_\Sigma$, \mbox{$\gamma \in \Sigma(k)$}, $1 \leq k \leq n-1$, the intersection numbers computed in $\V(\gamma)$ and the ones computed in $X_\Sigma$ do not always coincide. Let us explicit the relation between them.

Let us denote by $\iota: \V(\gamma) \hookrightarrow X_\Sigma$ the inclusion. The associated linear map between latticed vector spaces is the surjection $\bar{\iota} : N_\RR \rightarrow N(\gamma)_\RR=N(\gamma)\otimes \RR$ where $N(\gamma)=N/N_\gamma$. 
Then $\V(\gamma)$ is the toric variety defined by the fan
$$
\Star(\gamma)=\left\lbrace \bar{\iota}(\sigma) \mid \gamma \preceq \sigma \in \Sigma \right\rbrace
\; \subseteq N(\gamma)_\RR.
$$

The 1-dimensional cones of $\Star(\gamma)$ are the images through $\bar{\iota}$ of the 1-dimensional cones of $\Sigma$ adjacent to $\gamma$~:
$$
\Star(\gamma)(1)=\left\lbrace \bar{\iota}(\rho) \mid \rho \in \Adj(\gamma) \right\rbrace
\ \text{where} \ 
\Adj(\gamma)=\left\lbrace \rho \in \Sigma(1) \setminus \gamma(1), \gamma + \rho \in \Sigma \right\rbrace
$$
The tricky point is that, when the variety $X_\Sigma$ is singular, the minimal generator $u_{\bar{\rho}}$ of $\bar{\iota}(\rho)$ in the lattice $N(\gamma)$ is not necessarily the image $\bar{\iota}(u_\rho)$ of the minimal generator of $\rho$.  
What we have in general is $\bar{\iota}(u_\rho) = m_\rho u_{\bar{\rho}}$ \, with \, 
$m_\rho = \dfrac{\mult(\gamma + \rho)}{\mult(\gamma)} \in \NN^*$.

Hence if $D_{\bar{\rho}}$ denotes the toric divisor in $\V(\gamma)$ corresponding to the 1-dimensional cone $\bar{\rho} \in \Star(\gamma)(1)$, then the pullbacks of toric divisors are
\begin{equation*}
\iota^* D_\rho = \begin{cases} 
\frac{1}{m_\rho} D_{\bar{\rho}} & \text{if} \ \rho \in \Adj(\gamma), \\
0 & \text{else}.
\end{cases}
\end{equation*}

It follows that for every 1-cycle $C$ on $\V(\gamma)$ and every $\rho \in \Adj(\gamma)$ we have
\begin{equation*}
C \cdot D_{\bar{\rho}} = C \cdot m_\rho \, \iota^* D_{\rho} = m_\rho \, \iota_*C \cdot D_{\rho}
\geq \iota_*C \cdot D_{\rho}.  
\end{equation*}
and the inequality may be strict when $X_\Sigma$ is singular.

For later use, notice that we have in particular
\begin{equation}\label{locglob3}
C \cdot (-K_{\V(\gamma)}) \geq \iota_*C \cdot \sum_{\rho \in \Adj(\gamma)} D_{\rho}.  
\end{equation}

\subsection{Projective toric varieties, divisors and polytopes}\label{Polytopes}

The fan $\Sigma$ defining a projective toric variety is always the normal fan of a lattice polytope $P$. This means that the cones of $\Sigma$ are generated by the inner normal facets of 
$P$. In this case we write $\Sigma=\Sigma_P$. Such fans are said to be strongly polytopal.

In this subsection we review some results from \cite[Ch. 2, 4, 5 and 6]{CLS} linking lattice polytopes with projective toric varieties and their divisors.

\subsubsection{Polytopes and projective toric varieties}\label{PolPro}

At a very basic level, projective toric varieties are easy to construct out of lattice polytopes since they appear as Zariski closures of the image of a map defined by such a polytope.
Let us recall some of the important results of \cite[Ch. 2]{CLS}.
\smallskip

Let $P \subset M_\RR$ be a full dimensional convex polytope with vertices in the character group $M$ of the $n$-dimensional torus $T$, and let $P \cap M = \lbrace m_1, \ldots, m_s \rbrace$ be the set of all lattice points contained in $P$.  Suppose, for beginning, that $P$ is normal, which means that for all $k\geq 1$ we have
$$
kP \cap M = \underset{k \ \text{times}}{\underbrace{P \cap M + \cdots + P \cap M}}
$$
We have a map
$$
\begin{array}{rccc}
\Phi : & T            &    \rightarrow      & \PP^{s-1} \\
       & (t_1, \ldots, t_n)  & \mapsto     & \left( \chi^{m_1}, \ldots, \chi^{m_s} \right).  
\end{array}
$$
And the Zariski closure $Y_{P \cap M}=\bar{\Phi(T)}$ of the image of $\Phi$ is a projective toric variety having $\Phi(T) \simeq T$ as dense open subset. 
\bigskip

Now consider the normal fan $\Sigma_P$ of the polytope $P$. It is the fan whose maximal cones are the dual cones $C_v ^\vee$ where 
$$
C_v=\Cone(P \cap M -v), \quad  v \ \text{vertex of} \ P.
$$
The toric variety $X_{\Sigma_P}$ defined from this fan is (isomorphic to) $Y_{P \cap M}$.
\smallskip

In fact, even in the case where the full dimensional polytope $P$ is not a normal polytope, its multiples $kP$ are normal for all $k\geq n-1$. 
The toric variety associated to $P$ in this case is therefore 
$$
X_P := X_{\Sigma_P} \simeq Y_{(n-1)P \cap M}.
$$ 
Note that since, by Sumihiro's Theorem, any (normal separated) toric variety comes from a fan (\cite[Cor. 3.1.8]{CLS}), it follows that any such variety comes from a full dimensional lattice polytope as soon as it is projective.
But still, several polytopes may have the same normal fan and hence correspond to the same projective toric variety.

\subsubsection{Polytopes and torus-invariant divisors}
\bigskip

Let $D=\sum_{\rho \in \Sigma(1)} a_\rho D_\rho$ be a torus-invariant divisor on a complete simplicial toric variety $X_\Sigma$. For each 1-dimensional cone $\rho \in \Sigma(1)$, the inequality $\langle m_\sigma(D), u_\rho \rangle \geq -a_\rho$ defines an affine half space of $M_\RR\simeq\RR^n$. The intersection of all these half spaces is the convex polytope 
$$
P_D := \left\lbrace m \in M \mid \langle m, u_\rho \rangle \geq -a_\rho \quad  
\text{for all} \ \rho \in \Sigma(1) \right\rbrace
$$

This polytope characterizes the set of torus-invariant divisors $D'=\sum_{\rho \in \Sigma(1)} b_\rho D_\rho$ that are  effective and linearly equivalent to $D$. Indeed for such a divisor we have
$$
D' \ \text{is effective  and } D' \in [D] \quad \Leftrightarrow \quad
\exists m \in P_{D}, D'=D+\div(\chi^m),
$$
where $[D]$ denotes the linear equivalence class of $D$ and $\div(\chi^m)=\sum_{\rho \in \Sigma(1)}\langle m, u_\rho \rangle D_\rho$ is the divisor of the character $\chi^m, m \in M$.  

Since $\Sigma$ is simplicial, for each $n$-dimensional cone $\sigma \in \Sigma(n)$ the system of linear equations
$$
\langle m, u_\rho \rangle = -a_\rho \ , \ \rho \in \sigma(1)
$$
has a unique solution denoted $m_\sigma(D)$.

The set $\lbrace m_\sigma(D) \mid \sigma \in \Sigma(n) \rbrace$ is called the Cartier data of the torus invariant divisor $D$ because $D$ is a Cartier divisor if and only if all the $m_\sigma(D)$ are lattice points in $M$. 

Indeed for each $\sigma \in \Sigma(n)$ the fact that $m_\sigma(D)$ is a character of $T$ is equivalent to the existence of a torus-invariant divisor
$$
D^\sigma= D+\div(\chi^{m_\sigma}) = \sum_{\rho \in \Sigma(1)} a^\sigma_\rho D_\rho, \quad
\text{with} \quad a^\sigma_\rho = a_\rho+\langle m, u_\rho \rangle = 0 \quad 
\text{for all} \ \rho \in \sigma(1), 
$$ 
which means that the linear equivalence class is trivial on the toric affine open $U_\sigma$, and this happens for all $\sigma \in \Sigma(n)$ if and only if $D$ is locally principal. 
Note that the divisors $D^\sigma$ are uniquely determined by the class $[D]$, we call them local representatives of $[D]$.
\bigskip

Let us suppose that $D=\sum_{\rho \in \Sigma(1)} a_\rho D_\rho$ is a Cartier divisor. 
By definition, each character $m_\sigma(D)$ lie at the intersection of $n=\dim(M_\RR)$ supporting hyperplanes of the polytope $P_D$. 
Hence $m_\sigma(D)$ belongs to $P_D$ if and only if it is one of its vertices, and this happens for all $\sigma \in \Sigma(n)$ if and only if all the local representatives $D^\sigma$ are effective. 
Moreover the basepoints of $D$ are necessarily torus-invariant points, that is points of the form $\gamma_\sigma = \bigcap_{\rho \in \sigma(1)} D_{\rho}$  for some $\sigma \in \Sigma(n)$. 
Such a point is by definition avoided by the support of the local representative $D^\sigma$. 
It follows that we have 
\begin{equation}\label{Dbpf}
P_D=\Conv( m_\sigma(D), \sigma \in \Sigma(n))  \ \Leftrightarrow \ 
\forall \sigma \in \Sigma(n), D^\sigma \geq 0  \ \Leftrightarrow \ 
D \text{ is basepoint free}.
\end{equation} 

And when $D$ is base point free we have (see \cite[6.1]{CLS} for a proof) 
\begin{equation}\label{Dample}
m_\sigma(D)\neq m_{\sigma'}(D) \ \text{ for }  \sigma \neq \sigma' 
\ \Leftrightarrow \ 
D^\sigma\neq D^{\sigma'} \ \text{ for }  \sigma \neq \sigma' 
\ \Leftrightarrow \ 
D \text{ is ample}.
\end{equation} 

In this case, for $\sigma \in \Sigma(n)$, the character $m_\sigma$ is the only one verifying  for all $\rho \in \Sigma(1)$  
$$
\langle m, u_\rho \rangle \geq -a_\rho, \quad  
\text{with equality if and only if} \ \rho \in \sigma(1). 
$$
It follows on the one hand that  the coefficients of $D$ may be recovered from the polytope $P_D$ by the formula
$$
a_\rho = - \min_{m \in P_D} \langle m, u_\rho \rangle = - \min_{m \in P_D \cap M} \langle m, u_\rho \rangle \, , \quad 
\text{for all } \rho \in \Sigma(1).
$$

On the other hand, the normal fan $\Sigma_{P_D}$ is in fact the starting fan $\Sigma$. 
Indeed, its maximal cones $C_{m_\sigma} ^\vee$ verify for all $\rho \in \Sigma(1)$  
$$
\begin{array}{rcl}
u_\rho \in C_{m_\sigma} ^\vee & \Leftrightarrow & 
\forall m \in \Cone(P \cap M -m_\sigma), \langle m, u_\rho \rangle \geq 0  \\
& \Leftrightarrow & 
\forall p \in P\cap M, \langle \, p, u_\rho \rangle \geq \langle m_\sigma, u_\rho \rangle  \\
& \Leftrightarrow & \langle m_\sigma, u_\rho \rangle = -a_\rho
\ \Leftrightarrow \ u_\rho \in \sigma(1). 
\end{array}
$$

Finally it turns out that full dimensional lattice 
polytopes parametrize the set of pairs formed by a projective toric variety and a torus invariant ample divisor on it~: 

\begin{theo}[\cite{CLS} Th. 6.2.1]\label{invmaps}
The sets
$$
\left\lbrace (X_\Sigma , D) \mid \Sigma \text{ is a complete fan in } N_\RR, 
D \text{ is a torus-invariant ample divisor on } X_\Sigma \right\rbrace
$$
and 
$$
\left\lbrace P \subset M_\RR  \mid  P \text{ is a full dimensional lattice polytope} \right\rbrace
$$
are in one-to-one correspondence by the inverse maps 
$$
(X_\Sigma, D) \mapsto P_D \quad \text{and} \quad P \mapsto (X_P, D_P),  
$$
where $D_P$ is the torus invariant divisor 
$\ \sum_{\rho \in \Sigma_P(1)} \, (- \min_{m \in P} \langle m, u_\rho \rangle) \, D_\rho$. 
\end{theo}
 
\subsubsection{Polytopes, morphisms and basepoint free divisors}\label{Polbpf}

When we have a couple $(X_\Sigma, D)$ where $D=\sum_{\rho \in \Sigma(1)} a_\rho D_\rho$ is a torus-invariant Cartier divisor which is not ample but only nef (which in this setting is equivalent to basepoint free), applying the inverse maps of Theorem \ref{invmaps} we get a new couple $(X_{\bar{\Sigma}}, \bar{D})$ where
$X_{\bar{\Sigma}}$ is a projective variety and $\bar{D}$ an ample 
divisor on it, together with a toric morphism $\phi : X_{\Sigma} \longrightarrow X_{\bar{\Sigma}}$ such that $D$ is linearly equivalent to the pullback $\phi^*\bar{D}$. 
Let us give more details. 
\bigskip

A toric morphism $\varphi : X_{\Sigma_1} \rightarrow X_{\Sigma_2}$ is a morphism between toric varieties whose restriction to the open tori is a morphism of algebraic groups. In particular to such a morphism corresponds a $\ZZ$-linear map $\bar{\varphi} : N_1 \rightarrow N_2$ that is compatible with the fans $\Sigma_1$ and $\Sigma_2$, i.e. such that for any cone $\sigma_1 \in \Sigma_1$, there is a cone $\sigma_2 \in \Sigma_2$ containing the image $\bar{\varphi}(\sigma_1)$.
\smallskip

As stated in \eqref{Dbpf} and \eqref{Dample}, when the torus-invariant Cartier divisor $D$ is nef but not ample, the elements of the Cartier data $\lbrace m_\sigma \mid \sigma \in \Sigma(n) \rbrace$ are still the vertices of the polytope $P_D$ but they are not distinct any more.
In this case the maximal cones of the normal fan $\Sigma_{P_D}$ are unions of maximal cones of $\Sigma$ (\cite[Prop. 6.2.5]{CLS})~:
$$
C_v^\vee = \bigcup_{\substack{\sigma\in \Sigma(n) \\ 
m_\sigma=v}} \sigma. 
$$
It may happen (when $D$ is not big) that the lattice polytope $P_D$ fail to be full dimensional which means that the $C_v^\vee$ are not strongly convex cones. In this case we say that $\Sigma_{P_D}$ is a degenerate fan. 
All the cones of $\Sigma_{P_D}$ contain a linear space $U_0$ called the minimal cone of the degenerate fan (it is the orthogonal of the span of $P_D$).

If we put $N_0=U_0 \cap N$ then the linear map 
$$
\bar{\phi} : N \longrightarrow \bar{N}=N/N_0
$$
is compatible with the fan $\Sigma \subset N_\RR$ and the nondegenerate fan
\begin{equation}\label{Sigmabar}
\bar{\Sigma} = \lbrace \overline{\sigma}:= \sigma / U_0  \mid 
\sigma \in \Sigma_{P_D} \rbrace 
\subset \bar{N}_\RR=N_\RR/ U_0.
\end{equation}
\medskip

In other words the fan $\Sigma$ refines the possibly degenerate fan $\Sigma_{P_D}$ and hence we have a toric morphism between $X_\Sigma$ and the projective toric variety associated to $\Sigma_{P_D}$ (after having got rid of the degeneracy).   

Here is the precise result :
\begin{theo}[\cite{CLS} Th. 6.2.8]\label{phi}
Let $D$ be a basepoint free Cartier divisor on a complete toric variety $X_{\Sigma}$, and let $X_{\bar{\Sigma}}$ be the toric variety of the fan $\bar{\Sigma} \subset \bar{N}=N/N_0$ defined by \eqref{Sigmabar}. 
Then the linear map $\bar{\phi} : N \longrightarrow \bar{N}$ induces a proper toric morphism
$$
\phi : X_{\Sigma} \longrightarrow X_{\bar{\Sigma}}.
$$
Furthermore, $X_{\bar{\Sigma}}$ is a projective toric variety and $D$ is linearly equivalent to the pullback of an ample divisor $\bar{D}$ on $X_{\bar{\Sigma}}$.
\end{theo}
\bigskip

In other words, if $D$ is basepoint free, the image of the morphism associated to the complete linear system $|D|$ is the toric variety defined by the simplicial fan $\bar{\Sigma}$ in $\bar{N}_\RR$.

Note that the couple $(X_{\bar{\Sigma}}, \bar{D})$ is the one given by the second map of Theorem \ref{invmaps} applied to $P_D$ seen as a full dimensional lattice polytope in its real span.

\subsubsection{Newton polytopes}\label{Newton}

Let us first remark that the fans  $\Sigma_{P_D}$ and $\bar{\Sigma}$ do not depend on the precise torus-invariant divisor $D$, but only on its linear equivalence class $[D]$. Indeed, starting with another torus invariant divisor $D' \in [D]$ in the same class would give a polytope $P_{D'}$ equal to the translate of $P_{D}$ by the only character $m' \in M$ such that 
$$
D-D'=\div(m')= \sum_{\rho \in \Sigma(1)} \left\langle m', u_{\rho} \right\rangle D_\rho,
$$  
and hence give the same normal fan $\Sigma_{P_{D'}}=\Sigma_{P_{D}}$. The same is true for the nondegenerate fan $\bar{\Sigma}$ and since they do not depend on the choice of a representative of $[D]$, we will denote them by 
$$
\Sigma_{[D]} = \Sigma_{P_{D}} \quad \text{and} \quad \bar{\Sigma}_{[D]}=\bar{\Sigma}. 
$$
This explains why in Theorem \ref{phi} the divisor $D$ is not supposed to be torus-invariant~: 
we can start with any nef Cartier divisor $D$ up to choose a torus-invariant one $D' \in [D]$ in the same class (which is always possible, see next subsection).

The following result shows that we can even define a ``canonical model" of the polytopes $P_{D'}$ for $D' \in [D]$, depending only on the class $[D]$. 

\begin{lemdef} \label{hD'}
Let $D \in \Div(X_{\Sigma})$ be a Weil divisor and
$D'=\sum_{\rho \in \Sigma(1)} a'_{\rho} D_{\rho} \in [D]$ be a torus invariant divisor in the class of $D$. The map
$$
h_{D'} : 
M_{\RR} \longrightarrow \RR^{\Sigma(1)} , \quad
      m  \longmapsto      \left( \langle m, u_{\rho} \rangle + a'_{\rho} \right)_{\rho \in \Sigma(1)}
$$
embeds $M_{\RR}$ as an affine subspace of $\RR^{\Sigma(1)}$. 

The image $P_{[D]}:= h_{D'}(P_{D'})$ of the convex polytope $P_{D'}$ by this map   
is independent of the choice of $D'\in [D]$ and we have 
$$
P_{[D]} =h_{D'}(M_\RR) \cap \RR_+^{\Sigma(1)},
$$
\end{lemdef}

\begin{preuve}
It is Exercice 5.4.2 of \cite{CLS}.
The map $h_{D'}$ is injective because $\Sigma$ is full dimensional, so that
$$
\langle m, u_{\rho} \rangle=\langle m', u_{\rho} \rangle  \ \text{for all} \ \rho \in \Sigma(1)
\quad \Rightarrow m=m'.
$$

The equality comes from the equivalence
$$
m \in P_{D'} \Leftrightarrow
\forall  \rho \in \Sigma(1), \langle m, u_\rho \rangle \geq -a'_\rho \Leftrightarrow
\forall  \rho \in \Sigma(1), h_{D'}(m)_\rho \geq 0. 
$$ 

In order to prove the independence let us consider another torus invariant divisor $D''=\sum_{\rho \in \Sigma(1)} a''_{\rho} D_{\rho}$ in $[D]$. There exists a character $m''$ such that
$$
D'-D'' = \div(\chi^{m''}) 
= \sum_{\rho \in \Sigma(1)} \langle m'', u_\rho \rangle D_\rho,
$$
that is, $a'_{\rho}=a''_{\rho}+\langle m'', u_\rho \rangle$ for all $ \rho \in \Sigma(1)$. 
In particular for all $m \in M$ we have 
$$
h_{D''}(m)=\left( \langle m, u_{\rho} \rangle + a''_{\rho} \right)_{\rho \in \Sigma(1)}
=\left( \langle m-m'', u_{\rho} \rangle + a'_{\rho} \right)_{\rho \in \Sigma(1)}
=h_{D'}(m-m'')
$$
and 
\begin{equation*}
\begin{split}
m \in P_{D''} & \Leftrightarrow
\forall  \rho \in \Sigma(1), \langle m, u_\rho \rangle \geq -a''_\rho \\ 
&\Leftrightarrow  \forall  \rho \in \Sigma(1), \langle m-m'', u_\rho \rangle \geq -a'_\rho  \ \ \Leftrightarrow \, m-m'' \in P_{D'} 
\end{split}
\end{equation*} 
so that 
\begin{equation*}
\begin{aligned}
h_{D''}(P_{D''}) =& \lbrace h_{D''}(m) \mid m \in P_{D''}\rbrace 
= \lbrace h_{D'}(m-m'') \mid m-m'' \in P_{D'}\rbrace \\
= &\lbrace h_{D'}(m') \mid m' \in P_{D'}\rbrace = h_{D'}(P_{D'}).
\end{aligned}
\end{equation*} 
\end{preuve}

Let us now make the link with homogeneous coordinates on $X_\Sigma$. 
In \cite[§ 5.4]{CLS} the authors call $D'$-homogenization the process associating to a given character $\chi^m$ of the torus $T \subset X_\Sigma$, the character of the big torus $T' \subset \AA^{\Sigma(1)}$ (or the Laurent monomial)  
$$
x^{h_{D'}(m)} = \prod_{\rho \in \Sigma(1)} x_\rho^{\langle m, u_{\rho} \rangle + a_{\rho}} 
\in K\bigl(x_\rho, \rho \in \Sigma(1)\bigr)
$$
where $D'=\sum_{\rho \in \Sigma(1)} a_{\rho} D_{\rho}$.

This provides a nice description of the set of monomials that may appear in the homogeneous polynomial defining an effective divisor on $X_\Sigma$. 
Let us introduce some notations~:

\begin{nota}
Let $D \in \Div(X_{\Sigma})$ be an effective Weil divisor. We put 
$$
L_{[D]}=\biggl\lbrace (a_{\rho})_{\rho \in \Sigma(1)} \in \ZZ^{\Sigma(1)} \mathrel{\bigg|}
 \sum_{\rho \in \Sigma(1)} a_{\rho} D_{\rho} \in {[D]} \biggr\rbrace
\quad \text{and} \quad \AAA_{[D]}=\NN^{\Sigma(1)} \cap L_{[D]},
$$
so that we have 
$$
S_{[D]} = \bigoplus_{a \in \AAA_{[D]}} K x^a,
$$
where $x^a$ denotes the monomial $\;\prod_{\rho \in \Sigma(1)} x_{\rho}^{a_{\rho}}$ for $a=(a_{\rho})_{\rho \in \Sigma(1)} \in \NN^{\Sigma(1)}$.
\end{nota}

The crucial fact is that the monomials generating the graded piece $S_{[D]}$ 
are exactly the $D'$-homogenizations of the characters belonging to $P_{D'}$ for any torus-invariant divisor $D'$ in the class $[D]$ (see \cite[Prop. 5.4.1]{CLS}). 
In other words, for any torus-invariant $D' \in [D]$, the image of the map $h_{D'}$ is the affine hull of the lattice $L_{[D]}$. This means that the elements of $\AAA_{[D]}$ are exactly the points of $L_{[D]}$ contained in $P_{[D]}$~:
$\AAA_{[D]} = P_{[D]} \cap L_{[D]}$.
In particular, for any homogeneous polynomial $f  \in S_{[D]}$, the Newton polytope $\Delta_f$ of $f$ is contained in the polytope $P_{[D]}$.
\bigskip

Finally let us assume that $D$ is a nef and Cartier divisor. 
For any torus-invariant $D'\in[D]$ we have $P_D=\Conv( m_\sigma(D), \sigma \in \Sigma(n))$ in $M_\RR$. Hence in $\RR^{\Sigma(1)}$ we have 
$$
P_{[D]}=\Conv( a^\sigma([D]), \sigma \in \Sigma(n) ), \ \ \text{where} \ 
a^\sigma([D])=h_{D'}(m_\sigma(D')).
$$ 
In particular each $a^\sigma([D])$ belongs to $\AAA_{[D]}$ and $P_{[D]}$ is a lattice polytope : it is the Newton polytope of the general member of the graded piece $S_{[D]}$. 
Notice that for each $\sigma \in \Sigma(n)$, the vertex $a^\sigma([D])=(a^\sigma_{\rho})_{\rho \in \Sigma(1)}$ of $P_{[D]}$ is the coefficient vectors of the local representative of $[D]$ 
$$
D^\sigma =\sum_{\rho \in \Sigma(1)} a^\sigma_{\rho} D_\rho. 
$$   
In particular we have $a^\sigma_{\rho}=0$ for all $\rho \in \sigma(1)$ and $a^\sigma([D])$ is the unique element of $\AAA_{[D]}$ satisfying these vanishings.

\subsection{Toric Mori theory}\label{TorMor}

In this subsection we recall some of the tools from toric birational geometry that are essential for our results.   

\subsubsection{Toric Minimal Model Program}\label{TorMMP}

The aim of the Minimal Model Program (MMP for short) is the classification of projective algebraic varieties up to birational equivalence, with a special interest given to some special representatives in each class (the minimal models).

This program has been carried out for surfaces by the ``Italian school'' of Algebraic Geometry during the first half of the 20\textsuperscript{th} century, yielding the Enriques classification of algebaric surfaces. 

In higher dimensions it is also called Mori program and is a very active area of research since more than three decades. 
It relies in an essential manner on the properties of the Mori cone, that is, the cone of numerical classes of effective curves. 
More specifically, a crucial fact for a $\QQ$-factorial projective variety $X$, is that the $K_X$-negative half of the Mori cone $\NE(X)$ (i.e. the set of effective classes $[C]$ such that $C \cdot K_X<0$) is polyhedral. 
This allows to perform contractions of the rational curves contained in extremal rays of this polyhedral part, with the aim of simplifying the variety towards a minimal model. 

Repeating this yields a step by step procedure, again called MMP, consisting in a finite number of divisorial contractions or flips (see below) until the dimension of the variety drops, or until there is no more curve in the $K_X$-negative part of $\NE(X)$ (that is, the canonical divisor $K_X$ is nef).

An extremal contraction $\psi : X \to X'$ appearing in the MMP can be of three different types~:

\begin{enumerate}
\item[(1)] A divisorial contraction. It contracts a divisor to a lower dimensional subvariety, making the Picard number to drop by 1. It is the converse of a weighted blow-up. In this case we replace $X$ by $X'$ and continue the process.
\item[(2)] A flipping (or small) contraction. It contracts a subvariety of codimension $\geq2$, giving rise to a non $\QQ$-factorial singularity in $X'$. Since such a singularity is usually considered too bad, we take the flip map $\nu : X \dashrightarrow X^+$ which is an isomorphism in codimension 1 and continue with the $\QQ$-factorial variety $X^+$. A typical example of a flip is the composition of a blow-up and a blow-down in ``another direction''. 
\item[(3)] A fibering contraction. It contracts the whole variety onto a variety of smaller dimension. If $\psi$ is fibering we consider $X$ as minimal and stop the process.
\end{enumerate}
Note that the existence of the flip $\nu$ is in general hard to prove (see \cite{BCHM}) but easy for simplicial toric varieties. 
\medskip

The choice of the extremal rays to contract can obey to various constraints :

\begin{enumerate}
\item[-] Instead of contracting the rays that are $K_X$-negative, we can choose a so-called boundary $\QQ$-divisor $B$ and contract the rays that are $(K_X+B)$-negative. 
This is called a logarithmic MMP (log-MMP for short).
\item[-] Instead of taking the negative rays in the whole Mori cone $\NE(X)$,  we can restrict the choices to the rays in a proper face $\NE(X/Y)$ of $\NE(X)$ consisting of classes of curves contained in fibers of a given morphism $\phi : X \to Y$. 
The cone $\NE(X/Y)$ is then called the relative Mori cone and we say that we run the MMP relative to $Y$.
\end{enumerate}
These two restrictions on the choice of rays can be combined to give a relative log-MMP. Such a procedure stops either when it gets to a Mori fiber space or when the image of the divisor $(K_X+B)$ becomes $\phi$-nef. 
An instance of such a procedure is used in the proof of Proposition \ref{descent}.

Note that the canonical class of a toric variety is the class of the torus invariant divisor 
$K_{X_\Sigma} = - \sum_{\rho \in \Sigma(1)} D_\rho$ which is never effective, and a fortiori never nef.
In particular every classical MMP ends on a fibering contraction. The corresponding minimal model is called a Mori fiber space.  
\bigskip

In contrast with the general case, which yields extremely difficult issues, the toric MMP is quite simple to describe and it is relatively easy to prove that all its instances described above work for any simplicial projective toric variety (see \cite[Ch. 14]{Matsuki} for a proof or \cite[Ch. 15]{CLS} for the non-relative case).
 
Indeed for projective toric varieties, it is not just one half but the whole Mori cone which is polyhedral. 
This is one of the manifestations of the fact that projective toric varieties are Mori Dream Spaces (see \cite{HK00} for an introduction). 
In particular the MMP works particularly well for such varieties. 
Another manifestation of this fact is the existence of homogeneous coordinates, that is, the finite generation of the Cox ring (see paragraph \ref{Homcoor}). 
The link between the two relies on the theory of variation of geometric invariant theory (VGIT) and can be understood via the study of the secondary fan (see \cite[Ch. 14 and 15]{CLS}). 
This is a very beautiful subject that provides a convex-geometric insight on the MMP but we do not address it here, for lack of space reasons.
\medskip 

Let us begin with the cone theorem.

\begin{theo}
Let $X_{\Sigma}$ be a complete toric variety. The Mori cone $\NE(X_{\Sigma})$ is generated by classes of toric curves : 
$$
\NE(X_{\Sigma})=\operatorname{NE}(X_{\Sigma})=\sum_{\tau \in \Sigma(n-1)} \RR_{+} [C_{\tau}].
$$
In particular, it is a rational polyhedral cone in $N_1(X_{\Sigma})$. 
\end{theo}

The fact that extremal rays are generated by classes of toric curves allow a very explicit description of the steps of a toric MMP in terms of the fan of the variety.
Let $X_\Sigma$ be a simplicial projective variety of dimension $n$ and $\psi : X_\Sigma \rightarrow X_{\Sigma'}$ be the contraction of the extremal ray $\RRR$ in $\NE(X_\Sigma)$. 
The fan defining $X_{\Sigma'}$ is obtained from $\Sigma$ by removing all the walls $\tau \in \Sigma(n-1)$ such that $[C_\tau] \in \RRR$. 
In order to be more precise, we first need to describe some properties of these ``extremal walls''. 

For any $\tau \in \Sigma(n-1)$ such that $[C_\tau] \in \RRR$, there are exactly two $n$-dimensional  cones $\sigma_1, \sigma_2$ such that $\tau=\sigma_1 \cap \sigma_2$. Let us put 
$$
\Delta(\tau):=\sigma_1+\sigma_2 = \sum_{\substack{\sigma \in \Sigma(n) \\ \tau \prec \sigma}} \sigma \, ,
$$
and for any $\rho \in \Sigma(1)$ such that $\rho \subset \Delta(\tau)$ let us put 
$$
\Delta_\rho(\tau):=\Cone(u_{\rho'} \mid \rho' \subset \Delta(\tau), \rho' \neq \rho ).
$$
We have the following properties~:
 
\begin{lem}[\cite{Reid}]\label{Reid}
Let $X_{\Sigma}$ be a simplicial projective toric variety of dimension $n$ and $\RRR$ be an extremal ray of $\NE(X_\Sigma)$.  
For all $\tau \in \Sigma(n-1)$ such that $[C_\tau] \in \RRR$ we have
$$
\bigcup_{\rho \in J_{C_\tau}^+}\Delta_\rho(\tau)=\Delta(\tau)=\bigcup_{\rho \in J_{C_\tau}^-}\Delta_\rho(\tau)
$$ 
and 
$$
\Delta_\rho(\tau) \in \Sigma(n) \quad \text{for all } \rho \in J_{C_\tau}^+ .
$$
\end{lem}

Note that for two toric curve classes $[C_\tau],[C_{\tau'}] \in \RRR$ the intersection numbers $C_{\tau} \cdot D_\rho$ and $C_{\tau'} \cdot D_\rho$ have the same sign for all $\rho \in \Sigma(1)$. In particular the sets $J_{C_{\tau'}}^+ =J_{C_\tau}^+ $ and $J_{C_{\tau'}}^- =J_{C_\tau}^-$ depend only on the ray $\RRR$, hence we can denote them by $J_{\RRR}^+$ and $J_{\RRR}^-$.
\smallskip

Here is a combinatorial description of the different steps of the toric MMP~:

\begin{enumerate}
\item[(1)] If $J_\RRR^-=\{\rho_0 \}$ then $\psi$ is a divisorial contraction and $D_{\rho_0}$ is the divisor in $X_\Sigma$ contracted by $\psi$. In this case we have $\Sigma'(1) = \Sigma(1)\setminus \{\rho_0 \}$  
and 
$$
\Sigma'(n) = \Sigma(n) \setminus \lbrace \sigma \in \Sigma(n) \mid \rho_0 \in \sigma(1) \rbrace
\cup \lbrace \Delta(\tau) \mid [C_\tau] \in \RRR \rbrace.
$$
The fan $\Sigma$ is called the star subdivision of $\Sigma'$ at $\rho_0$.
\item[(2)] If $\Card(J_\RRR^-)\geq 2$ then $\psi$ is a small contraction. 
If $\nu : X_\Sigma \dashrightarrow X_{\Sigma^+}$ is the corresponding flip, then we have $\Sigma^+(1)=\Sigma'(1) = \Sigma(1)$, and 
$$
\Sigma^+(n)= \Sigma(n) \setminus
\lbrace \Delta_\rho(\tau) \mid [C_\tau] \in \RRR, \rho \in J_{C_\tau}^+ \rbrace
\cup \lbrace \Delta_\rho(\tau) \mid [C_\tau] \in \RRR, \rho \in J_{C_\tau}^- \rbrace. 
$$
Notice that if we denote by $\rho_0$ the one dimensional cone of $N_\RR$ generated by 
$$
\sum_{\rho \in J_\RRR^+} C_\tau \cdot D_\rho \: u_\rho = \sum_{\rho \in J_\RRR^-} (-C_\tau \cdot D_\rho) u_\rho
$$ 
then we have 
$$
\rho_0=\bigcap_{\substack{\rho \in J_\RRR \\ [C_\tau] \in \RRR}} \Delta_\rho(\tau)
$$
and the fans obtained by star subdividing $\Sigma$ and $\Sigma'$ at $\rho_0$ are the same. In particular the flip $\nu$ can be factored as a weighted blow-up followed by a divisorial contraction.

\item[(3)] If $J_\RRR^- = \emptyset$ then $\psi$ is a fibering contraction. 
In this case we have
$$
\Sigma'(n) =  \lbrace \Delta(\tau) \mid [C_\tau] \in \RRR \rbrace,
$$
which is a degenerate fan.
Indeed the cones $ \Delta(\tau)$ for  $[C_\tau] \in \RRR$ are not strongly convex for they all contain the linear space 
$$
U_0=\Cone(u_\rho, \rho \in J_\RRR^+)=\Span(u_\rho, \rho \in J_\RRR^+)
$$ 
As in Paragraph \ref{Polytopes}, we get a nondegenerate fan $\bar{\Sigma'}$ by taking the quotient by $U_0$. 
In particular the linear map $\bar{\psi}$ associated to $\psi$ is the quotient of real vector spaces $N_\RR \to N_\RR / U_0$ 
and the variety $X_{\Sigma'}=X_{\bar{\Sigma'}}$ has dimension $n'=n-n_0<n$ where $n_0=\dim(U_0)$. 

The fiber bundle structure of $X_{\Sigma}$ can be read from the fan $\Sigma$ itself~: there exist two subfans $\widehat{\Sigma}$ and $\Sigma_0$ of $\Sigma$ such that $\Sigma_0=\lbrace \sigma \in \Sigma \mid \sigma \subset U_0 \rbrace$  and $\bar{\psi}$ induces a bijection from $\widehat{\Sigma}$ to $\bar{\Sigma}_1$ such that we have
$$
\Sigma = \lbrace \sigma_0 + \widehat{\sigma} \mid \sigma_0 \in \Sigma_0, \widehat{\sigma} \in \widehat{\Sigma} \rbrace.
$$
In particular it follows from the description of the fan of a toric subvariety given in Paragraph \ref{InterSubvar} that for every maximal cone $\sigma_0 \in \Sigma_0(n_0)$ the toric subvariety $\V(\sigma_0)$ is isomorphic to $X_{\bar{\Sigma'}}$ and the inclusion $\V(\sigma_0) \hookrightarrow X_{\Sigma}$ induces a section $X_{\bar{\Sigma'}} \hookrightarrow X_{\Sigma}$ of the fibration $\psi$ whose image is $\V(\sigma_0)$.

Similarly, for every maximal cone $\widehat{\sigma} \in \widehat{\Sigma}(n')$ the toric subvariety $\V(\widehat{\sigma})$ is a fiber of $\psi$, isomorphic to $X_{\Sigma_0}$. 
It follows from Lemma \ref{Reid} that $\Sigma_0(1)=J_\RRR^+$ has cardinality $1+n_0$, that is, $X_{\Sigma_0}$ is a complete simplicial toric variety of Picard number 1.

\end{enumerate}
\medskip

In the divisorial and small cases the contracted locus in $X_\Sigma$ is the subvariety \mbox{$V_-=\V(\gamma_-)$} where $\gamma_-=\Cone(u_\rho, \rho \in J_\RRR^-) \in \Sigma$ and its image  by $\psi$ is $V_+=\V(\gamma'_+) \subset X_{\Sigma'}$ where $\gamma'_+=\Cone(u_\rho, \rho \in J_\RRR^+) \in \Sigma'$.
The restriction 
$\psi_{|V_-}: V_- \to V'_+$ is a fibering extremal contraction. 
In particular there exist a complete simplicial toric subvariety  $\V(\widehat{\sigma_-}) \subseteq V_- \subset X_\Sigma$ of Picard number 1, which is a fiber of $\psi$.  

As pointed by K. Matsuki in Remark 14-2-3 of \cite{Matsuki}, these subvarieties need not be weighted projective spaces in general and may instead be quotient of those by a finite abelian group, often called fake weighted projective spaces.

This leads to the following definition~:

\begin{defi}
Well formed weighted projective spaces together with fake weighted projective spaces are exactly the simplicial complete toric varieties of Picard number one. We call them \textit{generalized weighted projective spaces}. 
\end{defi} 

The fact that non trivial fibers of toric extremal contractions are generalized weighted projective spaces play a determinant role in the interplay between the different notions of low toric degree.
 
\subsubsection{Cones of positive numerical classes } 
 
In addition to the Mori cone $\NE(X_\Sigma)$ who plays a central role in the MMP, several cones appear naturally when studying the birational geometry of a projective variety.
First the cone of nef divisors, or nef cone, $\Nef(X_\Sigma)$ plays also an important role because it is the dual cone of $\NE(X_\Sigma)$ : nef divisors are those divisors $D$ having nonnegative intersection number with every effective curve.

It is worth noticing that the contraction of an extremal ray $\RRR$ of $\NE(X_\Sigma)$ is a particular case of the proper toric morphism $\phi : X_{\Sigma} \longrightarrow X_{\bar{\Sigma}_{[D]}}$ of Theorem \ref{phi} : the one associated to any nef divisor $D$ such that the hyperplane 
$D^\bot:=\lbrace [C] \in N_1(X_\Sigma) \mid C \cdot D =0 \rbrace$ intersects $\NE(X_\Sigma)$ exactly in $\RRR$ (such a divisor exists precisely because $\RRR$ is extremal).

Another interesting couple of dual cones of numerical classes is formed by the effective cone $\Eff(X_\Sigma)$ generated by divisor classes $[D]$ containing at least one effective representative and the cone of mobile curves, or mobile cone, $\Mob(X_\Sigma)$ generated by curve classes $[C]$ whose members move in families covering the whole variety $X_\Sigma$.

Since the effective cone is generated by the classes of toric divisors $D_\rho$, mobile curves are those curves $C$ on $X_\Sigma$ verifying $J_C^-=\emptyset$ and $J_C=J_C^+$.

The generators of the mobile cone have been characterized by D. Mons\^{o}res in his Ph.D thesis \cite{Monsores} : they are the mobile classes $[C]$ such that $\Card(J_C)=\dim(\Span(J_C)) + 1$, and are numerical pullbacks of mobile and extremal curves on a minimal model of $X_\Sigma$. 

\section{First results} \label{ResEx}

In this section we show that a Weil divisor $D$ of a simplicial complete toric variety $X_\Sigma$ has automatically a rational point over $\ZZ$ if the Newton polytope of its defining homogeneous polynomial $f \in S_{[D]}$ is too small, which is in  particular the case of divisors that are not Cartier and nef.
Then we give a first criterion of existence of a rational point based on Theorem~\ref{Kollar}~: 
the divisor $D$ admits a $K$-point as soon as its restriction $D_{|W}$ has lower degree than the anticanonical divisor $-K_W$ of a toric subvariety $W \subseteq X_\Sigma$ isomorphic to a generalized weighted projective space.

\subsection{Trivial rational points} \label{Triv}

It is a trivial fact that a hypersurface of degree $d$ of the projective space $\PP^n$ whose homogeneous equation does not feature the monomial $x_0^d$, automatically contains the $\ZZ$-point $(1 : 0 : \cdots : 0)$. 
More generally, as soon as the Newton polytope of a homogeneous polynomial does not contain all the monomials of the form $x_i^d$, it admits trivial integer solutions. 
This phenomenon also occurs for a hypersurface $D$ of a complete simplicial toric variety  $X_\Sigma$ endowed with homogeneous coordinate (i.e. satisfying hypothesis (H3), see \ref{Homcoor}) with the role of the special points $(0 : \cdots : 0 : 1 : 0 : \cdots : 0)$ played by the distinguished points 
$$
\gamma_\sigma := \bigcap_{\rho \in \sigma(1)} D_{\rho} \quad \text{for} \;
\sigma \in \Sigma(n).
$$
and the role of the extremal monomials $x_i^d$ played by the $D'$-homogenization of the Cartier data $\lbrace m_{\sigma}(D') \in M_\RR \mid  \sigma \in \Sigma(n) \rbrace$ for any torus-invariant divisor $D'$ linearly equivalent to $D$.

Here is the precise result

\begin{lem}\label{smallNewton} 
Let $X_{\Sigma}$ be a simplicial complete toric split variety satisfying hypothesis (H3) and $D$ be an effective divisor on $X_{\Sigma}$, zero locus of the homogeneous polynomial $f \in S_{[D]}$. 
Let $D'=\sum_{\rho \in \Sigma(1)} a'_{\rho} D_{\rho} \in [D]$ be a torus-invariant divisor linearly equivalent to $D$.
\smallskip

If the Newton polytope of $f$ does not contain the element
$$
a^\sigma([D])=h_{D'}(m_\sigma(D'))=\left( \langle m_\sigma(D'), u_{\rho} \rangle + a'_{\rho} \right)_{\rho \in \Sigma(1)} \in \RR^{\Sigma(1)}
$$ 
for some $\sigma \in \Sigma(n)$, then $D$ contains the corresponding torus-invariant point $\gamma_\sigma = \bigcap_{\rho \in \sigma(1)} D_{\rho}$.
\end{lem}

\begin{preuve}
Let us write
$$
f=\sum_{a \in \AAA_{[D]}} \lambda_a \prod_{\rho} x_{\rho}^{a_{\rho}} \in S_{[D]},
$$
and let us denote by $\AAA_f$ the set of monomials in $\AAA_{[D]}$ such that $\lambda_a \neq 0$. By definition the Newton polytope of $f$ is the convex hull $\Delta_f$ of $\AAA_f$ in  $\RR^{\Sigma(1)}$.
Now let us suppose that there exists $\sigma \in \Sigma$ such that $a^\sigma([D])$ does not belong to $\Delta_f$. 
In particular we have $a^\sigma([D]) \notin \AAA_f$. 

Since $a^\sigma([D])$ is the only element of $\AAA_{[D]}$ satisfying $a^\sigma_{\rho}=0$ for all $\rho \in \sigma(1)$, 
for all $a \in \AAA_f$ there exists at least one $\rho \in \sigma(1)$ with $a_{\rho} \neq 0$.

It follows that putting $\alpha_{\rho}=0$ for all $\rho \in \sigma(1)$ and $\alpha_{\rho}=1$ for all $\rho \not\in\sigma(1)$ gives a root $\alpha \in \lbrace 0,1 \rbrace^{\Sigma(1)}$ of $f$ and since $\pi(\alpha) = \gamma_\sigma$ this gives the result.
\end{preuve}

The homogenization of the Cartier data $\lbrace a^\sigma([D]) \mid \sigma \in \Sigma(n) \rbrace$ is a set of lattice points in $\ZZ^{\Sigma(1)}$ if and only if $D$ is a Cartier divisor and it is contained in the polytope $P_{[D]}$ if and only if $D$ is a nef divisor. The consequence is the following

\begin{prop}\label{trivial} 
Let $X_{\Sigma}$ be a simplicial complete split toric variety and $D$ be an effective divisor on $X_{\Sigma}$.
If one of the following conditions is satisfied :
\begin{enumerate}
\item $D$ is not a nef divisor, 
\item $X_\Sigma$ satisfies hypothesis (H3) and $D$ is not a Cartier divisor,
\end{enumerate}
then $D$ contains a torus invariant point
$
\gamma_\sigma = \bigcap_{\rho \in \sigma(1)} D_{\rho} \quad \text{for some} \;
\sigma \in \Sigma(n),
$
and hence a rational point over $\ZZ$. 
\end{prop}

\begin{preuve}
If $D$ is not nef then there exists a toric curve $C_\tau$, $\tau \in \Sigma(n-1)$, such that $C_\tau \cdot D<0$. This implies that $C_\tau$ is included in $D$, and for all $\sigma \in \Sigma(n)$ containing $\tau$ we have $\gamma_\sigma \in C_\tau$ which gives the result. If (H3) is verified and $D$ is not Cartier then $D$ is the zero locus of some homogeneous polynomial $f \in S_{[D]}$ and the hypotheses of the previous lemma are verified. The proposition is proved.
\end{preuve}

In the sequel we assume $D$ to be a nef Cartier divisor.

\subsection{A simple criterion (Restricted low toric degree) } \label{RLTD}

In order to find rational points in a Cartier divisor $D$, a first simple idea is to consider the restriction $D_{|W}$ of $D$ to a weighted projective space $W$ contained in the toric ambient $X_{\Sigma}$ and apply Theorem \ref{Kollar}.   

For the convenience of the reader, we reproduce here Theorem \ref{Kollar} with its proof.
\begin{theorem*}[Th. 6.7 p.232 \cite{Kollar}]
Let $K$ be a $C_1$ field admitting normic forms of arbitrary degree. Let $K[x_0, \ldots, x_n]$ be the polynomial algebra graded by $\deg(x_i)=a_i \in \NN^*$. 
Let $f_1, \ldots, f_r \in K[x_0, \ldots, x_n]$ be homogeneous polynomials for this graduation. 
If
$$
\sum_{j=1}^{r} \deg(f_j) < \sum_{i=0}^{n} a_i,
$$
then the system of equations $f_1 = \cdots = f_r=0$ has a non trivial solution in $K^{n+1}$.
\end{theorem*}

\begin{preuve}
For $i=0, \ldots, n$ let $g_i(y_{i,1}, \ldots ,y_{i,a_i})$ be normic forms of degree $a_i$ over $K$. The $g_i$'s are homogeneous polynomials of the graded algebra $A=K[y_{i,k_i}, 0 \leq i \leq n, 1 \leq k_i \leq a_i]$ (with the usual grading coming from $\PP^{a_0 + \cdots + a_n-1}$). 

Now for $j=1, \ldots, r$ let us put
$$
F_j:=f_j(g_0, \ldots, g_n) \in K[y_{i,k_i}, 0 \leq i \leq n, 1 \leq k_i \leq a_i].
$$
By construction, $F_j$ is a homogeneous polynomial in $A$, of degree $\deg(f_j)$ and its number of variables $y_{i,k_i}$ is $a_0 + \cdots + a_n$.
Since $\sum_{j=1}^{r} \deg(f_j) < \sum_{i=0}^{n} a_i$, it follows by Lang's Theorem 4 in \cite{Lang} that the equation $F_1= \ldots = F_r=0$ has a solution $\left(b_{i,k_i}\right)$ such that $b_{s,t} \neq 0$ for at least one couple $(s,t)$. As $g_s$ is normic, $c_s:=g_s(b_{s,1}, \ldots, b_{s,a_s})$ is non zero. 
Hence putting
$$
c_i=g_i(b_{i,1}, \ldots, b_{i,a_i}) \quad i=0, \ldots, n
$$  
yields a non trivial solution of the equation $f_1= \ldots = f_r=0$. 
\end{preuve}

In fact Theorem \ref{Kollar} applies as well to generalized weighted projective spaces admitting homogeneous coordinates because in this case any non trivial solution automatically avoids the exceptional set (see \ref{Homcoor}). Here is the precise result~:

\begin{cor}\label{Kollarbis}
Let $K$ be a $C_1$ field satisfying hypothesis (H1) and $W$ be a generalized weighted projective space over $K$ satisfying hypotheses (H2) and (H3).
Let $H_1, \ldots, H_k \subset X$ be hypersurfaces in $W$. For any effective 1-cycle class $[C] \in \NE(W)$,  if 
$$
C \cdot (H_1 + \cdots + H_k) < C \cdot (-K_W),
$$
then the intersection $H_1 \cap \cdots \cap H_k$ has a rational point over $K$.
\end{cor}

By focusing on generalized weighted projective spaces that are toric subvarieties of the ambient $X_\Sigma$ we get an easily handled first notion of low toric degree for divisors that implies the existence of a rational point.

\begin{defi}\label{dltd}
Let $X_\Sigma$ be a normal toric variety. 
A Cartier divisor $D$ on $X_\Sigma$ is said to have restricted low toric degree if there exists a complete simplicial toric subvariety $W \subseteq X_\Sigma$ of Picard number one, such that the restriction $D_{|W}$ has low degree in the sense of Corollary~\ref{Kollarbis}~: there exists a 1-cycle $C$ on $W$ such that
$$
C \cdot D_{|W} < C \cdot (-K_W).
$$
\end{defi} 

\begin{Remarks}
\begin{enumerate}
\item This definition does not require the ambient toric variety $X_\Sigma$ to be simplicial nor projective, and not even to satisfy hypothesis (H3). 

\item A toric subvariety is an irreducible torus invariant subvariety $\V(\alpha), \alpha \in \Sigma$, we could have considered other subvarieties that are also toric varieties but the $\V(\alpha)$ are sufficient because of their link with the extremal classes (see \ref{TorMMP}).  
\end{enumerate}
\end{Remarks}

\begin{cor}
Let $K$ be a $C_1$ field satisfying hypothesis (H1) and
$X_\Sigma$ be a simplicial complete toric variety over $K$ satisfying hypothesis (H2) and such that any toric subvariety of Picard number one satisfies hypothesis (H3), for example $X_\Sigma$ smooth.

A restricted low degree divisor on $X_\Sigma$ has a rational point over $K$.
\end{cor} 

Notice that every toric curve $C_\tau$ is a normal separated toric variety of dimension 1 and hence is isomorphic to $\PP^1$. In particular a divisor $D$ has restricted low toric degree as soon as its intersection number with one of the toric curves is 1.  
More generally, on a simplicial variety $X_\Sigma$, a divisor $D$ has restricted low toric degree as soon as $C_\tau \cdot D \leq C_\tau \cdot \sum_{\rho \in J_{C_\tau} \setminus \gamma(1)} D_\rho$ where $\gamma$ is any subcone of $\tau$ such that $\V(\gamma)$ has Picard number 1. 
In particular when $C_\tau$ is an extremal curve, such a $\V(\gamma)$ is given by any toric fiber of the corresponding extremal contraction. 
It follows that a divisor has restricted low toric degree as soon as it has low degree along an extremal curve, in the precise sense of the following Lemma. 

\begin{lem}\label{extrgral}
Let $X_\Sigma$ be a simplicial projective toric variety.  
If there exist an extremal class $[C_\tau]$ on $X_\Sigma$ such that 
\begin{equation*}
0< C_\tau \cdot D <  \sum_{\rho \in J_{C_\tau}^{+}} C_\tau \cdot D_{\rho}
\end{equation*}
then $D$ has restricted low toric degree.
\end{lem}

\begin{preuve}
It follows from \ref{TorMMP} that for $[C_\tau]$ extremal and $\gamma= \bigcap_{\rho \in J_{C_\tau}^+} \Delta_\rho(\tau)$, we have $J_{C_\tau}^- \subseteq \gamma(1)$ $\subseteq \tau(1)$, $J_{C_\tau}^+ = \Adj(\gamma)$ and $W=\V(\gamma)$ is a fiber of the contraction associated to $[C_\tau]$, hence a generalized weighted projective space.  
The 1-cycle $C_\tau$ is thus the pushforward $j_*\bar{C}$ of a 1-cycle $\bar{C}$ on $W$ and using the projection formula, we have 
$$
\bar{C} \cdot D_{|W} = j_*\bar{C} \cdot D = C_\tau \cdot D.  
$$
Since this number is positive, $D_{|W}$ is a non trivial effective (and hence ample) divisor of $W$.
It is then sufficient to prove that $C \cdot D < \bar{C} \cdot (-K_W)$.
We recall from  \eqref{locglob3} that we have $\bar{C} \cdot (-K_W) \geq j_*\bar{C} \cdot \sum_{\rho \in \Adj(\gamma)} D_{\rho}$, so that by hypothesis  
$$
C_\tau \cdot D <  \sum_{\rho \in J_{C_\tau}^{+}} C_\tau \cdot D_{\rho} 
= \sum_{\rho \in \Adj(\gamma)} C_\tau \cdot D_{\rho}
\leq \bar{C} \cdot (-K_W).
$$
The lemma is proved.
\end{preuve}

One important property of the restricted low toric degree condition for $D$ is its stability through toric morphisms $\phi$ such that $\phi^*\phi_*[D]=[D]$. More precisely we have the following result. 

\begin{prop}\label{dltdbar}
Let $X_{\Sigma}$ be a simplicial projective toric variety
and $\phi : X_{\Sigma} \rightarrow X_{\bar{\Sigma}}$ be a contracting toric morphism. 
Let $\bar{D}$ be a Cartier divisor on $X_{\bar{\Sigma}}$.
If $D=\phi^*\bar{D}$ has restricted low toric degree then so has $\bar{D}$.
\end{prop}

\begin{preuve}
This is a consequence of the effect of contracting toric morphisms on generalized weighted projective spaces~: they are either contracted to a point or left unchanged.  
If $D$ has restricted low toric degree there exist a toric subvariety $W \subset X_\Sigma$ of Picard number 1 and a 1-cycle $C$ on $W$ such that
$$
0< C \cdot D_{|W} < C \cdot (-K_W).
$$  
If we denote by $\iota : W \hookrightarrow X_{\Sigma}$ the inclusion map, we have 
$$
\iota_*C \cdot D = C \cdot \iota^*D = C \cdot D_{|W} >0,
$$
and using the projection formula, this implies that 
$$
\phi_* \iota_*C \cdot \bar{D} = \iota_*C \cdot \phi^*\bar{D}  >0.
$$ 

It follows that $\iota_*C$ is not contracted by $\phi$ and that the restriction of $\phi$ to $W$ is an isomorphism
$$
\phi_{|W} : W \overset{\sim}{\rightarrow} \bar{W}.
$$ 
In particular in $\bar{W}$ we again have 
$$
0< \bar{C} \cdot \bar{D}_{|\bar{W}} < C \cdot (-K_{\bar{W}})
$$ 
and $\bar{D}$ has restricted low toric degree, which proves the proposition.
\end{preuve}

The main limit of the notion of restricted low toric degree is that it does not account for the specificities of non ample divisors, as illustrated by the following example.

\begin{ex}\label{L3}
The Losev-Manin space $\bar{L_3}$ is a del Pezzo surface of degree 6 obtained from $\PP^2$ by blowing-up the three torus-invariant points (or any three non aligned points up to projectivity).
Here are one possible fan for this variety together with the intersection table of its 6 toric divisors :

\begin{minipage}[c]{6cm}
\includegraphics[scale=0.8]{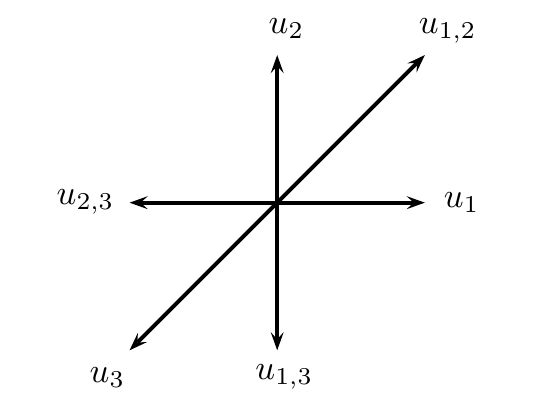} 
\end{minipage}
\hfill
\begin{minipage}[c]{8cm}
\vspace*{1cm}
$ \begin {array}{c|cccccc}  
\pmb{\cdot} & D_1 & D_{1,2} & D_2 & D_{2,3} & D_3 & D_{1,3}  \\
\hline 
D_{1} 	\quad &  -1&1&0&0&0&1  \\ 
D_{1,3}  \quad &  1&-1&1&0&0&0 \\  
D_{2}	 \quad &  0&1&-1&1&0&0 \\ 
D_{2,3}	 \quad &  0&0&1&-1&1&0 \\ 
D_{3}    \quad &  0&0&0&1&-1&1 \\
D_{1,3}	 \quad &  1&0&0&0&1&-1 \\
\end {array} $  
\vspace*{1.5cm}
\end{minipage}

\vspace{-0.3cm}
Let $C \subset \PP^2$ be a conic that contains none of the three blown-up points, and let $\phi : \bar{L_3} \to \PP^2$ be the contraction of the three exceptional curves.
Since $\phi$ is an isomorphism outside these three curves, the pullback $D=\phi^*C$ is sent isomorphically to $C$ by $\phi$ and in particular admits a rational point over any $C_1$ field.
Though, it has intersection number 2 with each strict transform of a toric line of $\PP^2$, and 0 with each exceptional curve. Indeed if $[L]$ denotes the class of the pullback of any line in $\PP^2$ we have $[D] = 2 [L]$ and 
$$
 [L] = [D_{1,3}+D_1+D_{1,2}] =[D_{1,3}+D_1+D_{1,2}] =[D_{1,3}+D_1+D_{1,2}] 
$$ 
$$
\Rightarrow  \  D_1 \cdot L = D_2 \cdot L = D_3 \cdot L = 1 \quad \text{and} 
\quad D_{1,2} \cdot L = D_{1,3} \cdot L = D_{2,3} \cdot L = 0.
$$

Since every toric subvariety of $\bar{L_3}$ with Picard number 1 is a $\PP^1$, this implies that $D$ is not a restricted low toric degree divisor.

Obviously one would expect a good notion of low toric degree to consider this $D$ as a low toric degree divisor. This motivates the search for a notion more stable by birational modifications.
\end{ex} 

\section{Low degree along a 1-cycle (Global low toric degree)}\label{LTD}

The goal of this section is to give some combinatorial conditions on the numerical class of the divisor $D$ implying the existence of rational points in $D$. Contrary to the results of last section, we try here to take the most possible advantage of the specificities of $D$. 
In particular we deepen the description of the links between the morphism $\phi : X_{\Sigma} \rightarrow X_{\bar{\Sigma}_{[D]}}$ of Theorem \ref{phi} and the Newton polytope $P_{[D]}$ of the homogeneous polynomial $f$ defining $D$, and use it to determine when a root of $f$ in $K$ gives a $K$-rational point. This yields a notion of low degree which may be lifted by toric birational contractions~: the notion of global low toric degree.

\subsection{Rational points in the total space}\label{RPTotal}

In order to generalize Lemma \ref{extrgral}, one can try to apply Corollary \ref{Kollarbis} from an inequality of the form
$$
0< C \cdot D <  \sum_{\rho \in I} C \cdot D_{\rho}
$$
with an arbitrary 1-cycle $C$. 
As soon as each $ D_{\rho}$ for $\rho \in I$ has positive degree along $C$ and there are at least two of them, this yields a root of the homogeneous polynomial defining $D$ but not a rational point in general.

\begin{lem}\label{fund}
Let $D$ be a nef Cartier divisor on $X_{\Sigma}$ and $f \in S_{[D]}$ a homogeneous polynomial such that $D=\V(f)$ as in \ref{Homcoor}. If there exist a 1-cycle $C$ and a subset $I \subseteq J_{C}^{+}$ with cardinal at least 2 such that
\begin{equation}\label{titdeg}
0< C \cdot D <  \sum_{\rho \in I} C \cdot D_{\rho}
\end{equation}
then $f$ has a root $\alpha \in K^{\Sigma(1)}$ such that
$$
\left\lbrace \rho \in \Sigma(1) \mid \alpha_{\rho} = 0 \right\rbrace 
\subseteq J_C \setminus \lbrace \rho_1 \rbrace 
$$ 
for some $\rho_1 \in I$.  
\end{lem}
  
\begin{preuve}
Since $f$ is homogeneous for the graduation by $\Cl(X_{\Sigma})$, it is a fortiori homogeneous for the graduation by $\ZZ$ defined by 
\begin{equation}\label{degC}
\deg(x_{\rho})=m C \cdot D_{\rho}
\end{equation}
where $m$ is any positive integer such that $m C \cdot D_{\rho} \in \ZZ$ for all $\rho \in \Sigma(1)$. Let $d\in \NN^*$ denote the degree of $f$ for this graduation.

Let us consider the polynomial $g \in S'=K[x_{\rho}, \rho \in I]$ defined from $f$ by putting $x_{\rho}=1$ for all $\rho \in \Sigma(1) \setminus J_{C}$ (for which $\deg(x_{\rho})=0$) and $x_{\rho}=0$ for all $\rho \in J_{C} \setminus I$. 

The polynomial $g$ is again homogeneous of degree $d$ for the graduation of $S'$ induced by the graduation of $S$ defined by \eqref{degC}, which is the graduation of a weighted projective space. 
Up to multiplying by $m$, we can thus apply Theorem \ref{Kollar} to inequality \eqref{titdeg} to get a non trivial zero $\alpha' \in K^{I}$ of $g$. This in turn yields a zero $\alpha \in K^{\Sigma(1)}$ of $f$ with $\alpha_{\rho}=\alpha'_{\rho}$  for $\rho \in I$, $\alpha_{\rho}=0$ for $\rho \in J_{C} \setminus I$ and $\alpha_{\rho}=1$ for $\rho \in \Sigma(1) \setminus J_{C}$. 
By construction, we have
$
\left\lbrace \rho \in \Sigma(1) \mid \alpha_{\rho} = 0 \right\rbrace 
\subset J_C  
$ 
and the fact that $\alpha'$ is not trivial exactly means that $\alpha_{\rho_1} \neq 0$ for some $\rho_1 \in I$, which terminates the proof. 
\end{preuve}
\bigskip

Lemma \ref{fund} does not automatically provide a $K$-rational point in $D$ since for a root 
\mbox{$\alpha \in K^{\Sigma(1)}$} of $f$ we have 
\begin{equation}\label{condB0}
\begin{array}{rcl}
\pi(\alpha) \in X_{\Sigma} & \ \Leftrightarrow \ & \alpha \notin Z(\Sigma) \\
  & \ \Leftrightarrow \ & \left\lbrace \rho \in \Sigma(1) \mid \alpha_{\rho} = 0 \right\rbrace  \subseteq \sigma(1) 
\ \text{for some cone} \ \sigma \in \Sigma. 
\end{array}
\end{equation}

If we require that for each $\rho \in I$ there exists a cone $\sigma \in \Sigma$ such that $J_C \setminus \{ \rho \} \subseteq \sigma(1)$, then obviously \eqref{condB0} is satisfied. But this turns out to be a too restrictive condition and, as shown in next subsection, we only need that each $J_C \setminus \{ \rho \}$ is contained in a cone of the (possibly degenerate) fan obtained from $\Sigma$ by removing all the walls $\tau \in \Sigma(n-1)$ such that $C_{\tau} \cdot D = 0$.

\subsection{The Newton polytope of basepoint free divisors}\label{NonAmple} 

As we have seen in the preliminaries, to a torus-invariant divisor $D=\sum_{\rho \in \Sigma(1)} a_\rho D_\rho$ can be associated a polytope $P_D=\lbrace m\in M_{\RR} \mid \left\langle m, u_{\rho} \right\rangle \geq -a_{\rho} \rbrace$ in the space of characters $M_{\RR}$.
To such a polytope are turn associated two objects that only depend on the numerical class $[D]$~: 
\begin{enumerate}
\item[-] The normal fan $\Sigma_{P_D}=\Sigma_{[D]}$ whose cones are in duality with the faces of $P_D$
\item[-] The ``canonical" polytope $P_{[D]}$, image of $P_D$ through the embedding \begin{equation}\label{hD} 
h_{D} : 
M_{\RR} \longrightarrow \RR^{\Sigma(1)} , \quad
      m  \longmapsto      \left( \langle m, u_{\rho} \rangle + a_{\rho} \right)_{\rho \in \Sigma(1)}
\end{equation}
If $D$ has no trivial rational point (see Subsection \ref{Triv}) then $P_{[D]}$ is naturally identified with the Newton polytope of the homogeneous polynomial defining $D$. 
In other words the lattice points of this polytope are the exponent vectors of the monomials belonging to the graded piece $S_{[D]}$ in the total coordinate ring. 
\end{enumerate}
 
Let us describe more precisely the relation between these two objects.
 
The correspondence between cones of the generalized fan $\Sigma_{P_D}=\Sigma_{[D]}$ and faces of the polytope $P_{D}$ is one-to-one and dimension reversing :
$$
\alpha' \in \Sigma_{P_{D}} \longleftrightarrow  
Q_{D}(\alpha')=
\left\lbrace m \in P_{D} \mid  \langle m, u_{\rho} \rangle = -a_\rho \ \text{for all} \  \rho \subseteq \alpha' \right\rbrace \subset P_D 
$$ 
Through the embedding $h_{D}$ of \eqref{hD}, this induces the contravariant correspondence 
\begin{equation}\label{corresp}
\begin{array}{rcl}
\lbrace \text{Cones of}\ \Sigma_{[D]} \rbrace & \longleftrightarrow & \lbrace \text{Faces of}\ P_{[D]} \rbrace  \\
\alpha'   & \longmapsto     & Q_{[D]} (\alpha') 
\end{array}
\end{equation}

where 
$$
Q_{[D]} (\alpha') = \lbrace (a_{\rho}) \in P_{[D]} \mid a_{\rho}=0 \ \text{for all} \ \rho \subseteq \alpha' \rbrace.
$$

Now let us assume that $D$ is Cartier and nef. 
In particular the Cartier data $\lbrace m_\sigma=m_\sigma(D) \mid \sigma \in \Sigma(n) \rbrace$ is the set of vertices of the polytope $P_D$ and it follows that the set of vertices of the image $P_{[D]}=h_D(P_D)$ is the set $\lbrace  a^{\sigma}=a^{\sigma}(D)=h_D(m_\sigma(D)) \in P_{[D]} \mid  \sigma \in \Sigma(n)  \rbrace$.
Moreover, since we have $m_{\sigma_1}=m_{\sigma_2}$ if and only if there exists $\sigma' \in \Sigma_{[D]}$ such that $\sigma_1 \cup \sigma_2 \subseteq \sigma'$, the set of vertices of $P_{[D]}$ is also the set  $\lbrace a^{\sigma'} \in P_{[D]} \mid  \sigma' \in \Sigma_{[D]}(n)  \rbrace$ where
\begin{equation}\label{asigma'}
\lbrace \rho \in \Sigma(1) \mid  a^{\sigma'}_{\rho}=0  \rbrace
=
\lbrace \rho \in \Sigma(1) \mid \rho \subseteq \sigma' \rbrace.
\end{equation} 
\bigskip

Since $\Sigma$ refines $\Sigma_{[D]}$, to each cone $\alpha$ of $\Sigma$ is naturally associated a cone of $\Sigma_{[D]}$ : the smallest that contains $\alpha$. 
We generalize this notion with the following definition :

\begin{defi}\label{barphi}
Let $\phi : X_{\Sigma} \longrightarrow X_{\bar{\Sigma}_{[D]}}$ be the toric morphism of Theorem \ref{phi}.
Let $c \subset \Sigma(1)$ be a set of 1-cones of $\Sigma$ all contained in a single cone $\alpha$ of $\Sigma_{[D]}$. We call $\phi$-hull of $c$ and denote by $\barphi{c}$ the smallest cone of $\Sigma_{[D]}$ containing $c$ :
$$
\barphi{c} := \bigcap_{\substack{\alpha \in \Sigma_{[D]} \\ c \subset \alpha}} \alpha \subset N_\RR. 
$$
\end{defi}

For any cone $\alpha \in \Sigma_{[D]}$, the face $Q_{[D]} (\alpha)$ of the polytope $P_{[D]}$ is by definition the intersection of $P_{[D]}$ with the coordinate hyperplanes corresponding to all the 1-dimensional cones $\rho \in \Sigma(1)$ contained in $\alpha$.
The following result shows that we can characterize the face $Q_{[D]} (\alpha)$ by intersections with the coordinate hyperplanes corresponding to any set $c$ of 1-dimensional cones $\rho \in \Sigma(1)$ such that 
$\barphi{c}=\alpha$~:

\begin{lemdef}\label{ouf}
Let $D$ be a nef Cartier divisor on $X_{\Sigma}$ and $c \subset \Sigma(1)$ be a set of 1-cones of 
$\Sigma$ all contained in a single cone of $\Sigma_{[D]}$. Then the set
$$
Q_{[D]}(c) := \lbrace (a_{\rho}) \in P_{[D]} \mid 
\forall \rho \in c, a_{\rho}=0  \rbrace
$$ 
is equal to the face $Q_{[D]} (\barphi{c})$ of $P_{[D]}$.
\end{lemdef}

\begin{Remark}
We use the same notation $Q_{[D]} (\cdot)$ for cones in $\Sigma_{[D]}$ and for sets of 1-dimensional cones in $\Sigma$. Of course for any $\alpha \in \Sigma_{[D]}$, by putting (with a slight abuse of notation)  
$\alpha(1)= \lbrace \rho \in \Sigma(1) \mid \rho \subseteq \alpha \rbrace$, we have $\alpha=\Cone(\rho, \rho \in \alpha(1))$ and by the above definition  $Q_{[D]} (\alpha(1))=Q_{[D]} (\alpha)$. What the lemma asserts is that we can in some cases choose a set $c$ strictly smaller than $\alpha(1)$ but still verifying $Q_{[D]} (c)=Q_{[D]} (\alpha)$. 
\end{Remark}

\begin{preuve}
Since $c \subseteq \barphi{c}$ we clearly have the inclusion
$$
Q_{[D]} (\barphi{c}) = \lbrace (a_{\rho}) \in P_{[D]} \mid 
 a_{\rho}=0 \ \text{for all} \ \rho \subseteq \barphi{c}  \rbrace \subseteq Q_{[D]}(c).
$$
For the reverse inclusion, it follows from the correspondence \eqref{corresp} that the face $Q_{[D]} (\barphi{c})$ is the convex hull of those vertices $v^{\sigma'}$ of $P_{[D]}$ verifying $\barphi{c} \subset \sigma'$.
But for all $\sigma' \in \Sigma_{[D]}$ we have by Definition \ref{barphi} :  
$$
c \subseteq \sigma' \Leftrightarrow \barphi{c} \subseteq \sigma',
$$
so that in fact we have 
\begin{equation}\label{QD}
Q_{[D]} (\barphi{c})=\Conv(v^{\sigma'} \mid \sigma' \in \Sigma_{[D]}(n), c \subseteq \sigma').
\end{equation}

Now take $a=(a_{\rho}) \in Q_{[D]}(c) \subset P_{[D]}$. 
Since $P_{[D]} = \Conv(v^{\sigma'} \mid \sigma' \in \Sigma_{[D]}(n))$ there exist $\sigma'_1, \ldots, \sigma'_k \in \Sigma_{[D]}(n)$ such that
$$
a = \sum_{i=1}^{k} \lambda_i v^{\sigma'_i} \quad \text{with} \ \lambda_i \in \QQ_{+}^{*} \; \text{for all} \; 1\leq i \leq k.
$$
Then for all $\rho \in c$ we have 
$0=a_\rho=\sum_{i=1}^{k} \lambda_i v^{\sigma'_i}_{\rho}$ with $v^{\sigma'_i}_{\rho} \in \NN$. Since the $\lambda_i$ are positive, it follows that for all $\rho \in c$ and all $1\leq i \leq k$, $v^{\sigma'_i}_{\rho}=0$.
By \eqref{asigma'} this implies that for $i=1, \ldots, k$ every 1-dimensional cone of $c$ is contained in $\sigma'_i$ and hence by \eqref{QD} $v^{\sigma'_i}$ is contained in the face $Q_{[D]}(\barphi{c})$. 
It follows that $a \in  Q_{[D]}(\barphi{c})$ and the lemma is proved.
\end{preuve}

\begin{cor}\label{CD0}
Let $D$ be a nef Cartier divisor and $C$ be a 1-cycle on $X_{\Sigma}$.  We have
$$
J_C \, \text{is contained in a single cone of } \, \Sigma_{[D]} \ \Rightarrow \ C \cdot D = 0, 
$$ 
and the implication is an equivalence if $[C]$ is the numerical class of an irreducible curve.
\end{cor}

\begin{preuve}
Suppose there exists a cone $\sigma \in \Sigma_{[D]}(n)$ containing all the 1-dimensional cones of $J_C$. Then by \eqref{asigma'}, the corresponding vertex $a^{\sigma}$ of the polytope $P_{[D]}$ verifies $a^{\sigma}_{\rho}=0$ for all $\rho \in J_C$. 
In particular the associated local representative $D^\sigma = \sum_{\rho \in \Sigma(1)} a^{\sigma}_{\rho} D_{\rho} \in [D]$ verifies 
$$
C \cdot D = C \cdot D^\sigma = C \cdot \sum_{\rho \in J_C} a^{\sigma}_{\rho} \, D_{\rho} =0. 
$$

Now suppose that $C \cdot D = 0$ and that $C$ is irreducible. By \cite[Prop. 1]{Payne} there exists a cone $\alpha \in \Sigma$ such that $J_C^- \subseteq \alpha(1)$ and $J_C^+ \subseteq \Adj(\alpha)$.  
Let $\sigma \in \Sigma_{[D]}(n)$ containing $\alpha$ and consider the corresponding local representative  $D^\sigma \in [D]$. We have $a^{\sigma}_{\rho}=0$ for all $\rho \in J_C^-$ and the equality 
$$
0 =  C \cdot D^\sigma =  \sum_{\rho \in J_C} a^{\sigma}_{\rho} \, C \cdot D_{\rho} 
= \sum_{\rho \in J_C^+} a^{\sigma}_{\rho} \, C \cdot D_{\rho}. 
$$
implies that $a^{\sigma}_{\rho}=0$ for all $\rho \in J_C^+$ also.
Then it follows from \eqref{asigma'} that $J_C \subset \sigma$ and the corollary is proved.
\end{preuve}

The main consequence of Lemma \ref{ouf}, from the arithmetic point of view is that a root of the homogeneous polynomial $f$ yields a rational point in $D=\V(f)$ as soon as the set $\lbrace \rho \in \Sigma(1) \mid \alpha_{\rho}=0 \rbrace$ indexing its zero coordinates is contained in a single cone of the (eventually degenerate) fan $\Sigma_{[D]}$.

\begin{prop}\label{easierpoint}
Let $K$ be a $C_1$ field satisfying hypothesis (H1) and $X_{\Sigma}$ be a simplicial projective toric variety satisfying hypotheses (H2) and (H3).
Let $D$ be a Cartier nef divisor on $X_{\Sigma}$ and let $f \in S_{[D]}$ be the homogeneous polynomial such that  
$$
D=\V(f)=\lbrace \pi(x) \in X_{\Sigma} \mid f(x)=0 \rbrace 
$$
where $\pi : \AA^{\Sigma(1)} \setminus Z(\Sigma) \to X_{\Sigma}$ is the good geometric quotient.
The following conditions are equivalent : 
\begin{enumerate}
\item[(1)] There exists a root $\alpha \in K^{\Sigma(1)}$ such that the set 
$\lbrace \rho \in \Sigma(1) \mid \alpha_{\rho}=0 \rbrace$ is contained in a single cone of $\Sigma_{[D]}$. 
\item[(2)] There exists a root $\alpha \in K^{\Sigma(1)}$ such that the set 
$\lbrace \rho \in \Sigma(1) \mid \alpha_{\rho}=0 \rbrace$ is contained in a single cone of $\Sigma$. 
\end{enumerate}
If these conditions are satisfied then $D$ has a rational point over $K$.  
\end{prop}

\begin{preuve}
The facts that the first condition implies the second and the existence of a rational point are trivial.
We have to show that if $f$ has a root $\alpha$ with all zero coordinates indexed by 1-cones of $\Sigma$ contained in a single cone $\sigma'$ of $\Sigma_{[D]}$ then it has a root $\beta$ with all zero coordinates indexed by generators of a cone $\sigma$ of $\Sigma$.

Let us write 
$$
f=\sum_{a \in \AAA_{[D]}} \lambda_a x^a \in S_{[D]}.
$$ 
and suppose $f$ has a root $\alpha$ such that the set $c=\lbrace \rho \in \Sigma(1) \mid \alpha_{\rho}=0 \rbrace$ is contained in some cone of $\Sigma_{[D]}$.    
Let us consider the polynomial $g_1$ obtained from $f$ by putting $x_{\rho}=0$ for all $\rho \in c$ :
$$
g_1= \sum_{a \in Q_{[D]}(c)} \lambda_a x^a \in K[x_{\rho}, \rho \not\in c].
$$
Since $\alpha_{\rho}=0$ for all $\rho \in c$ the root $\alpha \in K^{\Sigma(1)}$ gives a root $\alpha' \in (K^*)^{\Sigma(1)\setminus c}$ of $g_1$.

Since $c$ is contained in a single cone of $\Sigma_{[D]}$, by Definition \ref{ouf} we can put $\sigma'=\barphi{c} \in \Sigma_{[D]}$ and  consider the polynomial $g_2$ obtained from $f$ by putting $x_{\rho}=0$ for all $\rho \subseteq \sigma'$ :
$$
g_2=\sum_{a \in Q_{[D]}(\sigma')} \lambda_a x^a \in K[x_{\rho}, \rho \in J_{\sigma'}], \quad J_{\sigma'}:=\lbrace \rho \in \Sigma(1) \mid \rho \not\subseteq \sigma' \rbrace
$$
By Lemma \ref{ouf} we have $Q_{[D]}(c) = Q_{[D]}(\sigma')$ so that $g_1$ can be considered as an element of $K[x_{\rho}, \rho \in J_{\sigma'}]$ equal to $g_2$. By forgetting the $\alpha_{\rho}$ for $\rho \subseteq \sigma'$ we then get a root $\alpha'' \in (K^*)^{J_{\sigma'}}$ of $g_2$.

Now let us remark that for every $\sigma \in \Sigma$ of the same dimension as $\sigma'$ and such that $\sigma \subseteq \sigma'$, we have $\barphi{\sigma}=\sigma'$. Let us fix such a $\sigma$ and
consider the polynomial $g_3$ obtained from $f$ by putting $x_{\rho}=0$ for all $\rho \in \sigma(1)$~:
$$
g_3=\sum_{a \in Q_{[D]}(\sigma(1))} \lambda_a x^a \in K[x_{\rho}, \rho \not\in \sigma(1)].
$$
Again by Lemma \ref{ouf}, $g_3$ can be considered as an element of $K[x_{\rho}, \rho \in J_{\sigma'}]$ equal to $g_2$ and by putting $\beta'_{\rho}= \alpha''_{\rho}$ for all $\rho \in J_{\sigma'}$ and $\beta'_{\rho}=1$ for all $\rho \in \Sigma(1)$ such that $\rho \subset (\sigma' \setminus \sigma)$,
we get a root $\beta' \in (K^*)^{\Sigma(1) \setminus \sigma(1)}$ of $g_3$.
Finally by putting $\beta_{\rho}= \beta'_{\rho}$ for all $\rho \in \Sigma(1) \setminus \sigma(1)$ and $\beta_{\rho}= 0$ for all $\rho \in \sigma(1)$, we get the desired root $\beta \in K^{\Sigma(1)}$ : 
$$
f(\beta) = g_3(\beta') = g_2(\alpha'') = g_1(\alpha') =f(\alpha) = 0.
$$
The proposition is proved. 
\end{preuve}

\subsection{Global low toric degree}\label{Glob}

In this subsection we define the notion of global low toric degree, based on Lemma \ref{fund} and Proposition \ref{easierpoint}, and describe its most important properties. In particular we show that it can be lifted through toric birational contractions (Proposition \ref{open2}).

\begin{defi}\label{gltd}
Let $X_{\Sigma}$ be a complete simplicial toric variety. 
Let $D$ be a nef Cartier divisor on $X_{\Sigma}$, let $\phi : X_{\Sigma} \rightarrow X_{\bar{\Sigma}_{[D]}}$ be the toric morphism of Theorem \ref{phi}, contracting all the rational curves having zero intersection with $D$
and \mbox{$\Sigma_{[D]}$} be the corresponding (eventually degenerate) fan in $N_\RR$. 

We say that $D$ has global low toric degree if there exist a 1-cycle $C$ on $X_\Sigma$ and a subset 
$I \subseteq J_{C}^{+}$ with cardinality at least 2 such that the following conditions hold~:
\begin{equation}\label{titdeggen} \tag{A}
0< C \cdot D <  \sum_{\rho \in I} C \cdot D_{\rho}
\end{equation}
and   
\begin{equation}\label{adhoc2} \tag{B}
\forall \rho \in I, \ \exists \sigma \in \Sigma_{[D]}(n), \ 
J_{C} \setminus \lbrace \rho \rbrace \subset \sigma.
\end{equation}
\end{defi}

\begin{Remarks}
\begin{enumerate}
\item The notion of global low toric degree depends only on the numerical class of the divisors. For this reason, we will sometimes tell that the classes itself have (or not) global low toric degree.
 
\item When $C=C_\tau$ is a toric curve there is a minimal subset $I$ of cardinality 2 satisfying condition (B) : the set $\sigma_1(1)\cup \sigma_2(1) \setminus \tau(1)$ where $\sigma_1,\sigma_2 \in \Sigma (n)$ are such that $\tau = \sigma_1 \cap \sigma_2$. Indeed it is a subset of every \textit{primitive collection}  contained in $\sigma_1(1)\cup \sigma_2(1)$ (see \cite[5.1]{CLS} for a definition). 
\end{enumerate}
\end{Remarks}

\begin{theo}\label{gltd=>ratpt}
Let $K$ be a $C_1$ field satisfying hypothesis (H1) and $X_{\Sigma}$ be a complete simplicial toric variety satisfying hypotheses (H2) and (H3).
Let $D$ be an effective nef Cartier divisor on $X_{\Sigma}$ defined over $K$. 
If $D$ has global low toric degree then it has a rational point over $K$.
\end{theo}

\begin{preuve}
Let $f \in S_{[D]}$ be a homogeneous polynomial of degree $[D]$ such that 
$$
D=\V(f)=\lbrace \pi(x) \in X_{\Sigma} \mid f(x)=0 \rbrace
$$
where $\pi$ is the good geometric quotient. 
Using \eqref{titdeggen}, Lemma \ref{fund} gives us a root $\alpha \in K^{\Sigma(1)}$ such that
$$
\left\lbrace \rho \in \Sigma(1) \mid \alpha_{\rho} = 0 \right\rbrace 
\subseteq J_C \setminus \lbrace \rho \rbrace 
$$ 
for some $\rho \in I$. 
Then by use of \eqref{adhoc2}, Proposition \ref{easierpoint} allows us to construct a rational point of $D$ out of $\alpha$.  This proves the theorem. 
\end{preuve}
\medskip

\begin{prop}\label{gltdcool}
In Definition \ref{gltd} we may assume that the class $[C]$ satisfies the following additional properties :
\begin{enumerate}
\item There exist a cone $\alpha \in \Sigma$ such that $J_C^- \subseteq \alpha(1)$. In particular
the class $[C]$ is effective. 
\item The set $J_C$ is minimally linearly dependent, i.e., the space of linear relations between the minimal generators $u_\rho$ for $\rho \in J_C$ has dimension one. 
In other words we have~: 
$$
\forall [C_1] \in \NE(X_\Sigma), \  J_{C_1} \subseteq J_C \Rightarrow [C_1] \in \RR_+ [C].
$$ 
\end{enumerate}

In particular the global low toric degree condition can be checked on a finite number of 1-cycles.  
\end{prop}

In order to prove Proposition \ref{gltdcool} we need the following technical lemma.

\begin{lem}\label{JCflex}
Let $D$ be a global low toric degree divisor on $X_\Sigma$ and let $C$ be a 1-cycle and $I \subseteq J_C^+$ a subset satisfying conditions (A) and (B) of Definition \ref{gltd}.
Let $\gamma = \barphi{J_C \setminus I} \in \Sigma_{[D]}$ be the minimal cone of $\Sigma_{[D]}$ containing $J_C \setminus I$ and let $k$ be its dimension. 
For any subset $\BBB \subset \Sigma(1)$ of cardinality $k$ such that $\lbrace u_\rho \mid \rho \in \BBB \rbrace$ is linearly independent and included in $\gamma$, there exists a 1-cycle $C'$ such that $(C',I)$ satisfies conditions (A) and (B) and $J_{C'} \setminus I \subseteq \BBB$.
\end{lem}

\begin{preuve}
By induction on the cardinality of $(J_C \setminus I) \setminus \BBB$ it is sufficient to prove that for any  $\rho_0 \in (J_C \setminus I) \setminus \BBB$ there exists a 1-cycle $C'$ such that $(C',I)$ satisfies conditions (A) and (B) and $J_{C'} \setminus I \subseteq J_C \setminus \{\rho_0\} \cup \BBB$.
Let $\rho_0 \in (J_C \setminus I) \setminus \BBB$. Since $\rho_0 \in \gamma$ and the real span of $\BBB$ is equal to the one of $\gamma$, we have a relation
$$
\sum_{\rho \in \BBB \cup \{\rho_0\}} \lambda_\rho u_\rho =0, \quad \lambda_\rho \in \QQ, \ \lambda_{\rho_0} >0. 
$$
Let $C_0$ be a 1-cycle whose class is given by this relation (i.e. such that $C_0 \cdot D_\rho = \lambda_\rho$ if $\rho \in \BBB \cup \{\rho_0\}$ and $C_0 \cdot D_\rho =0$ otherwise).
By Corollary \ref{CD0}, $J_{C_0} \subseteq \gamma$ implies $C_0 \cdot D=0$. 
Similarly $C \cdot D>0$ implies that $J_C$ in not contained in a single cone of $\Sigma_{[D]}$. In particular for all $\rho \in I$ we have $J_C \not \subseteq \barphi{J_C \setminus  \{\rho\}}$ and hence $\rho \notin \barphi{J_C \setminus  \{\rho\}}$. A fortiori we have $\rho \not\subseteq \gamma$ and hence $\rho \notin J_{C_0}$, that is $C_0 \cdot D_\rho=0$.

It follows that putting $C'=\lambda_{\rho_0} C - (C \cdot D_{\rho_0}) \ C_0$ we have 
$$
0<C' \cdot D = \lambda_{\rho_0} C \cdot D 
< \lambda_{\rho_0} \sum_{\rho \in I} C \cdot D_\rho=\sum_{\rho \in I} C' \cdot D_\rho,
$$
which shows that condition (A) is satisfied by $(C',I)$.

Moreover, by construction we have $J_{C'} \subseteq J_C \cup J_{C_0}$ and
$$
C' \cdot D_{\rho_0} = 
\lambda_{\rho_0} C \cdot D_{\rho_0} - (C \cdot D_{\rho_0}) \ C_0 \cdot D_{\rho_0} =0.
$$
It follows that $J_{C'} \subseteq J_C \cup J_{C_0} \setminus  \{\rho_0\}
\subseteq J_C \setminus \{\rho_0\}  \cup \BBB$ as expected.

Finally for all $\rho \in I$ we have 
$J_{C'} \setminus \{\rho\} \subseteq J_{C} \setminus \{\rho\} \cup \BBB 
\subseteq \barphi{J_C \setminus  \{\rho\}} \cup \gamma = \barphi{J_C \setminus  \{\rho\}}$
and $(C',I)$ satisfies also condition (B). The lemma is proved.
\end{preuve}

\begin{demo}{Proof of Proposition \ref{gltdcool}}
Let $C$ be a 1-cycle and $I \subseteq J_C^+$ a subset satisfying conditions (A) and (B).
Let $\gamma = \barphi{J_C \setminus I} \in \Sigma_{[D]}$ be the smallest cone of $\Sigma_{[D]}$ containing $J_C \setminus I$ and let $k$ be its dimension. Let $\alpha \in \Sigma(k)$ be a cone of $\Sigma$ contained in $\gamma$. By Lemma \ref{JCflex} there exists a 1-cycle $C'$ such that $(C',I)$ satisfies conditions (A) and (B) and $J_{C'} \setminus I \subseteq \alpha(1)$.

In particular $[C']\in \NE(X_\Sigma) = \Nef(X_\Sigma)^\vee$ because any nef class $[D']$ contains an effective divisor whose support does not contain $V(\alpha)$ and hence verifies $[C'] \cdot [D'] \geq 0$. This proves part $(i)$. 

It remains to prove that for all $[C_1] \in \NE(X_\Sigma)$ such that $J_{C_1} \subseteq J_{C'}$  we have $[C_1] \in \RR_+ [C']$. 

Let $[C_1] \in \NE(X_\Sigma)$ be an effective class of curve such that $J_{C_1} \subseteq J_{C'}$ and let us put $\gamma_\rho=\barphi{J_{C'} \setminus \{\rho\}}$ for the smallest cone of $\Sigma_{[D]}$ containing $J_{C'} \setminus \{\rho\}$.

Let us suppose first that there exists $\rho_1 \in I$ such that $\rho_1 \notin J_{C_1}$.  

The fact that $J_{C_1} \subseteq J_{C'} \setminus \{\rho\}$ implies that the sets $J_{C_1}^+$ and $J_{C_1}^-$ are each contained in a single cone of the (possibly degenerate) fan $\Sigma_{[D]}$.
Let us denote by $\beta_{\pm}=\barphi{J_{C_1}^{\pm}}$ the smallest of these cones.
Since we have $\sum_{\rho \in J_{C_1}} C_1 \cdot D_\rho \ u_\rho =0$, we can consider the vector
$$
v = \sum_{\rho \in J_{C_1}^+} \underset{>0}{\underbrace{ C_1 \cdot D_\rho}} \ u_\rho
= \sum_{\rho \in J_{C_1}^-} \underset{>0}{\underbrace{(- C_1 \cdot D_\rho})} \ u_\rho
$$

By construction, $v$ lies in the relative interior of $\beta_{+}$ and $\beta_{-}$. 
But even if the complete fan $\Sigma_{[D]}$ is degenerate, the relative interior of its cones form a partition of $N_\RR$ (with the convention that the minimal cone is its own relative interior). 
This implies that $\beta_{+}=\beta_{-}$ and we can denote this cone by $\beta$.
Note that by definition of $\beta_{+}$ and $\beta_{-}$ we have $\beta=\barphi{J_{C_1}}$.

Now let us take $\rho \in I$. We have either $\rho \notin J_{C_1}^+$ or $\rho \notin J_{C_1}^-$, that is either $J_{C_1}^+ \subseteq J_{C'} \setminus \{\rho\}$ or $J_{C_1}^- \subseteq J_{C'} \setminus \{\rho\}$. 
In either case this implies $\beta \subseteq \gamma_\rho$. 
As we have seen in the proof of Lemma \ref{JCflex}, condition (A) implies that $\rho \not\subseteq \gamma_\rho$ because otherwise we would have $J_{C'} \subset \gamma_\rho$ and $C' \cdot D = 0$ by Corollary \ref{CD0}. It follows that $\rho \not \subseteq \beta$ and hence $\rho \notin J_{C_1}$.
Since this is true for all $\rho \in I$, we have $I \cap J_{C_1}= \emptyset$ or equivalently $J_{C_1} \subseteq J_{C'} \setminus I \subseteq \alpha(1)$, which implies $[C_1]=0$ since $\alpha$ is a simplicial cone.

This proves that any class $[C_1] \in \NE(X_\Sigma)$ such that $J_{C_1} \subseteq J_{C'}$ and $I \not \subseteq J_{C_1}$ is trivial.

Now let us consider $[C_2] \in \NE(X_\Sigma)$ with $J_{C_2} \subseteq J_{C'}$ and $I \subseteq J_{C_2}$. If we had $J_{C_2}^+ \subseteq J_{C'} \setminus I \subseteq \alpha(1)$ it would imply $[-C_2] \in \NE(X_\Sigma)$ by the argument we used to prove that $[C'] \in \NE(X_\Sigma)$. 
Since we cannot have $[C_2] \in \NE(X_\Sigma)$ and $[-C_2] \in \NE(X_\Sigma)$ it follows that there exists $\rho_2 \in J_{C_2}^+ \cap I$. 
Let us put $\lambda_2 = \dfrac{C_2 \cdot D_{\rho_2}}{C' \cdot D_{\rho_2}} \in \QQ_+^*$ and $C_1=C_2-\lambda_2 C'$. By construction we have $J_{C_1} \subseteq J_{C'} \cup J_{C_2} = J_{C'}$ and
$$
C_1 \cdot D_{\rho_2}=C_2 \cdot D_{\rho_2} - \lambda_2 \ C' \cdot D_{\rho_2} = 0.
$$ 
It follows that $\rho_2 \notin J_{C_1}$ and in particular $I \not\subseteq J_{C_1}$, which implies $[C_1]=0$ as above, that is $[C_2]= \lambda_2 [C']$.
This proves part $(ii)$.

The final assertion follows from the finiteness of the set of minimal linear relations and the invariance of conditions (A) and (B) by rescaling of $C$.
\end{demo}

One of the important properties of the global low toric degree condition is its behavior under a toric birational contraction over the variety $X_{\bar{\Sigma}_{[D]}}$ defined out of $[D]$. Let us first precise this notion.  

\begin{defi}\label{over}
Let $X_{\Sigma}$ be a complete simplicial toric variety and $D$ be a nef Cartier divisor on $X_\Sigma$. 
Let $\phi : X_{\Sigma} \rightarrow X_{\bar{\Sigma}_{[D]}}$ be the toric morphism of Theorem \ref{phi}, contracting all the rational curves having zero intersection with $D$ and let
$\bar{D}$ be the ample divisor on $X_{\bar{\Sigma}_{[D]}}$ such that $[D]=\phi^*[\bar{D}]$.

A toric rational map $\xi : X_{\Sigma} \dashrightarrow  X_{\Sigma'}$ is a \textit{map over $X_{\bar{\Sigma}_{[D]}}$} if there exists a proper toric morphism $\phi' : X_{\Sigma'} \rightarrow  X_{\bar{\Sigma}}$ making the following diagram commutative~:
$$
\xymatrix{
X_{\Sigma}\ar[rd]_{\phi} \ar@{-->}[rr]^{\xi}&  & X_{\Sigma'}\ar[ld]^{\phi'} \\
& X_{\bar{\Sigma}} & }
$$
\end{defi}

\begin{lem}\label{SigmabarD'}
Let $X_{\Sigma}$, $D$, $\phi : X_{\Sigma} \rightarrow X_{\bar{\Sigma}_{[D]}}$, and $\bar{D}$ be as in Definition \ref{over}. 
Let additionally $\xi : X_{\Sigma} \dashrightarrow X_{\Sigma'}$ be a toric map
over $X_{\Sigma_{[D]}}$, and $\phi' : X_{\Sigma'} \rightarrow  X_{\Sigma_{[D]}}$ be the corresponding proper morphism.

If $D'$ is a divisor on $X_{\Sigma'}$ such that $[D']=\phi'^*[\bar{D}]$ then \mbox{$\phi' : X_{\Sigma'} \rightarrow X_{\Sigma_{[D']}}$} is the toric morphism contracting all the curves having zero intersection with $D'$ and we have  
$$
X_{\bar{\Sigma}_{[D']}} =X_{\bar{\Sigma}_{[D]}}.
$$
Moreover, if $\xi$ is birational, then the result holds for all $D'$ such that $[D']=\xi_*[D]$.
\end{lem}

\begin{preuve}
If $[D']=\phi'^*[\bar{D}]$, by the projection formula we have for all curve $C'$ on $X_{\Sigma'}$  
$$
C' \cdot D'=0 \Leftrightarrow  \phi'_*C' \cdot \bar{D}=0 \Leftrightarrow  \phi'_*C'=0 
$$
because $\bar{D}$ is ample.
It follows that $\phi'$ is the morphism of Theorem \ref{phi} and $X_{\bar{\Sigma}_{[D']}} =X_{\bar{\Sigma}_{[D]}}$.

For the final assertion it suffices to remark that if $\xi$ is birational then the composition $\xi_* \circ \xi^*$ is the identity of $N^1(X_{\Sigma'})$. 
Indeed, since $\phi = \phi' \circ \xi$, this gives 
$$
\xi_*[D] \, =\, \xi_*\phi^*[\bar{D}] \, =\, \xi_*\xi^*\phi'^*[\bar{D}] \, =\, \phi'^*[\bar{D}]
$$
and the result follows.
\end{preuve}

\begin{prop}\label{open2}
Let $X_{\Sigma}$ be a complete simplicial toric variety and $D$ be a nef Cartier divisor on $X_\Sigma$. 
Let $\xi : X_{\Sigma} \dashrightarrow X_{\Sigma'}$ be a toric birational contraction over $X_{\bar{\Sigma}_{[D]}}$ such that $X_{\Sigma'}$ is simplicial.
If $\xi_*[D]$ has global low toric degree, then so has $[D]$. 
\end{prop}

\begin{preuve}
We may assume that $\Sigma$ and $\Sigma'$ belong to the same real vector space $N_\RR$ so that every 1-dimensional cone $\rho \in \Sigma'(1)$ can be considered as an element of $\Sigma(1)$.
The associated toric divisors of $X_{\Sigma'}$ are denoted $D'_{\rho}$ while those of $X_{\Sigma}$ are denoted $D_{\rho}$.
  
By hypothesis the conditions of Definition \ref{gltd} are satisfied by any divisor $D'$ on $X_{\Sigma'}$ such that $[D']=\xi_*[D]$~: there exists a 1-cycle $C'$ on $X_{\Sigma'}$ and a subset $I \subseteq J_{C'}^{+}$ with cardinality at least 2 such that 
$$
\forall \rho \in I, \ \exists \bar{\sigma} \in \Sigma_{[D']}, \ 
J_{C'} \setminus \lbrace \rho \rbrace \subseteq \bar{\sigma}
$$
and
$$
0< C' \cdot D' <  \sum_{\rho \in I} C' \cdot D'_{\rho}.
$$

But since $\xi$ is birational over $X_{\Sigma_{[D]}}$, it follows from Lemma \ref{SigmabarD'} that $\Sigma_{[D']}=\Sigma_{[D]}$. 
Then we only have to remark that by definition of the numerical pullback of curves, if $C$ is any 1-cycle on $X_{\Sigma}$ such that $[C]=\xi_{num}^{*}[C']$, then we have $J_C=J_{C'}$, $[C] \cdot [D] = [C'] \cdot \xi_*[D] \:$ and $ \: \sum_{\rho \in I} C \cdot D_{\rho} =\sum_{\rho \in I} C' \cdot D'_{\rho}$ so that conditions of Definition \ref{gltd} are also satisfied by $D$. The proposition is proved.
\end{preuve}

\subsection{Global vs. Restricted low toric degree}\label{GlobDir}

It follows from Proposition \ref{open2} that in Example \ref{L3}, the pullback of the conic has global low toric degree while it has not restricted low toric degree. 
In fact, by definition, the restricted low toric degree condition can be viewed as a particular case of the global low toric degree one~:

\begin{lem}\label{rgltd}
Let $X_\Sigma$ be a normal toric variety. 
A Cartier divisor $D$ on $X_\Sigma$ has restricted low toric degree if and only if there exists a complete simplicial toric subvariety $W \subseteq X_\Sigma$ of Picard number one, such that the restriction $D_{|W}$ has global low toric degree.
\end{lem}  

It is nevertheless important to observe that a divisor may have restricted low toric degree without having global low toric degree. 
Indeed the intersection numbers considered in Definition \ref{gltd} are computed in the ambient variety $X_\Sigma$ and can differ from the intersection numbers computed in a proper subvariety if $X_\Sigma$ is singular, as noticed in paragraph \ref{InterSubvar}. This is where the name global low toric degree comes from. Let us give an example. 

\begin{ex}
Let $\Sigma$ be the 2-dimensional fan with the following minimal generators :
$$
u_1=\vv{1}{0}, u_2=\vv{1}{2}, u_3=\vv{-1}{2}, u_4=\vv{0}{-1}
$$
Here is a drawing of $\Sigma$ together with the intersection matrix of its 4 toric divisors :

\hfill
\mbox{
\begin{minipage}[c]{6cm}
\includegraphics[scale=0.8]{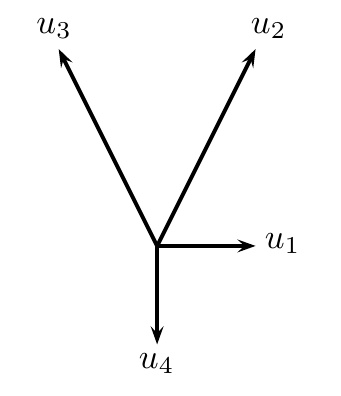} 
\end{minipage}
\begin{minipage}[c]{7cm}
\vspace*{1cm}
\renewcommand{\arraystretch}{1.5}
$ \begin {array}{c|ccccc}  
\pmb{\cdot} & D_1 & D_{2} & D_3 & D_{4}  \\
\hline 
D_{1} \quad & -\frac{1}{2} & \frac{1}{2}  &   0         &    1  \\ 
D_{2} \quad & \frac{1}{2}  & -\frac{1}{4} & \frac{1}{4} &    0  \\ 
D_{3} \quad &     0        &  \frac{1}{4} & \frac{1}{4} &    1  \\ 
D_{4} \quad &     1        &      0       &      1      &    2 \\
\end {array} $  
\vspace*{1.0cm}
\end{minipage}
} 

\vspace{-0.3cm}
Let $D$ be a divisor linearly equivalent to $4 D_3 + 2 D_4$. We have 
$$
D_1 \cdot D = 2,\quad  D_2 \cdot D = 1, \quad  D_3 \cdot D = 3 \quad \text{and} \quad D_4 \cdot D = 8,
$$
which shows that $D$ is an ample divisor with restricted low toric degree since $D_2\simeq \PP^1$ and the degree of the restriction $D_{|D_2}$ is 1.
But according to Proposition \ref{tresbeau} below, $D$ has not global low toric degree since the extremal classes of $X_\Sigma$ are $[D_1]$ and $[D_2]$ and we have
$$
D_1 \cdot D = 2 > \frac{3}{2} = D_1 \cdot (D_2 + D_4) \quad \text{and} \quad
 D_2 \cdot D = 1 > \frac{3}{4} = D_2 \cdot (D_1 + D_4). 
$$ 
\end{ex}

Of course, when $X_\Sigma$ is smooth, the intersection numbers in the ambient and in its subvarieties coincide such that having restricted low toric degree implies having global low toric degree, the converse being false in general.

The following result shows that for ample divisors, global low toric degree implies restricted low toric degree.

\begin{prop}\label{tresbeau}
Let $X_\Sigma$ be a simplicial projective toric variety and $D$ be an ample divisor on $X_{\Sigma}$.  
The following conditions are equivalent :
\begin{enumerate}
\item[(1)] There exist an effective 1-cycle $C$ and a subset $I \subseteq \Sigma(1)$ such that
$$
C \cdot D < \sum_{\rho \in I} C \cdot D_{\rho}. 
$$
\item[(2)] There exists an extremal class $[C]$ such that 
$$
C \cdot D < \sum_{\rho \in J_{C}^+} C \cdot D_{\rho}.
$$
\end{enumerate}
\end{prop}

\begin{preuve}
The implication (2) $\Rightarrow$ (1) is trivial. 
Let us show that (1) implies (2).
    
Let $C$ be an effective 1-cycle and $I \subseteq \Sigma(1)$ a subset such that (1) holds. 
Since $X_{\Sigma}$ is projective there exist extremal classes $[C_1], \ldots, [C_m]$ such that 
$$
[C]= \sum_{i=1}^{l} \lambda_i [C_i] 
\; \text{with} \; \lambda_1, \ldots, \lambda_l \in \QQ_{+}^{*}.
$$
For all $1 \leq i \leq l$, we have by definition of $J_{C_i}^{+}$ : 
$$
\sum_{\rho \in J_{C_i}^{+}}  C_i \cdot D_{\rho} \geq
\sum_{\rho \in J_{C_i}^{+} \cap I}  C_i \cdot D_{\rho} \geq
\sum_{\rho \in  I}  C_i \cdot D_{\rho}.
$$
It follows that
$$
\sum_{i=1}^{l} \lambda_i \, C_i \cdot D = C \cdot D  < 
\sum_{\rho \in I} C \cdot D_{\rho} = 
\sum_{i=1}^{l} \lambda_i \sum_{\rho \in I}  C_i  \cdot D_{\rho} \leq 
\sum_{i=1}^{l} \lambda_i \sum_{\rho \in J_{C_i}^{+}}  C_i \cdot D_{\rho},
$$ 
so that we must have 
$$
C_i \cdot D \; < \sum_{\rho \in J_{C_i}^{+}}  C_i \cdot D_{\rho}
$$
for at least one $i \in \lbrace 1, \ldots, l \rbrace$.
This shows that (2) holds and the proposition is proved.    
\end{preuve}

The very complementarity of the notions of restricted and global low toric degree comes from their behavior under toric maps (Propostions \ref{dltdbar} and \ref{open2}) and appears clearly in Theorem \ref{gral}.

\section{Low toric degree~: general definition}\label{Def}

This section is devoted to state and prove the main theorem of the article, on which is based our definition of a low toric degree.

\begin{theo}\label{gral}
Let $X_{\Sigma}$ be a simplicial projective toric variety.
Let $D$ be a nef Cartier divisor on $X_{\Sigma}$.
Let $\phi : X_{\Sigma} \rightarrow X_{\bar{\Sigma}_{[D]}}$ be the toric morphism of Theorem \ref{phi}, contracting all the rational curves having zero intersection with $D$
and let $\bar{D}$ be an ample divisor on $X_{\bar{\Sigma}_{[D]}}$ such that $D$ is linearly equivalent to the pullback $\phi^*\bar{D}$.

The following conditions are equivalent~:
\begin{enumerate}
\item[(i)] The ample divisor $\bar{D}$ has restricted low toric degree. 

\item[(ii)] There exist a subvariety $V \subset X_{\Sigma}$  
such that the restriction $D_{|V}$ has global low toric degree.
 
\item[(iii)] There exists a desingularization $\eta : X_{\tilde{\Sigma}} \rightarrow X_{\Sigma}$ such that the pullback $\tilde{D}=\eta^* D$ has global low toric degree.
\end{enumerate}
\end{theo}

\begin{defi}\label{ltd}
Let $X_{\Sigma}$ be a simplicial projective toric variety. 
A Weil divisor $D$ on $X_\Sigma$ is said to have \textit{low toric degree} if it is not nef and Cartier or if it is nef and Cartier and satisfies the conditions of Theorem \ref{gral}. 
\end{defi}

\begin{Remarks} 
\begin{enumerate}
\item It is easily seen that if $X_\Sigma$ is a weighted projective space then the divisor $D$ has low toric degree if and only if it has low weighted degree in the sense of Theorem \ref{Kollar}. 
In particular for $X_\Sigma=\PP_K^n$, $D$ has low toric degree if and only if its degree $d$ verifies $d \leq n$ as in Definition \ref{defC1}.
\item For a toric variety $X_{\Sigma}$ simplicial and complete but not projective, we can take condition $(i)$ of Theorem \ref{gral} as a definition for low toric degree. In this case the implications $(i) \Rightarrow (ii) \Rightarrow (iii)$ still hold but the proof of $(iii) \Rightarrow (i)$ we give below relies on the toric Minimal Model Program which is known only in the projective case.
\end{enumerate}
\end{Remarks}
\bigskip

Before proving Theorem \ref{gral}, let us show that the low toric degree condition ensures as expected the existence of rational points over $C_1$ fields.

\begin{cor}\label{ratpt}
Let $K$ be a $C_1$ field satisfying hypothesis (H1) and 
$X_{\Sigma}$ be a simplicial projective toric variety satisfying hypothesis (H2). 

Every low toric degree divisor on $X_{\Sigma}$ admits a rational point over $K$.
\end{cor}

\begin{preuve}
Let $D$ be an effective divisor on $X_{\Sigma}$. If $D$ is not Cartier and nef then it has a $K$-rational point by Proposition \ref{trivial}. If it is Cartier and nef then by hypothesis condition (iii) of Theorem \ref{gral} is satisfied. 
In particular the pullback $\tilde{D}$ of $D$ by $\eta$ admits a rational point $p \in \tilde{D}$ by Theorem \ref{gltd=>ratpt}. 
Then since the morphism $\eta$ is toric and every toric map is defined over $K$, the image of $p$ by $\eta$ is a rational point $\eta(p)$ in $D$. 
This proves the corollary.   
\end{preuve}

The proof of the equivalences of Theorem \ref{gral} goes circularly according to the following scheme~:

$$
\xymatrix{
& &  (ii)  \ar@{=>}[lld]_{
\substack{\text{Lem. \ref{Wlift}}\\ \text{Prop. \ref{open2}} \quad}
} & &  \\
(iii)\ar@{=>}[rrrr]_{
\substack{\text{Cor. \ref{gltd=>dltd}, Lem. \ref{SigmabarD'}, Prop. \ref{dltdbar}} \\ \text{Prop. \ref{descent}} }
} & & & &
\ar@{=>}[llu]_{
\substack{\text{Lem. \ref{rgltd}, \ref{Wlift}}\\ \text{Prop. \ref{open2}}}
} (i) 
}
$$

Let us give some details.
\medskip

\begin{demo}{Sketch of the proof of Theorem \ref{gral}}\label{skpr}
Once we have remarked that by Lemma \ref{rgltd} condition $(i)$ is equivalent to the existence of a generalized weighted projective space $\bar{W} \subseteq X_{\bar{\Sigma}_{[D]}}$ such that the restriction $\bar{D}_{|\bar{W}}$ has global low toric degree, the implication $(i) \Rightarrow (ii)$ amounts to lift this global low toric degree condition to the restriction  $D_{|V}$ of $D$ to a well chosen subvariety $V$ of $X_\Sigma$.

The lifting is done in two steps~: Lemma \ref{Wlift} proves the existence of a birational toric morphism between such a $V \subseteq X_\Sigma$ and $\bar{W} \subseteq X_{\bar{\Sigma}_{[D]}}$, and then Proposition \ref{open2} gives the result.
\smallskip

$$
\xymatrix{
X_{\tilde{\Sigma}} \ar[d]_{\eta} & \ar@{_{(}->}[l]  \ \tilde{V} \ar[d]^{\eta_{|\tilde{V}}}  \\
X_{\Sigma}  \ar[d]_{\phi}  & \ar@{_{(}->}[l]  \ V  \ar[d]^{\phi_{|V}} \\
X_{\bar{\Sigma}_{[D]}}& \ar@{_{(}->}[l]  \ \bar{W}
}
$$

The implication $(ii) \Rightarrow (iii)$ is similar. 
By Proposition \ref{open2} we lift the global low toric degree condition 
of $D_{|V}$ to the restriction $\tilde{D}_{|\tilde{V}}$ of the pullback $\tilde{D}=\phi^*D$ to a subvariety $\tilde{V} \subseteq X_{\tilde{\Sigma}}$ given by Lemma \ref{Wlift}.
Then the smoothness of $X_{\tilde{\Sigma}}$ implies that $\tilde{D}$ itself has global low toric degree.
\smallskip

Finally $(iii) \Rightarrow (i)$ is a direct application of the toric MMP. 
Since restricted low toric degrees can be descended by contracting toric morphisms by Proposition \ref{dltdbar}, it is sufficient to prove that some pullback of $\bar{D}$ by such a morphism has restricted low toric degree to get $(i)$.
The idea is to mimic the proof of Proposition \ref{tresbeau} to prove that the global low toric degree of $\tilde{D}$ implies the restricted low toric degree of its image by some birational contraction.  
Namely, starting with the inequalities 
\begin{equation*} \tag{A}
0< C \cdot \tilde{D} <  \sum_{\rho \in I} C \cdot \tilde{D}_{\rho}
\end{equation*}  
we consider a decomposition $[C]= \sum_{i=1}^{l} \lambda_i [C'_{i}]$ of $[C]$ as a positive linear combination of extremal classes and show that we can deduce inequalities similar to (A) with some $C'_i$ instead of $C$. 
This is in general not possible in $X_{\tilde{\Sigma}}$ since some of the $C'_i$'s may fail to have positive intersection with $\tilde{D}$. 
The solution is therefore to get rid of the superfluous $C'_i$'s either by contracting them or by flipping them, which can be realized by running a logarithmic MMP with boundary divisor $B=\sum_{\rho \in \tilde{\Sigma} \setminus I} \tilde{D}_{\rho}$ relatively to the variety $X_{\bar{\Sigma}_{[D]}}$ (Proposition \ref{descent}). 
This yields a birational contraction $\xi : X_{\tilde{\Sigma}} \dashrightarrow X_{\Sigma'}$ over $X_{\bar{\Sigma}_{[D]}}$ such that the divisor $D'=\xi_*\tilde{D}$ has restricted low toric degree. Then it follows from Lemma \ref{SigmabarD'} that $D'$ is naturally a pullback of $\bar{D}$ and we get the result.  
\end{demo}
\bigskip

\begin{lem}\label{Wlift}
Let $\phi : X_{\Sigma_1} \rightarrow X_{\Sigma_2}$ be a contracting toric morphism between complete toric varieties.  
For all toric subvariety $W_2 \subseteq X_{\Sigma_2}$ there exists a toric subvariety $W_1 \subseteq X_{\Sigma_1}$ such that the restriction $\phi_{|W_1} : W_1 \rightarrow W_2$ is a birational toric morphism.
\end{lem}

\begin{preuve}
Let us denote by $N_1$ and $N_2$ the respective lattices of one-parameter subgroups of the dense tori
$T_1 \subset X_{\Sigma_1}$ and $T_2 \subset X_{\Sigma_2}$, and by $N_0$ the kernel of the map of lattices $\bar{\phi} : N_1 \to N_2$ associated to $\phi$. 
We have $N_2 \simeq N_1/N_0$.

Let $\alpha_2$ be the cone of $\Sigma_2$ such that $W_2=\V(\alpha_2)$. 
The preimage $\bar{\phi}^{-1}(\alpha_2)$ is a union of cones of $\Sigma_1$ with dimension $\dim(\alpha_2)+\rk(N_0)$. In particular there exists a cone $\alpha_1 \in \Sigma_1$ such that $\bar{\phi}(\alpha_1) \subseteq \alpha_2$ and the dimension of $W_1=\V(\alpha_1)$ is $\dim(W_2)$.
 
For such a cone $\alpha_1$ we have $\Span(\alpha_2)=\bar{\phi}(\Span(\alpha_1))$ and hence $(N_2)_{\alpha_2}=\bar{\phi}((N_1)_{\alpha_1})$. Since $\bar{\phi}$ is surjective from $N_1$ to $N_2$ it follows that the map of lattices $\bar{\phi_{|W_1}} : N_1(\alpha_1) \to N_2(\alpha_2)$ associated to $\phi_{|W_1}$ is surjective. Since these lattices have the same rank $\dim(W_1)=\dim(W_2)$, the surjection is in fact a bijection, which implies that $\phi_{|W_1}$ is birational and proves the lemma.
\end{preuve}

\begin{prop}\label{descent}
Let $X_{\Sigma_1}$ be a simplicial projective toric variety. Let $D_1$ be a nef Cartier divisor on $X_{\Sigma_1}$, $\phi_1 : X_{\Sigma_1} \rightarrow X_{\bar{\Sigma}_{[D_1]}}$ be the toric morphism of Theorem \ref{phi} and \mbox{$\Sigma_{[D_1]}$} be the corresponding (possibly degenerate) fan in $N_\RR$. 

Suppose that there exist a 1-cycle $C_1$ on $X_{\Sigma_1}$ and a subset 
$I_1 \subset J_{C_1}^{+}$ such that the following conditions are satisfied~:  
\begin{enumerate}
\item[(a)] \begin{equation*}
0< C_1 \cdot D_1 <  \sum_{\rho \in I_1} C_1 \cdot D_{\rho},
\end{equation*}

\item[(b)] \begin{equation*}
\exists \sigma \in \Sigma_{[D_1]}, \ J_{C}^- \subset \sigma \quad 
\text{and} \quad
\forall \rho \in I_1, \ \rho \not\subset \sigma.
\end{equation*}
\end{enumerate}

Then there exist a toric birational contraction $\xi : X_{\Sigma_1} \dashrightarrow X_{\Sigma_2}$ 
over $X_{\bar{\Sigma}_{[D_1]}}$ such that $\xi_*D_1$ has restricted low toric degree.
\end{prop}

\begin{preuve}
It is essentially an application of the relative logarithmic MMP for simplicial projective toric varieties. Let us consider the logarithmic MMP with boundary divisor $B=\sum_{\rho \in \Sigma_1(1) \setminus I} D_{\rho}$ relative to the variety $X_{\bar{\Sigma}_{[D_1]}}$. 
To begin with, we show that conditions (a) and (b) are invariant through each step of the MMP. For this 
let us describe the full algorithm, keeping track of what happens to $[D_1]$, $[C_1]$ and $I_1$.
\medskip

We start with $X_\Sigma=X_{\Sigma_1}$, $[B]=[B_1]$, $\phi=\phi_1$, $[D]=[D_1]$, $[C]=[C_1]$ and  $I=I_1$ and do the following steps~:
\begin{enumerate}
\item[(1)] If $K_{X_\Sigma} + B$ is $\phi$-nef, then stop.
\item[(2)] If $K_{X_\Sigma} + B$  is not $\phi$-nef, then there is an extremal ray $R$ in the relative Mori cone $\NE(X_\Sigma/X_{\bar{\Sigma}_{[D_1]}})$ with $(K_{X_\Sigma} + B) \cdot R < 0$.
Let $\psi : X_\Sigma \rightarrow X_{\Sigma'}$ be the corresponding extremal contraction.
\item[(3)] If $\psi$ is fibering, then stop.
\item[(4)] If $\psi$ is divisorial, then replace $(X_\Sigma, B, \phi, [D], [C], I)$ with 
$(X_{\Sigma'}, \psi_*B, \phi', \psi_*[D], \psi_*[C], I')$ where $\phi' : X_{\Sigma'} \rightarrow X_{\bar{\Sigma}_{[D_1]}}$ is the unique morphism satisfying $\phi = \phi' \circ \psi$ and $I'=I \cap \Sigma'(1)$. Return to step (1) and continue.
\item[(5)] If $\psi$ is flipping and $\nu : X_\Sigma \dashrightarrow X_\Sigma^+$ is the corresponding flip, replace $(X_\Sigma, B, \phi, [D], [C], I)$ with $(X_\Sigma^+, \nu_*B, \phi^+, \nu_*[D], \nu_*[C], I)$ where $\phi^+ : X_{\Sigma}^+ \rightarrow X_{\bar{\Sigma}_{[D_1]}}$ is the unique morphism satisfying $\phi = \phi^+ \circ \nu$.
Return to step (1) and continue.
\end{enumerate}

Notice that every toric map that appears during the process is a map over $X_{\bar{\Sigma}_{[D_1]}}$, i.e., at each step we have a toric morphism $\phi : X_{\Sigma} \rightarrow X_{\bar{\Sigma}_{[D]}}=X_{\bar{\Sigma}_{[D_1]}}$ contracting the curves that have zero intersection with $D$. Moreover since these maps are birational we always have $\Sigma_{[D]}=\Sigma_{[D_1]}$ in $N_\RR$, therefore what we show is that conditions (a) and (b) are invariant \textit{for a fixed} $\sigma \in \Sigma_{[D_1]}$.
\medskip

Let us begin by the flipping case. The flip $\nu : X_\Sigma \dashrightarrow X_\Sigma^+$ is an isomorphism in codimension 1 and we have $\Sigma(1) =\Sigma^+(1)$. If we denote by $D_\rho$ and $D^+_\rho$ the toric divisor corresponding to $\rho \in \Sigma(1)$ respectively in $X_\Sigma$ and $X_\Sigma^+$, then the pushforward of divisor classes is the isomorphism  $\nu_* : N^1(X_\Sigma) \rightarrow N^1(X_\Sigma^+)$ induced by $\nu_*D_\rho=D^+_\rho$. 
By duality, the numerical pullback of curves is the isomorphism  $\nu_{num}^* : N_1(X_{\Sigma^+}) \rightarrow N_1(X_\Sigma)$ defined by
$$
\nu_{num}^*[C']\cdot [D_\rho] = [C']\cdot [D^+_\rho] \ \text{for all} \ \rho \in \Sigma(1).
$$
Its inverse isomorphism is the pushforward of curve classes so that we also have 
$$
\nu_*[C]\cdot [D^+_\rho] = [C] \cdot [D_\rho] \ \text{for all} \ \rho \in \Sigma(1).
$$
It follows that conditions (a) and (b) are satisfied by $(D,C,I)$ if and only if they are satisfied by $(\nu_*D,\nu_*C,I)$.
\medskip

Then let us turn to the divisorial case.
The morphism $\psi : X_\Sigma \rightarrow X_{\Sigma'}$ is surjective and we have $\Sigma'(1) = \Sigma(1)\setminus \{\rho_0 \}$ where $D_{\rho_0}$ is the exceptional divisor in $X_\Sigma$. If we denote by $D_\rho$ and $D'_\rho$ the toric divisor corresponding to $\rho \in \Sigma'(1)$ respectively in $X_\Sigma$ and $X_{\Sigma'}$, then the pushforward of divisor classes is the epimorphism  $\psi_* : N^1(X_\Sigma) \rightarrow N^1(X_\Sigma')$ induced by $\psi_*D_\rho=D'_\rho$. 
By duality, the numerical pullback of curves is the monomorphism  $\psi_{num}^* : N_1(X_\Sigma') \rightarrow N_1(X_\Sigma)$ defined for $[C'] \in N_1(X_\Sigma')$ by
$$
\psi_{num}^*[C']\cdot [D_\rho] = \begin{cases} 
[C']\cdot [D'_\rho] & \text{if} \ \rho \in \Sigma'(1) \\
0  & \text{if} \quad \rho = \rho_0.
\end{cases}
$$
If $[C_0]$ is a generator of the extremal ray contracted by $\psi$, then the pushforward of curve classes is the unique epimorphism  (with abuse of notation) $\psi_* : N_1(X_\Sigma) \rightarrow N_1(X_\Sigma')$ sending $[C_0]$ to 0 and admitting $\psi_{num}^*$ as a section. 
In particular there exists a coefficient $\lambda_{[C]} \in \RR$ such that 
$$
\psi_{num}^*\psi_*[C] = [C] + \lambda_{[C]} [C_0].
$$

Since $[D]$ is the pullback of an ample class $[\bar{D}]$ on $X_{\bar{\Sigma}_{[D_1]}}$ and $\psi :X_\Sigma \rightarrow X_{\Sigma'}$ is a toric morphism over $X_{\bar{\Sigma}_{[D_1]}}$ we also have $\psi_*[D]=\phi'^*[\bar{D}]$ and hence $[D]=\psi^*\psi_*[D]$.  
In particular, by the projection formula we have 
\begin{equation}\label{projform}
\psi_*[C] \cdot \psi_*[D] = [C] \cdot \psi^*\psi_*[D] = [C] \cdot [D]. 
\end{equation}

Now let $V'=\V(\alpha'), \alpha' \in \Sigma'$, be the image of $D_{\rho_0}$ by $\psi$. 
We may assume that the class $[C_0]$ is given by the relation  
$$
-u_{\rho_0} + \sum_{\rho \in \alpha'(1)} c_\rho u_\rho, \quad 
\text{with} \ c_\rho \in \QQ_+^* \ \text{for all} \ \rho \in \alpha'(1),
$$
so that $\psi_{num}^* \psi_*[C] \cdot [D_{\rho_0}] = 0$ implies $\lambda_{[C]}=[C] \cdot [D_{\rho_0}]$.

Note that in this situation, by the definition of $\phi$-hull (Definition \ref{barphi}) we have 
\begin{equation}\label{JC0}
\barphi{\{\rho_0 \}} = \barphi{\alpha'(1)} = \barphi{J_{C_0}}.
\end{equation}

To simplify notations let us put $[C']=\psi_*[C]$ and  $[C] \cdot [D_{\rho}]= b_\rho$ for all $\rho \in J_C$ (in particular $\lambda_{[C]}=b_{\rho_0}$).
Since $\psi_{num}^*[C'] = [C] + b_{\rho_0} [C_0]$, for $\rho \in \Sigma'(1)$ we have \begin{equation}\label{xiCD}
[C'] \cdot [D'_{\rho}] =  \psi_{num}^*[C'] \cdot [D_{\rho}] = 
\begin{cases} 
b_\rho & \text{if} \quad  \rho \in J_C \setminus J_{C_0} \\
b_\rho + b_{\rho_0} c_\rho & \text{if} \quad \rho \in J_C \cap \alpha'(1)  \\
b_{\rho_0} c_\rho & \text{if} \quad \rho \in \alpha'(1) \setminus J_{C} \\ 
0  & \text{if} \quad \rho = \rho_0 \ \text{or} \ \rho \notin  J_C \cup J_{C_0} .
\end{cases}
\end{equation}

In particular we have $J_{C'} \subset J_C \setminus \{\rho_0 \} \cup \alpha'(1)$. 

We then distinguish three cases~:
\begin{enumerate}
\item[$1^{st}$ case] If $\rho_0 \notin J_C$ then $[C'] \cdot [D'_{\rho}] =  [C] \cdot [D_{\rho}]$ for all $\rho \in \Sigma'(1)$ and conditions  (a) and (b) are satisfied by $([D],[C],I)$ if and only if they are satisfied by $(\psi_*[D],[C'],I)$.

\item[$2^{nd}$ case] If $\rho_0 \in J_C^-$ we claim that $[C_0] \cdot \sum_{\rho \in I} [D_{\rho}] = 0$, or equivalently that $[C_0] \cdot (K_{X_\Sigma} + B) = 0$, so that this case never occurs.

Indeed, it follows from \eqref{JC0} and condition (b) that
$$
\barphi{J_{C_0}} \subset \barphi{J_C^-} \subset \bar{\Sigma}_{[D_1]}. 
$$
and for all $\rho \in I$, $\barphi{\rho} \notin \bar{\Sigma}_{[D_1]}$. 
In particular we have $I \cap J_{C_0} = \emptyset$ and $[C_0] \cdot \sum_{\rho \in I} [D_{\rho}] = 0$ as claimed.
 
\item[$3^{rd}$ case] If $\rho_0 \in J_C^+$ then it follows from \eqref{xiCD} that for all $\rho \in J_{C'}$ we have $[C'] \cdot [D'_{\rho}] \geq [C] \cdot [D_{\rho}]$.
In particular $J_{C'}^- \subset J_C^- \subset \sigma$ and since $I'=I \cap \Sigma'(1) \subset I$ we have for all $\rho \in I'$, $\rho \not\subset \sigma$. 
This shows that condition (b) is satisfied by $C'$ and $I'$.

For condition (a) we may assume that $\rho_0 \in I$ because otherwise we have $I'=I$ and condition (a) for $C'$ and $I'$ is the same as for $C$ and $I$. 
Recall that with our notations we have 
$$
[C_0] \cdot  [D_{\rho}] = 
\begin{cases} 
c_\rho & \text{if} \quad  \rho \in \alpha'(1) \\
-1 & \text{if} \quad \rho = \rho_0 \\ 
0  & \text{else}.
\end{cases}
$$
The condition $[C_0] \cdot [K_{X_\Sigma} + B] < 0$ coming from step (2) is equivalent to the condition $[C_0] \cdot \sum_{\rho \in I} [D_{\rho}] > 0$,  that is 
$$
\sum_{\rho \in \alpha'(1) \cap I} c_\rho > 1.
$$
Since $[C'] \cdot [D'_{\rho_0}]=0$ and 
$I \setminus J_{C_0} \subset J_C \setminus J_{C_0}$, we have by \eqref{xiCD}~:
$$
\begin{aligned}
\sum_{\rho \in I} [C'] \cdot [D'_{\rho}] 
& =   \sum_{\rho \in \alpha'(1) \cap  I} [C'] \cdot [D'_{\rho}]
     + \sum_{\rho \in I \setminus J_{C_0} } [C'] \cdot [D'_{\rho}]   \\
& =   \sum_{\rho \in \alpha'(1) \cap  I} (b_\rho + b_{\rho_0} c_\rho)
     + \sum_{\rho \in I \setminus J_{C_0} } b_\rho                   \\  
& =   \sum_{\rho \in I \setminus \{\rho_0 \}} b_\rho  
     + \sum_{\rho \in \alpha'(1) \cap  I } b_{\rho_0} c_\rho         \\   
& >   \sum_{\rho \in I \setminus \{\rho_0 \}} b_\rho +  b_{\rho_0} 
= \sum_{\rho \in I} b_\rho = \sum_{\rho \in I} [C] \cdot [D_{\rho}].
\end{aligned}
$$
Now with \eqref{projform} this gives 
$$
0< [C'] \cdot \psi_*[D] <  \sum_{\rho \in I} [C'] \cdot [D'_{\rho}],
$$
which shows that condition (a) is satisfied in this case too.
\end{enumerate}
\medskip

Finally, since steps (4) and (5) preserve conditions (a) and (b), we claim that step (3) cannot occur because every extremal class of fibering type in the relative  Mori cone $\NE(X_\Sigma/X_{\bar{\Sigma}_{[D_1]}})$ have zero intersection with $(K_{X_\Sigma} + B)$. 
Indeed let $[C_0]$ be such a class. We have $J_{C_0}=J_{C_0}^+\subset \ker \bar{\phi}$ so that by condition (b), $I \cap J_{C_0} = \emptyset$. 
In particular we have $[C_0] \cdot \sum_{\rho \in I} [D_{\rho}] = 0$, that is
$[C_0] \cdot [K_{X_\Sigma} + B] = 0$ as claimed.
\medskip

It follows that when the algorithm finishes, we get a variety  $X_{\Sigma_2}$, a divisor $D_2$, a morphism $\phi_2 : X_{\Sigma_2} \rightarrow X_{\bar{\Sigma}_{[D_1]}}$, a 1-cycle $C_2$ a set $I_2 \subset J_{C_2}^+$ and a boundary divisor $B_2$ such that conditions (a) and (b) are satisfied by $(D_2,C_2,I_2)$ 
and the divisor $-\sum_{\rho \in I_2} D_{\rho}=K_{X_2}+B_2$ is $\phi_2$-nef.
By construction, the composition $\xi$ of all the birational maps appearing at each step of the MMP is a birational toric contraction over $X_{\bar{\Sigma}_{[D_1]}}$ and we have $D_1=\xi^*D_2$. 
\medskip

It remains to prove that $D_2=\xi_*D_1$ has restricted low toric degree.
\smallskip

Since $J_{C_2}^-$ is contained in a single cone of $\Sigma_2$ by condition (b), we know that the class $[C_2]$ belongs to the Mori cone $\NE(X_{\Sigma_2})$. Let us consider its decomposition as a linear combination of extremal classes with positive rational coefficients~:
$$
[C_2]= \sum_{i=1}^{l} \lambda_i [C'_{i}] \quad
\lambda_1, \ldots, \lambda_l \in \QQ_{+}^{*}.
$$ 
Up to a change of numbering we may assume that there exists $k \in \lbrace 1, \ldots, l \rbrace$  such that \mbox{$C'_i \cdot D_2 > 0$} for $i=1, \ldots, k$ and $C'_i \cdot D_2 = 0$ for $i=k+1,\ldots, l$, that is, the extremal classes $[C'_{k+1}], \ldots, [C'_{l}]$ are contracted by $\phi_2$. 
Since $-\sum_{\rho \in I_2} D_{\rho}$ is $\phi_2$-nef, we have
$$
C'_i \cdot \sum_{\rho \in I_2} D_{\rho} \leq 0 \quad \text{for} \; i=k+1,\ldots, l.
$$ 
By condition (a) for $(C_2,I_2)$ this implies
$$
\sum_{i=1}^k  \lambda_i \ C'_i \cdot \sum_{\rho \in I_2}  D_{\rho} \geq
\sum_{i=1}^l \lambda_i \ C'_i \cdot \sum_{\rho \in I_2}  D_{\rho} = 
C_2 \cdot \sum_{\rho \in I_2} D_{\rho} >
C_2 \cdot D_2 = \sum_{i=1}^{k} \lambda_i \ C'_{i} \cdot D_2.  
$$ 
Since all the $\lambda_i$ are positive, it follows that there exists necessarily at least one index \mbox{$i \in \lbrace 1, \ldots k \rbrace$} such that 
$$
C'_{i} \cdot \sum_{\rho \in I_2}  D_{\rho} > C'_{i} \cdot D_2.
$$
Moreover, by definition of $J_{C'_{i}}^+$ we have
$$
 C'_{i} \cdot\sum_{\rho \in J_{C'_{i}}^+}  D'_{\rho} \geq
  C'_{i} \cdot\sum_{\rho \in  I' \cap J_{C'_{i}}^+ } D'_{\rho} \geq
 C'_{i} \cdot \sum_{\rho \in  I'}  D'_{\rho}.
$$
Finally, since $i \leq k$ we also have $C'_{i} \cdot D_2>0$ so that 
$$ 
0 < C'_{i} \cdot D_2 < \sum_{\rho \in  I'}  C'_{i} \cdot D'_{\rho}
$$
and since $C'_{i}$ is extremal, this implies that $D_2$ has restricted low toric degree by Lemma \ref{extrgral}. 
This achieves the proof of the proposition.
\end{preuve}

\begin{cor}\label{gltd=>dltd}
Let $X_{\Sigma}$ be a simplicial projective toric variety and $D$ be a nef Cartier divisor on $X_{\Sigma}$. If $D$ has global low toric degree then there exist a toric birational contraction \mbox{$\xi : X_\Sigma \dashrightarrow X_{\Sigma'}$} over $X_{\bar{\Sigma}_{[D]}}$ such that $\xi_*[D]$ has restricted low toric degree.
\end{cor}

\begin{preuve}
It is sufficient to show that whenever $D$ is a nef Cartier divisor having global low toric degree on $X_{\Sigma}$, there exists a 1-cycle $C$ and a subset $I \subseteq J_{C}^{+}$ satisfying conditions (a) and (b) of Proposition \ref{descent}.

It is clear that condition (a) is just condition (A) of Definition \ref{gltd}, minus the condition on the cardinality of the indices set $I$. 
Let us suppose that condition (B) of the definition is verified and prove that so is condition (b) of the proposition.

For $\rho \in I$, let us denote by $\gamma_\rho=\barphi{J_{C} \setminus \{\rho\}}$ the smallest cone of $\Sigma_{[D]}$ containing $J_{C} \setminus \{\rho\}$.
It does not contain $\rho$ because otherwise we would have $J_C \subseteq \gamma_\rho$ and hence $C \cdot D = 0$ by Corollary \ref{CD0}.
In particular $\rho$ is not contained in the cone $\gamma:=\bigcap_{\rho' \in I} \gamma_{\rho'} \in \Sigma_{[D]}$. 
Now it suffices to remark that $J_C^{-}$ is contained in $J_C \setminus \{\rho\}$ for all $\rho \in I$ and hence $J_C^{-} \subseteq \gamma$.
This shows that (b) is verified and the corollary is proved.   
\end{preuve}

\begin{demo}{Proof of Theorem \ref{gral}}
We prove the implications  circularly~: 
$(i) \Rightarrow (ii) \Rightarrow (iii) \Rightarrow (i)$.

\begin{enumerate}
\item[$(i) \Rightarrow (ii)$] Suppose that the ample divisor $\bar{D}$ has restricted low toric degree. 
By Lemma \ref{rgltd} there exists a complete simplicial toric subvariety $\bar{W} \subseteq X_{\bar{\Sigma}_{[D]}}$ of Picard number one, such that the restriction $\bar{D}_{|\bar{W}}$ has global low toric degree.
Then by Lemma \ref{Wlift} there exist a toric subvariety $V \subseteq X_{\Sigma}$ such that the restriction $\phi_{|V} : V \rightarrow \bar{W}$ is a birational toric morphism.
By Proposition \ref{open2} the fact that $\bar{D}_{|\bar{W}}$ has global low toric degree implies that so has $D_{|V}$, which is $(ii)$.

\item[$(ii) \Rightarrow (iii)$] Suppose that $(ii)$ is verified. 
By \cite[Th. 11.1.9]{CLS} there exists a smooth refinement $\tilde{\Sigma}$ of $\Sigma$ that can be obtained by a sequence of star subdivisions. 
The induced morphism \mbox{$\eta : X_{\tilde{\Sigma}} \rightarrow X_{\Sigma}$} is a sequence of divisorial extremal contractions, hence a contracting toric morphism.
In particular, by Lemma \ref{Wlift} there exists a subvariety $\tilde{V} \subseteq X_{\tilde{\Sigma}}$ such that the restriction 
$\eta_{\tilde{V}} : \tilde{V} \rightarrow V$ is a birational toric morphism. 
Then it follows from Proposition \ref{open2} that the restriction $\tilde{D}_{|\tilde{V}}$ of the divisor $\tilde{D}=\eta^* D$ to $\tilde{V}$ has global low toric degree.
Finally since $X_{\tilde{\Sigma}}$ is smooth this implies that $\tilde{D}$ itself has global low toric degree, which is $(iii)$.

\item[$(iii) \Rightarrow (i)$] Suppose that $(iii)$ is verified. 
First notice that since $[\tilde{D}]=\eta^*[D] = \eta^*\phi^*[\bar{D}]$, the morphism $\tilde{\phi} : X_{\tilde{\Sigma}} \rightarrow X_{\bar{\Sigma}_{[\tilde{D}]}}$ of Theorem \ref{phi}, contracting all the curves having zero intersection with $\tilde{D}$ is just $\phi \circ \eta$. In particular we have $X_{\bar{\Sigma}_{[D]}}=X_{\bar{\Sigma}_{[\tilde{D}]}}$.

Now by Corollary \ref{gltd=>dltd}, there exist a toric a birational contraction $\xi : X_{\tilde{\Sigma}} \dashrightarrow X_{\Sigma'}$ such that $D':=\xi_*\tilde{D}$ has restricted low toric degree.

It follows from Lemma \ref{SigmabarD'} (and its proof) that $X_{\bar{\Sigma}_{[D']}}=X_{\bar{\Sigma}_{[\tilde{D}]}}=X_{\bar{\Sigma}_{[D]}}$ and $[D']$ is the pullback of $[\bar{D}]$ by the morphism $\phi': X_{\Sigma'} \rightarrow X_{\bar{\Sigma}_{[D]}}$ contracting all the curves that have zero intersection with $D'$. In particular, by Proposition \ref{dltdbar}, the fact that $D'$ has restricted low toric degree implies that so has $\bar{D}$, which is $(i)$. 
\end{enumerate}

This concludes the proof of Theorem \ref{gral}.
\end{demo}

\section{More criteria}\label{Geom}

In this section we give the weakest and the most general criteria we could find to determine whether a given nef Cartier divisor has low toric degree. Then we give some more natural criteria, linking the notion of low toric degree to different notions of positivity for divisor classes. 
\bigskip

\subsection{Most general criterion}

Let us first remark that along the proof we found a numerical condition weaker to the global low toric degree condition that implies the low degree one~:

\begin{lem}\label{easiest}   
$X_{\Sigma}$ be a simplicial projective toric variety and $D$ be a nef Cartier divisor on $X_{\Sigma}$. 
Let $\Sigma_{[D]}$ be the (possibly  degenerate) fan obtained from $\Sigma$ by removing the codimension 1 cones $\tau \in \Sigma(n-1)$ such that $C_\tau \cdot D=0$.

Then $D$ has low toric degree has soon as there exist a 1-cycle $C$ on $X_\Sigma$ and a subset $I\subseteq J_{C}^+$ verifying
\begin{equation*}\label{titdeggen4} \tag{a}
0< C \cdot D <  \sum_{\rho \in I} C \cdot D_{\rho},
\end{equation*}
and
\begin{equation*}\label{adhoc4} \tag{b}
\exists \sigma \in \Sigma_{[D]}, \  J_{C}^- \subseteq \sigma \quad 
\text{and} \quad
\forall \rho \in I, \ \rho \not\subseteq \sigma.
\end{equation*}

\end{lem}

\begin{preuve}
It is sufficient to remark that conditions (a) and (b) allow to apply Proposition \ref{descent}. Then it follows from Proposition \ref{dltdbar} that $\bar{D}$ has restricted low toric degree (see the proof of the implication $(iii) \Rightarrow (i)$ of Theorem \ref{gral}).
\end{preuve}

It then follows from Theorem \ref{gral} and Definition \ref{ltd} that the low toric degree property of the nef Cartier divisor $D$ on $X_\Sigma$ is essentially a property of its ``ample model" $\bar{D}$ on $X_{\bar{\Sigma}_{[D]}}$. In particular we have the following consequence

\begin{theo}\label{mostgral}
$X_{\Sigma}$ be a simplicial projective toric variety and $D$ be a nef Cartier divisor on $X_{\Sigma}$. 
Let $\phi : X_{\Sigma} \rightarrow X_{\bar{\Sigma}_{[D]}}$ be the toric morphism of Theorem \ref{phi}, contracting all the rational curves having zero intersection with $D$
and let $\bar{D}$ be an ample divisor on $X_{\bar{\Sigma}_{[D]}}$ such that $D$ is linearly equivalent to the pullback $\phi^*\bar{D}$.

The divisor $D$ has low toric degree if and only if one of the following equivalent conditions is satisfied~:
\begin{enumerate}
\item[(1)] There exists a toric subvariety $V \subseteq X_{\Sigma}$ such that the restriction $D_{|V}$  has low toric degree.
\item[(2)] There exists a simplicial toric variety $X_{\Sigma'}$ and a toric rational map $\xi : X_\Sigma \dashrightarrow X_{\Sigma'}$ over $X_{\bar{\Sigma}_{[D]}}$ such that the pullback of $\bar{D}$ in $X_{\Sigma'}$  has low toric degree.
\item[(3)] There exist a simplicial toric variety $X_{\Sigma'}$ over $X_{\bar{\Sigma}_{[D]}}$, a subvariety $V' \subseteq X_{\Sigma'}$, a 1-cycle $C'$ on $V'$ and a subset $I'\subseteq J_{C'}^+$ verifying conditions (a) and (b) of Lemma \ref{easiest}.
\end{enumerate}
\end{theo}

\begin{preuve}
The three conditions are clearly necessary so we just have to see that they are sufficient. Let us show that they imply the first condition of Theorem \ref{gral}~: 
\begin{enumerate}
\item[(1)] If we denote by $\phi_{V} : V \rightarrow \bar{V}$ the morphism contracting all the curves that have zero intersection with $D_{|V}$ then we have $\phi_{V}=\phi_{|V}$ and $D_{|V}=\phi_{V}^*\bar{D}_{|\bar{V}}$. It follows that $D_{|V}$ has low toric degree if and only if 
$\bar{D}_{|\bar{V}}$ has restricted low toric degree which implies that $\bar{D}$ has restricted low toric degree.
\item[(2)] By hypothesis there exists a proper toric morphism $\phi' : X_{\Sigma'} \rightarrow  X_{\bar{\Sigma}_{[D]}}$. Let us put $D'=\phi'^*\bar{D}$.
It follows from Lemma \ref{SigmabarD'} that $X_{\bar{\Sigma}_{[D']}}=X_{\bar{\Sigma}_{[D]}}$ and hence $D'$ has low toric degree if and only if 
$\bar{D}$ has restricted low toric degree. 
\item[(3)] This one is just a combination of the other two and Lemma \ref{easiest}.
\end{enumerate}
The theorem is proved.
\end{preuve}

\subsection{Criteria in terms of positivity cones}

The  inequalities characterizing a low toric degree, of the form
$$
0 < C \cdot D < C \cdot \sum_{\rho \in I} D_{\rho}, 
$$  
can be interpreted in terms of the proximity of $[D]$ to the border of a cone of ``positive divisor classes''. Indeed if the hyperplane $H_{[C]} \subset N^1(X_\Sigma)$ of classes having zero intersection with $[C]$ is a supporting hyperplane of a given polyhedral cone $\KKK$ of divisor classes, then the first inequality holds as soon as $[D]$ belongs to the interior of $\KKK$ and the second exactly means that substracting a divisor of the form $[\sum_{\rho \in I} D_{\rho}]$ is enough to get out of $\KKK$. We detail two important examples when $\KKK=\Nef(X_{\Sigma})$ and $\KKK=\Eff(X_{\Sigma})$.

\begin{prop}\label{top1}
Let $X_{\Sigma}$ be a simplicial projective toric variety and $D$ be an ample divisor on $X_{\Sigma}$. 
If there exists a subset $I \subseteq \Sigma(1)$ such that 
\begin{equation}\label{horsnef}
\left[ D - \sum_{\rho \in I} D_{\rho} \right] \not\in \Nef(X_{\Sigma})
\end{equation}
then $D$ has restricted low toric degree.
\end{prop}

\begin{preuve}
It is a direct consequence of the fact that the Mori cone $\NE(X_\Sigma)$ is the dual of the nef cone $\Nef(X_\Sigma)$ and is generated by extremal classes.
Indeed the fact that $D - \sum_{\rho \in I} D_{\rho}$ is not nef means that there exists an extremal curve $C$ such that $C \cdot (D - \sum_{\rho \in I} D_{\rho}) <0$.
By Lemma \ref{extrgral} and Proposition \ref{tresbeau} this is sufficient to prove that $D$ has restricted low toric degree. 
\end{preuve}

\begin{Remarks} 
\begin{enumerate}
\item When $X_\Sigma$ is smooth, then condition \eqref{horsnef} is in fact necessary and sufficient. In the general case, as explained in paragraph \ref{InterSubvar}, we have to consider every toric subvariety.   
\item The same arguments work with the sum $\sum_{\rho \in  I} D_{\rho}$ replaced by any $\QQ$-divisor $D'$ satisfying $0 \leq D' \leq -K_{X_\Sigma}$. 
Moreover, it is easy to see that in any generalized weighted projective space $W$ of dimension $k$ different from $\PP^k$, there exists a toric curve $C$ such that \mbox{$C \cdot (-K_W)< k+1$}. 
Using inequality \eqref{locglob3}, this shows that if $D$ has restricted low toric degree then there exists an extremal curve $C_\tau$ in $X_\Sigma$ such that $C_\tau \cdot D < \dim(W)+1$ where $W$ is a general fiber of the contraction of $[C_\tau]$. In particular $C_\tau \cdot D \leq n$ and we recover the simplicial case of Fujita's freeness conjecture for toric varieties (see \cite{Payne2} and references therein). 
\item As Theorem \ref{Kollar} applies to any intersection of hypersurfaces, Proposition \ref{top1} extends immediately to any intersection of ample divisors. Such an intersection contains a $K$-rational point as soon as the sum of the divisors has restricted low toric degree.  
\end{enumerate}
\end{Remarks}

\begin{prop}\label{big}
Let $X_{\Sigma}$ be a simplicial projective toric variety and $D$ be a Cartier divisor on $X_{\Sigma}$ that is nef and big. 
If there exists a subset $I \subseteq \Sigma(1)$ such that 
$$
\left[ D - \sum_{\rho \in I} D_{\rho} \right] \not\in \Eff(X_{\Sigma})
$$
then $D$ has low toric degree.
\end{prop}

\begin{preuve}
Since the dual cone of the effective cone is the cone of classes of mobile curves, by reasoning like in the proof of the previous proposition we prove that there exists a mobile curve $C$ satisfying
$$
0 < C \cdot D < C \cdot \sum_{\rho \in I} D_{\rho}. 
$$ 
This is condition (a) of Lemma \ref{easiest}  (1), it remains to prove condition (b). 
Since $D$ is big, the fan $\Sigma_{[D]}$ is not degenerate, that is, its minimal cone is the trivial cone $\{ 0 \}$. Hence we can take $\sigma=\{ 0 \}$ and since $C$ is mobile, we have $J_C^-=\emptyset \subset \sigma$ so that (b) is satisfied and the proposition is proved.
\end{preuve}

\section{Relation to Rational Connectedness}\label{RC}

\subsection{Separably rationally connected hypersurfaces have low toric degree}

Let us recall some definitions.

A variety $Y$ of dimension $n$ over an algebraically closed field $\bar{K}$ is said to be uniruled (resp. separably uniruled) if there exists a projective variety $Z$ of dimension $n-1$ and a dominant (resp. dominant and separable) rational map  $\varphi : \PP^1 \times Z \dashrightarrow Y$.
Similarly a variety $Y$ over $\bar{K}$ is said to be rationally connected (resp. separably rationally connected) if there exists a quasi-projective variety $Z$ and a morphism $\varphi : \PP^1 \times Z \rightarrow Y$ such that the morphism 
$$
\PP^1 \times \PP^1 \times Z \rightarrow Y \times Y, \quad \ (t,t',z) 
\mapsto \left( \varphi(t,z),\varphi(t',z)\right)
$$
is dominant (resp. a dominant and separable).

A variety $Y$ over an arbitrary field $K$ is said to be uniruled, separably uniruled,
rationally connected or separably rationally connected, if the variety $Y \otimes_K \bar{K}$ is, where $\bar{K}$ in an uncountable algebraically closed extension of $K$. 

It is clear from the definitions that any rationally connected (resp. separably rationally connected) variety is also uniruled (resp. separably uniruled). It is also clear that the four notions are stable by any birational modification that is a separable map, i.e. generically smooth, in particular by any divisorial contraction.
Similarly it follows from the definition that (separable) rational connectedness is stable by any surjective morphism.  
\medskip

A variety $Y$ is separably uniruled (resp. separably rationally connected) if and only if it contains a free curve (resp. a very free curve), that is roughly speaking, a curve that can be deformed to pass by any point of $Y$ (resp. by any two points of $Y$). For precise definitions and proofs see \cite[Cours 2]{RatConBonav}. 

For our arithmetic purpose, the key property of free curves on a variety $Y$ is that they have negative intersection with the canonical divisor. When $Y$ is a smooth hypersurface of a smooth ambient $X$, this gives a low degree to $Y$ through the adjunction formula. 

Before stating the result let us recall the general hypotheses defining our framework (see \ref{TorC1} and \ref{Homcoor})~: we consider a simplicial projective toric variety $X_{\Sigma}$ defined over a $C_1$ field $K$ such that  

\begin{itemize}
\item[(H1)] $K$ admits normic forms of arbitrary degree.
\item[(H2)] The torus $T$ acting on $X_{\Sigma}$ is split over $K$.
\item[(H3)] The order of the torsion part of the class group $\Cl(X_{\Sigma})$ is prime to the characteristic of $K$.
\end{itemize} 

Then we have the following result~: 

Here is the main result of the paper~:

\begin{theo}\label{RCpoint}
Let $X_{\Sigma}$ be a simplicial projective toric variety defined over a $C_1$ field $K$ such that hypotheses (H1), (H2) and (H3) are verified (see \ref{TorC1} and \ref{Homcoor}). Let $D$ be a hypersurface of $X_{\Sigma}$ defined over $K$.  
If $D$ is smooth and separably rationally connected then it has low toric degree and in particular admits a rational point over $K$. 
\end{theo}

The proof of the theorem is split in two parts stated as propositions.

\begin{prop}\label{bigunir}
Let $K$ and $X_{\Sigma}$ be as in Theorem \ref{RCpoint}.  
Let $D$ be a big divisor of $X_{\Sigma}$ defined over $K$.  
If $D$ is separably uniruled and smooth, then it has low toric degree.
\end{prop} 

\begin{preuve}
First recall that we may assume that $D$ is Cartier and nef because otherwise it has low degree by Definition \ref{ltd}.

Then let us consider a desingularization $\psi : X_{\tilde{\Sigma}} \rightarrow X_{\Sigma}$ and the pullback $\tilde{D}=\psi^*D$ of $D$ by $\psi$. We may assume $\tilde{D}$ is also the strict transform of $D$ because otherwise this would mean that $D$ contains some component of the blown-up locus, which consist only of toric subvarieties and hence $D$ would again contain torus-invariant points.

Since $\psi$ is birational it follows that $\tilde{D}$ is also separably uniruled and  
smooth. Moreover it has low toric degree if and only if $D$ has by Theorem \ref{mostgral} (2), so that up to replacing $X_{\Sigma}$  by $X_{\tilde{\Sigma}}$ and $D$ by $\tilde{D}$, we may assume that $X_{\Sigma}$ is smooth.

Now since $D$ is separably uniruled, it contains a free curve $C$.

The key point is that the free curve $C$ have negative intersection with the canonical divisor of $D$~:

$$
C \cdot K_D <0.
$$
  
Let us consider the pullback of divisors $\iota^*$ and the pushforward of curves $\iota_*$ associated to the inclusion $\iota : D \hookrightarrow X_{\Sigma}$. 

Since $X_\Sigma$ and $D$ are smooth, we have by the adjunction formula 
$$
 K_D = \bigl(K_{X_{\Sigma}}+D \bigr)_{|D} =  \iota^*\bigl(K_{X_{\Sigma}} + D \bigr).
$$ 

On the other hand, by the projection formula we have~:

$$
 \iota_*C \cdot \bigl(K_{X_{\Sigma}}+D \bigr)
=  C \cdot \iota^*\bigl( K_{X_{\Sigma}}+D \bigr) 
=  C \cdot K_D  < 0.
$$

Using the fact that $K_{X_{\Sigma}}=-\sum_{\rho \in \Sigma(1)} D_{\rho}$ this gives
$
\iota_*C \cdot \biggl( D - \sum_{\rho \in \Sigma(1)} D_{\rho} \biggr) <0.      
$

Finally, since $C$ is a mobile curve on $D$, which is nef, $\iota_*C$ is also mobile on $X_\Sigma$. Since $D$ is a big this gives $\iota_*C \cdot D>0$ and we have  

$$
0 < \iota_*C \cdot D < \iota_*C \cdot \sum_{\rho \in \Sigma(1)} D_{\rho}
= \iota_*C \cdot  \sum_{\rho \in J_C} D_{\rho}, 
$$

which is condition (a) of Lemma \ref{easiest} with $I=J_C=J_C^+$. It remains to prove condition (b), that is, there exists a cone $\sigma$ in the fan $\Sigma_{[D]}$ containing all the set $J_C^-$ but no element of $J_C^+$. 
Since $C$ is mobile we have $J_C^-=\emptyset$, hence we can choose $\sigma=\{ 0 \}$ because $\Sigma_{[D]}$ since $D$ is big. As $\{ 0 \}$ contains no 1-cone of $\Sigma$ condition (b) is satisfied and $D$ has low degree. The lemma is proved.
\end{preuve}
\bigskip

The big hypothesis is necessary here because if the divisor $D$ is not big it may happen that for all free curve $C$ in $D$ we have $\iota_*(C) \cdot D =0$. 
For example if $E$ is a non singular elliptic curve in $\PP^2$ without rational point over $K$, then the nef but not big hypersurface
$$
Y=E \times \PP^1 \subset \PP^2 \times \PP^1
$$ 
has no rational point over the $C_1$ field $K$, though it is smooth and separably uniruled in the smooth complete toric variety $\PP^2 \times \PP^1$.
\bigskip

\begin{prop}\label{nonbigSRC}
Let $K$, $X_{\Sigma}$ be as in Theorem \ref{RCpoint}.  
Let $D$ be a non big divisor of $X_{\Sigma}$ defined over $K$.  
If $D$ is separably rationally connected and smooth, then it has low toric degree.
\end{prop}

\begin{preuve}
We use the fact that separable rational-connectedness is stable both by surjective morphisms and by birational maps.

First, since the morphism $\phi : X_\Sigma \to X_{\bar{\Sigma}_{[D]}}$ of Theorem \ref{phi} is surjective, the ample divisor $\bar{D}$ is also separably rationally connected.
Then by Lemma \ref{Wlift} there exist a toric subvariety $V \subseteq X_{\Sigma}$ such that the restriction $\phi_{|V} : V \rightarrow X_{\bar{\Sigma}_{[D]}}$ is a birational toric morphism.
This induces a birational morphism between the divisors $D_{|V}$ and $\bar{D}$ and it follows that  $D_{|V}$ is itself separably rationally connected and in particular separably uniruled.

Finally, the morphism $\phi_{|V} : V \rightarrow X_{\bar{\Sigma}_{[D]}}$ is the morphism contracting all the curves that have zero intersection with $D_{|V}$, so the fact that it is birational means that $D_{|V}$ is a big divisor in $V$. 
We may thus apply Proposition \ref{bigunir} to prove that $D_{|V}$ has low toric degree. This implies that $D$ has low toric degree by Theorem \ref{mostgral} which  
proves the proposition and achieves the proof of Theorem \ref{RCpoint}.
\end{preuve}

\subsection{A low degree condition implying Rational Chain-Connectedness}\label{Ltd=>RCC}

The question of the converse implication, from low toric degree to rational connectedness seems difficult to treat in full generality by our method but as suggested by C. Araujo, it is quite easy to find conditions ``of low toric degree type" implying a weak version of rational connectedness~: rational chain-connectedness.  

While any two points of a rationally connected variety $Y$ over a field $K$ can be joined by a rational curve over the algebraically closed extensions of $K$, two points of a rationally chain-connected variety can be joined only by a chain of such rational curves. 
The condition is weaker and implies less good properties. For instance rational chain-connectedness is stable by fibration but not by birational modification. 
Indeed for any variety $Y$, the cone over $Y$ is rationally chain-connected because any two points can be joined to the node by a line and hence joined together by a length 2 chain of rational curves. 
But if $Y$ is not rationally chain-connected then the blow-up of the cone at the node is not either rationally chain-connected. 

The following result gives a condition that imply rational chain-connectedness.
It is easy to see that it implies also (global \textit{and} restricted) low toric degree. 
Notice that in this case no requirement is done about the singularities of the hypersurface.

\begin{prop}\label{RCC}
Let $K$ and $X_{\Sigma}$ be as in Theorem \ref{RCpoint}.  
Let $D$ be a hypersurface of $X_{\Sigma}$ defined over $K$.  
If there exists a mobile extremal curve $C$ such that 
$$
0 < C \cdot D < C \cdot (-K_{X_{\Sigma}})
$$
then $D$ is rationally chain-connected.  
\end{prop} 

\begin{preuve}
Since the class $[C]$ is mobile and extremal its contraction is a Mori fibration
$$
\pi : X_{\Sigma} \rightarrow X_{\bar{\Sigma}}
$$ 
with general fiber $F$ isomorphic to a generalized weighted projective space. 
The condition $C \cdot D < C \cdot (-K_{X_{\Sigma}})$ means that the restriction $D_{|F}$ of $D$ to such a fiber is $\QQ$-Fano and hence rationally chain connected. 

On the other hand, the condition $0 < C \cdot D$ implies that $D$ is linearly equivalent to a torus invariant effective divisor $\sum a_\rho D_\rho$ with at least one $\rho \in \Sigma(1)$ such that $a_\rho > 0$ and $C \cdot D_\rho > 0$. Since such a $D_\rho $ is a section of $\phi$ it follows that $D$ dominates $X_{\bar{\Sigma}}$.

Putting this together we have that the restriction of $\phi$ to $D$ is a fibration onto the toric (hence rational) variety $X_{\bar{\Sigma}}$ with rationally chain connected fibers $D_{|F}$. Hence $D$ is rationally chain connected and the proposition is proved. 
\end{preuve}
 
If the hypersurface is smooth and the characteristic of $K$ is zero, then we can make two improvements. First, rational connectedness and rational chain-connectedness are equivalent for smooth projective varieties (see \cite[Cours 4]{RatConBonav}). 
Second, in characteristic zero every map is separable so that rational connectedness is also equivalent to separable rational connectedness. 
It follows that in this case we have
\begin{center}
Low degree along a mobile extremal curve $\Rightarrow$ Rationally connected $\Rightarrow$ Low toric degree. 
\end{center}
\smallskip

Let us remark that having low degree along a mobile extremal curve is quite a strong condition for $D$ since it implies that $D$ is a fiber bundle with toric base whose fibers are $\QQ$-Fano hypersurfaces of a generalized weighted projective space. 
\bigskip

Here is an example.

\begin{ex}\label{belex}
Let $Y$ be a hypersurface of $\PP^n$ of degree $d$ containing a projective subspace $Z$ of dimension $k$ with multiplicity $m$. Let $\eta : X_{\Sigma} \longrightarrow \PP^n$ denote the blowing-up of $\PP^n$ along $Z$ and let $D$ be the strict transform of $Y$ in $X_\Sigma$. Notice that $D$ may be singular.

Let us number the toric divisors from 1 to $n+2$ with $D_{n+2}$ the exceptional divisor of the blowing-up and denote by $D'_i=\eta(D_i)$ the image of $D_i$ in $\PP^n$. 
Up to projectivity we may assume that $Y=D_1 \cdot D_2 \cdots D_{n-k}$ is a torus-invariant subvariety of $\PP^n$. 
There are three different classes of toric curves in $X_\Sigma$~: the class $[C_1]$ of a line contracted by $\eta$, the class $[C_2]$  of the strict transform of a line in $\PP^n$ containing a point of $Z$ and the class $[C_3]=[C_1]+[C_2]$  of the strict transform of a line disjoint from $Z$.

The corresponding intersection table is given by

$$
 \begin {array}{c|ccccccc}  
\pmb{\cdot} & D_1 & \cdots & D_{n-k} & D_{n-k+1} & \cdots & D_{n+1} & D_{n+2}   \\
\hline 
C_{1} \quad &  1  & \cdots &1 & 0& \cdots &0&-1  \\ 
C_{2} \quad &  0  & \cdots &0 & 1& \cdots &1&1 \\  
C_{3} \quad &  1  & \cdots &1 & 1& \cdots &1&0 \\ 
\end {array} 
$$
\vspace{-0.6cm}
$$
 \underset{n-k}{\underbrace{\hspace{2.1cm}}}
\hspace{1.2cm} \underset{k+1}{\underbrace{\hspace{2.3cm}}} \hspace{0.3cm}
$$
The classes $[C_1]$ and $[C_2]$ are extremal and $[C_2]$ is moreover a mobile class, whose contraction is the fibration of $X_\Sigma$ onto $\PP^{n-k-1}$.

The hypersurface $D$ dominates $\PP^{n-k-1}$. Indeed by construction it is a fiber bundle with base $\PP^{n-k-1}$ and fibers isomorphic to hypersurfaces of $\PP^{k+1}$ of degree $d-m$. In particular we have 
$$
C_1 \cdot D = m, \ C_2 \cdot D = d - m \ \, \text{and} \, \ C_3 \cdot D = d. 
$$

It follows that whenever the multiplicity $m$ is greater than $d-k-2$ we have
$$
0< C_2 \cdot D = d - m < k+2 = C_2 \cdot (-K_{X_{\Sigma}}), 
$$ 
and proposition \ref{RCC} applies~: $D$ has low toric degree and is rationally chain-connected for the same reason.
\end{ex}

\thispagestyle{empty}
\bibliographystyle{alpha}
\bibliography{biblio}

\end{document}